\documentclass{amsart}
\usepackage{amsrefs}
\usepackage{amsmath,amssymb,stmaryrd}
\usepackage{graphicx}

\def \rr {\mathbb{R}}
\def \nn {\mathbb{N}}
\def \rn {\mathbb{R}^n}
\def \rnp {\mathbb{R}^n_+}

\def \crits {2^\star(s)}
\def \crit {2^\star}
\def \ue {u_\eps}

\def \eps {\varepsilon}
\def \ap {\alpha_+(\gamma)}
\def \am {\alpha_-(\gamma)}
\def \aps {\alpha_+}
\def \ams {\alpha_-}
\def \dundeuxr {D^{1,2}(\rn)}
\def \dundeuxrnp {D^{1,2}(\rnp)}
\def \dundeux {D^{1,2}(\Omega)}
\def \eucl {\hbox{Eucl}}

\def \tOmega {\tilde{\Omega}}
\def \tPhi {\tilde{\Phi}}
\def \tD {\tilde{D}}
\def \tH {\tilde{H}}

\newtheorem{theorem}{Theorem}[section]
\newtheorem{proposition}[theorem]{Proposition}
\newtheorem{coro}[theorem]{Corollary}
\newtheorem{lem}[theorem]{Lemma}
\newtheorem{rem}[theorem]{Remark}
\newtheorem{definition}[theorem]{Definition}

\date{October 14th, 2014}

\title[Hardy-Schr\"odinger operator]{On the Hardy-Schr\"odinger operator with a boundary singularity}

\author{Nassif Ghoussoub}
\address{Nassif Ghoussoub, Department of Mathematics, 1984 Mathematics Road,The University of British Columbia, BC, Canada V6T 1Z2}
\email{nassif@math.ubc.ca}

\author{Fr\'ed\'eric Robert}
\address{Fr\'ed\'eric Robert, Institut \'Elie Cartan, Universit\'e de Lorraine, BP 70239, F-54506 Vand{\oe}uvre-l\`es-Nancy, France}
\email{frederic.robert@univ-lorraine.fr}


\thanks{This work was carried out while F. Robert was visiting the Pacific Institute for the Mathematical Sciences (PIMS) at the University of British Columbia in 2013-2014, as a member of the Unit\'e Mixte Internationale of the French Centre National de la Recherche Scientifique (CNRS). He thanks the CNRS (INSMI) and UBC (PIMS) for this support. N. Ghoussoub was partially supported by a research grant from the Natural Science and Engineering Research Council of Canada (NSERC)}
\begin{document}

\begin{abstract} We investigate the Hardy-Schr\"odinger operator $L_\gamma=-\Delta -\frac{\gamma}{|x|^2}$ on domains $\Omega\subset\rn$, whose boundary contain the singularity $0$. The situation is quite different from the well-studied case when $0$ is in the interior of $\Omega$. For one, if $0\in\Omega$, then  $L_\gamma$ is positive if and only if $\gamma<\frac{(n-2)^2}{4}$, while  if 
$0\in\partial\Omega$ the operator  $L_{\gamma}$ could be positive for larger value of $\gamma$,  potentially reaching the maximal constant $\frac{n^2}{4}$ on convex domains.  

We prove optimal regularity and a Hopf-type Lemma for variational solutions of corresponding linear Dirichlet boundary value problems of the form $L_{\gamma} u=a(x)u$, but also for non-linear equations including $L_{_\gamma} u=\frac{|u|^{\crits-2}u}{|x|^s}$, where $\gamma <\frac{n^2}{4}$, $s\in [0,2)$ and $\crits:=\frac{2(n-s)}{n-2}$ is the critical Hardy-Sobolev exponent. We also provide a Harnack inequality and a complete description of the profile of all positive solutions --variational or not-- of the corresponding linear equation on the punctured domain. The value $\gamma=\frac{n^2-1}{4}$ turned out to be another critical threshold for the operator $L_\gamma$, and our analysis yields a corresponding notion of  ``Hardy singular boundary-mass" $m_\gamma(\Omega)$ of a domain $\Omega$ having $0\in \partial \Omega$, which could be defined whenever $\frac{n^2-1}{4}<\gamma<\frac{n^2}{4}$. 

As a byproduct, we give a complete answer to problems of existence of extremals for Hardy-Sobolev inequalities of the form 
\[
\hbox{$C\left(\int_{\Omega}\frac{u^{\crits}}{|x|^s}dx\right)^{\frac{2}{\crits}}\leq \int_{\Omega} |\nabla u|^2dx-\gamma \int_{\Omega}\frac{u^2}{|x|^2}dx$\quad  for all $u\in \dundeux$,}
\]
whenever $\gamma <\frac{n^2}{4}$, and in particular, for those of Caffarelli-Kohn-Nirenberg. These results extend previous contributions by the authors in the case $\gamma =0$, and by Chern-Lin for the case $\gamma<\frac{(n-2)^2}{4}$. Namely, if $0\leq \gamma \leq \frac{n^2-1}{4}$, then the negativity of the mean curvature of $\partial \Omega$ at $0$ is sufficient for the existence of extremals. This is however not sufficient for $\frac{n^2-1}{4}<\gamma <\frac{n^2}{4}$, which then requires the positivity of the Hardy singular boundary-mass of the domain under consideration. 
\end{abstract}
\maketitle
\newpage
\tableofcontents

\section{Introduction}\label{sec:intro}

For any smooth open domain $\Omega \subset \rn$, define the best constant in the corresponding Hardy inequality by, 
\begin{equation} 
\gamma_H(\Omega):=\inf\left\{\frac{\int_{\Omega}|\nabla u|^2\, dx}{\int_{\Omega}\frac{u^2}{|x|^2}\, dx}; \,  \; u\in D^{1,2}(\Omega)\setminus\{0\}\right\}, 
\end{equation}
where $D^{1,2}(\Omega)$ is the completion of $C^{\infty}_c(\Omega)$ with respect to the norm given by $||u||^2=\int_{\Omega}|\nabla u|^2 dx.$ 
It is well known that $\gamma_H(\Omega)=\frac{(n-2)^2}{4}$ for any domain $\Omega$ having $0$ in its interior, including $\rn$,  and that it is never attained by a function in $D^{1,2}(\Omega)$.  On the other hand, it has been noted by several authors (See for example Pinchover-Tintarev \cite{PT} Fall-Musina \cite{FM} or the book of Ghoussoub-Moradifam \cite{gm}) that the situation is quite different for the half-space $\rnp$, in which  case, $\gamma_H(\rnp)=\frac{n^2}{4}$. More generally,  if $0\in \partial \Omega$ the boundary of $\Omega$, then $\gamma_H(\Omega)$ can be anywhere in the interval $\left(\frac{(n-2)^2}{4}, \frac{n^2}{4}\right]$ (see Proposition \ref{prop:gamma}). Moreover, $\gamma_H(\Omega)$ is attained whenever $\gamma_H(\Omega)<\frac{n^2}{4}$ (See Section \ref{sec:hardy}). This already points to the fact that the Hardy-Schr\"odinger operator $L_\gamma=-\Delta -\frac{\gamma}{|x|^2}$ behaves differently when the singularity $0$ is on the boundary of a domain $\Omega$, than when $0$ is in the interior. The latter case has already been extensively covered in the literature. Without being exhaustive, we refer to Ghoussoub-Yuan \cite{GY}, Guerch-V\'eron \cite{GV}, Jaber \cite{jaber},  Kang-Peng \cite{KP}, Pucci-Servadei \cite{PS}, Ruiz-Willem \cite{RW}, Smets \cite{SmetsTAMS}, and the references therein. 

The study of nonlinear singular variational problems when $0\in\partial\Omega$ was initiated by Ghoussoub-Kang \cite{gk} and was studied extensively by Ghoussoub-Robert \cites{gr1,gr2, gr3}. For more recent contributions, we refer  to Attar-Merch\'an-Peral \cite{AMP}, D\'avila-Peral \cite{DP}, 
and Gmira-V\'eron \cite{GmVe}. We also learned recently --after a first version of this paper was  posted on arxive-- about a paper of Pinchover \cite{PinchoverAIHP}, and a more recent one by Devyver-Fraas-Pinchover \cite{DFP} that also treat the Hardy potential when $0\in \partial \Omega$.  

Our main goal in this paper is to show that the above noted discrepancy --
between the case when the singularity $0$ belongs to the interior of the domain and when it is on the boundary-- is only the tip of the iceberg. The differences manifest themselves in both linear problems of the form 
\begin{eqnarray} \label{zero}
\left\{ \begin{array}{llll}
-\Delta u-\gamma \frac{u}{|x|^2}&=&a(x)u\ \ &\text{on } \Omega\\
\hfill u&=&0 &\text{on }\partial \Omega,
\end{array} \right.
\end{eqnarray} 
and in nonlinear Dirichlet boundary value problems associated to $L_\gamma$, such as: 
\begin{eqnarray} \label{one}
\left\{ \begin{array}{llll}
-\Delta u-\gamma \frac{u}{|x|^2}&=&\frac{u^{\crits-1}}{|x|^s} \ \ &\text{on } \Omega\\
\hfill u&>&0 &\text{on } \Omega\\
\hfill u&=&0 &\text{on }\partial \Omega,
\end{array} \right.
\end{eqnarray} 
where $0\leq s < 2$ and $\crits:=\frac{2(n-s)}{n-2}$. Actually, Ghoussoub-Kang noted in \cite{gk} that even when $\gamma=0$,  the situation can already be quite different whenever $0$ belongs to the boundary of a bounded $C^2$-smooth domain $\Omega$ as long as $s>0$.  Ghoussoub and Robert \cites{gr1, gr2} eventually proved that if the mean curvature at $0$ of such domains is negative, and provided $s>0$, then minimizers for the functional 
\begin{equation} 
J^\Omega_{s}(u):=\frac{\int_{\Omega} |\nabla u|^2dx}{(\int_{\Omega}\frac{u^{\crit}}{|x|^s}dx)^{\frac{2}{\crit}}}
\end{equation}
exist in $ \dundeux\setminus \{0\}=H^1_0(\Omega)\setminus \{0\}$ and are solutions to equation (\ref{one}) in the case when $\gamma=0$. While this new phenomenon occured because of the presence of the singularity $|x|^{-s}$ in the nonlinear term, we shall show in this paper, that the differences also appear on the linear level, as soon as $\gamma>0$, but also as one varies $\gamma$ between $0$ and $\frac{n^2}{4}$.

\medskip\noindent Another motivation for this work came from the recent work of C.S. Lin and his co-authors \cites{CL4, CL5} on the existence of extremals for the Caffarelli-Kohn-Nirenberg (CKN) inequalities \cite{ckn}. These  inequalities state that in dimension $n\geq 3$, there is a constant $C:=C(a,b,n)>0$ such that for all $u\in C^\infty_c(\rn)$, the following inequality holds:
\begin{equation} \label{CKN}
\left(\int_{\rn}|x|^{-bq}|u|^q \right)^{\frac{2}{q}}\leq C\int_{\rn}|x|^{-2a}|\nabla u|^2 dx,
\end{equation}
where
\begin{equation}\label{cond1}
-\infty<a<\frac{n-2}{2}, \ \ 0 \leq b-a\leq 1 \ {\rm and} \ q=\frac{2n}{n-2+2(b-a)}.
\end{equation}
A proof and various extensions of (\ref{CKN}) will be given in section \ref{sec:hardy}.

\medskip\noindent For a domain $\Omega$ in $\rn$, we let $D^{1,2}_a(\Omega)$ be the completion of $C^{\infty}_c(\Omega)$ with respect to the norm $||u||_{a}^2=\int_{\Omega}|x|^{-2a}|\nabla u|^2 dx.$ 
Consider the best constant defined as:
 \begin{equation}\label{best}
 S(a,b,\Omega)= \inf \left\{\frac{\int_{\Omega}|x|^{-2a}|\nabla u|^2dx}{\left(\int_{\Omega}|x|^{-bq}|u|^q \right)^{\frac{2}{q}}\, dx}; u\in D^{1,2}_a(\Omega) \backslash \{0\}\right\}.
 \end{equation}
The extremal functions for $S(a,b,\Omega)$ are then the least-energy solutions of the corresponding Euler-Lagrange equations: 
\begin{eqnarray}\label{CKN1}
\left\{ \begin{array}{lll}
-{\rm div}(|x|^{-2a}\nabla u)&=|x|^{-bq}u^{q-1} \,\, &\text{on } \Omega\\
\hfill u&>0& \text{on } \Omega\\
\hfill u&=0& \text{on }\partial \Omega.
\end{array} \right.
\end{eqnarray} 
The existence or nonexistence of minimizers for (\ref{best}), when the domain $\Omega$ is the whole space $\rn$, have been extensively
studied for the past twenty years, see Catrina-Wang \cite{CW}, Chou-Chu \cite{CC}, Dolbeault-Esteban-Loss-Tarantello \cite{DELT}, Lin-Wang \cite{LW1} and references
therein. The result can be briefly summarized in the following:

\medskip\noindent {\bf Theorem A:} {\it Suppose $n\geq 3$ and that $a, b$ and $q$ satisfy condition (\ref{cond1}). Then minimizers exist for the best constant $S(a; b; \rn)$ if and only if $a, b$ satisfy
\begin{equation} 
\hbox{either $a < b < a + 1$ or  $b = a \geq 0$.}
 \end{equation} }
 \noindent If now $\Omega$ is any domain in $\rn$ that contains $0$ in its interior, one can easily see that 
 scale invariance yields for any $\lambda >0$, that $S(a; b; \lambda\Omega) = S(a; b; \Omega)$
where $\lambda \Omega =\{\lambda x; x\in \Omega\}$. It follows that if $0\in \Omega$, then 
 $S(a; b; \Omega)= S(a; b; \rn)$, which means that 
$S(a; b;\Omega)$ can never be achieved unless $\Omega = \rn$ (up to a set of capacity zero).
However, 
 as mentioned above, if $0$ belongs to the boundary of a smooth bounded domain $\Omega$ and if the mean curvature at $0$ of such domains is negative, then minimizers for the best constant $S(0; b; \Omega)$ were shown to be attained \cites{gk, gr1, gr2}.
 
\noindent This result was later extended by Chern and Lin \cite{CL5}, who eventually established existence of minimizers under the same negative mean curvature condition at $0$ provided   $a, b$ satisfy one of the following conditions:
 \begin{eqnarray}\label{almost}
 \left\{ \begin{array}{ll}
 \hbox{(i)\,\, $a<b<a+1$ and $n\geq 3$}\\
 \hbox{(ii)\, $a=b$ and $n\geq 4$.}
 \end{array} \right.  
 \end{eqnarray}
 They left open the case when $n=3$ and $0<a=b<\frac{n-2}{2}$ open, a problem that we address in Theorem \ref{th:dim3} (see also Section \ref{sec:n:3}).

\medskip\noindent To make the connection, we note that by making the substitution 
$w(x)=|x|^{-a}u(x)$ for $x\in \Omega$, 
one can see that  if $a<\frac{n-2}{2}$, then $u \in D^{1,2}_{a}(\Omega)$ if and only if $w\in \dundeux$ by the Hardy inequality,
and
\begin{equation}\nonumber
\frac{\int_{\Omega}|x|^{-2a}|\nabla u|^2dx}{\left(\int_{\Omega}|x|^{-bq}|u|^q \right)^{\frac{2}{q}}}=\frac{\int_{\Omega} |\nabla w|^2-\gamma \int_{\Omega}\frac{w^2}{|x|^2}dx}{(\int_{\Omega}\frac{w^{\crit}}{|x|^s}dx)^{\frac{2}{\crit}}},
\end{equation}
where $\gamma= a(n-2-a)$, $s=(b-a)q$ and $\crit=\frac{2n}{n-2+2(b-a)}$. 
This means that 
$u$ is a solution of (\ref{CKN1}) if and only if $w(x)$ is a solution of equation (\ref{one})
where $0\leq s < 2$ and $\crit:=\crit(s)=\frac{2(n-s)}{n-2}$. Therefore, instead of looking for solutions of (\ref{CKN1}) one can study equation (\ref{one}). To state the result of Chern-Lin in this context, we define  the functional 
\begin{equation} 
J^\Omega_{\gamma, s}(u):=\frac{\int_{\Omega} |\nabla u|^2-\gamma \int_{\Omega}\frac{u^2}{|x|^2}dx}{(\int_{\Omega}\frac{u^{\crit}}{|x|^s}dx)^{\frac{2}{\crit}}},
\end{equation}
and its infimum on $\dundeux\setminus \{0\}$, that is
\begin{equation}\label{ineq:sobo:s:g}
\mu_{\gamma,s}(\Omega):=\inf\left\{J^\Omega_{ \gamma, s}(u); u\in \dundeux\setminus \{0\}\right\}.
\end{equation}

 \begin{theorem}[Chern-Lin \cite{CL5}] \label{CL} Let $\Omega$ be a smooth bounded domain in $\rn$ ($n\geq 3$). Assume $\gamma <\frac{(n-2)^2}{4}$ and $0\leq s < 2$. If either $\{s>0\}$ or $\{n\geq 4$ and $\gamma>0\}$, then there are extremals for $\mu_{\gamma,s}(\Omega)$,  provided  the mean curvature of $\partial \Omega$ at $0$ is negative.

 \end{theorem}
The case when $n=3$, $s=0$ and $\gamma>0$ remained open. Moreover,  as we shall  see  in section \ref{sec:hardy}, the infimum $\mu_{\gamma, s}(\Omega)$ is finite for all $\gamma <\frac{n^2}{4}$,  whenever $0\in \partial \Omega$. 
This means that equation (\ref{one}) may have positive solutions for $\gamma$ beyond $\frac{(n-2)^2}{4}$ and all the way to $\frac{n^2}{4}$. This turned out to be the case as we shall establish in this paper.

\medskip\noindent We first note that standard compactness arguments \cites{gk, gm, CL5} --also described in section \ref{compactness}--  yield that for $\mu_{\gamma,s}(\Omega)$ to be attained it is sufficient to have that 
 \begin{equation}\label{compact}
 \mu_{\gamma,s}(\Omega)<\mu_{\gamma,s}(\rnp),
 \end{equation}
 where the latter is the corresponding best constant on $\rnp$. In order to prove the existence of such a gap, one tries to construct test functions for $ \mu_{\gamma,s}(\Omega)$ that are based on the extremals of $\mu_{\gamma,s}(\rnp)$ provided the latter exist. The cases where this is known are given by the following standard proposition (see for instance Bartsch-Peng-Zhang \cite{BPZ} and Chern-Lin \cite{CL5}). See Proposition \ref{prop:ext:rnp} in the appendix for a proof. 
 
  \begin{proposition}\label{bartsh} Assume $\gamma<\frac{n^2}{4}$, $n\geq 3$ and $0\leq s<2$. Then, there are extremals for $\mu_{\gamma, s}(\rnp)$ provided either  $\{s>0\}$ or $\{n\geq 4$ and $\gamma>0\}$.\\
  On the other hand, 
 \begin{enumerate}
\item If $\{s=0$ and $\gamma \leq 0\}$, then there are no extremals for $\mu_{\gamma, 0}(\rnp)$ for any $n\geq 3$.
. \item Furthermore, whenever $\mu_{\gamma, 0}(\rnp)$ has no extremals, then necessarily
\begin{equation}\label{def:K}
\mu_{\gamma, 0}(\rnp)=\inf_{u\in D^{1,2}(\rn)\setminus\{0\}}\frac{\int_{\rn}|\nabla u|^2\, dx}{\left(\int_{\rn}|u|^{\crit}\, dx\right)^{\frac{2}{\crit}}}=\frac{1}{K(n,2)^2},
\end{equation}
where the latter is the best constant in the Sobolev inequality and $\crit:=\crit(0)=\frac{2n}{n-2}$. 
\end{enumerate}
\end{proposition}
 The only unknown situation is again when $s=0$, $n=3$ and $\gamma>0$, which we address below (see Theorem \ref{th:dim3}) and in full detail in Section \ref{sec:n:3}. For now, we shall discuss the new ingredients that we bring to the discussion. 

\medskip\noindent Assuming first that an extremal for $\mu_{\gamma, s}(\rnp)$ exists and that one knows its profile at infinity and at $0$, then this information can be used to construct test functions for $\mu_{\gamma, s}(\Omega)$. This classical method has been used by Kang-Ghoussoub \cite{gk}, Ghoussoub-Robert \cites{gr1, gr2} when $\gamma=0$, and by  Chern-Lin \cite{CL5} for $0<\gamma <\frac{(n-2)^2}{4}$ in order to establish (\ref{compact}) under the assumption that $\partial \Omega$ has a negative mean curvature at $0$. Actually, the estimates of Chern-Lin \cite{CL5} extend directly to establish an analogue of Theorem \ref{CL} for all $\gamma < \frac{n^2-1}{4}$ under the same local geometric condition.
However, the case where $\gamma =  \frac{n^2-1}{4}$ already requires a much more refined analysis of the Hardy-Schr\"odinger operator  $L_\gamma:=-\Delta -\frac{\gamma}{|x|^2}$ on
  the half-space $\rnp$.

 \medskip\noindent The remaining range $(\frac{n^2-1}{4}, \frac{n^2}{4})$ for $\gamma$ turned out to be even more interesting for the operator $L_\gamma$. Indeed, the curvature condition at $0$ is not sufficient anymore to insure existence, as more global test functions are required. We therefore proceed to isolate a notion of  ``Hardy boundary-mass"  $m_\gamma(\Omega)$ for any bounded domain $\Omega$ (with $0\in \partial \Omega$)  that is associated to the operator $L_\gamma$. This is stated in Theorem \ref{def:mass} below and is reminiscent of the positive mass theorem of Schoen-Yau \cite{SY} that was used to complete the solution of the Yamabe problem.

\medskip\noindent In order to explain the new critical threshold that is $\frac{n^2-1}{4}$, we need first to consider the Hardy-Schr\"odinger operator $L_\gamma:=-\Delta -\frac{\gamma}{|x|^2}$ on $\rnp$. The most basic solutions for $L_\gamma u=0$, with $u=0$ on $\partial \rnp$ are of the form $u(x)=x_1|x|^{-\alpha}$, and a straightforward computation yields
$-\Delta (x_1|x|^{-\alpha})=\frac{\alpha(n-\alpha)}{|x|^2}x_1|x|^{-\alpha}$ on $\rnp$, 
which means that
\begin{equation}\nonumber
\hbox{$\left(-\Delta-\frac{\gamma}{|x|^2}\right)\left(x_1|x|^{-\alpha}\right)=0$ on $\rnp$,}
\end{equation}
for $\alpha\in \{\am,\ap\}$ where $\alpha_{\pm}(\gamma):=\frac{n}{2}\pm\sqrt{\frac{n^2}{4}-\gamma}.$ Actually, any non-negative solution of $L_\gamma u=0$ on $\rnp$ with $u=0$ on $\partial \rnp$ is a (positive) linear combination of these two solutions (
Proposition \ref{prop:funda:2} below).

\medskip\noindent Note that $\am<\frac{n}{2}<\ap$, which points to  the difference between the ``small" solution, namely $x\mapsto x_1|x|^{-\am}$, and the ``large one"  $x\mapsto x_1|x|^{-\ap}$. Indeed, the ``small" solution is ``variational", i.e. is locally in $D^{1,2}(\rnp)$, while the large one is not. This turned out to be a general fact since we shall show that $x\mapsto d(x, \partial \Omega)|x|^{-\am}$ is essentially the profile at $0$ of any variational solution --positive or not-- of equations of the form $L_\gamma u=f(x, u)$ on a domain $\Omega$, as long as the nonlinearity $f$ is dominated by $C(|v|+\frac{|v|^{\crits-1}}{|x|^s})$. Here $d(x, \partial \Omega)$ denotes the distance function to $\partial \Omega$. To state the theorem, we use the following terminology. Say that $u\in \dundeux_{loc, 0}$ if there exists $\eta\in C^\infty_c(\rn)$ such that $\eta\equiv 1$ around $0$ and $\eta u\in \dundeux$. Say that $u\in \dundeux_{loc, 0}$ is a weak solution to the equation 
\[-\Delta u=F\in \left(\dundeux_{loc, 0}\right)^\prime,
\]
 if for any $\varphi\in \dundeux$ and $\eta\in C^\infty_c(\rn)$, we have 
$\int_\Omega (\nabla u,\nabla (\eta\varphi))\, dx=\left\langle F ,\eta\varphi\right\rangle.$

\medskip\noindent The following theorem will be established in section \ref{sec:hopf}.

\begin{theorem}[Optimal regularity and Generalized Hopf's Lemma]Let $\Omega$ be a smooth domain in $\rn$ such that $0\in \partial \Omega$, and let $f: \Omega\times \rr\to \rr$ be a Caratheodory function such that
$$|f(x, v)|\leq C|v| \left(1+\frac{|v|^{\crits-2}}{|x|^s}\right)\hbox{ for all }x\in \Omega\hbox{ and }v\in \rr.$$
Assume $\gamma<\frac{n^2}{4}$ and  let $u\in \dundeux_{loc, 0}$ be such that for some $\tau>0$, 
\begin{equation}\label{regul:eq.0}
-\Delta u-\frac{\gamma+O(|x|^\tau)}{|x|^2}u=f(x,u)\hbox{ weakly in }\dundeux_{loc, 0}.
\end{equation}
 Then, there exists $K\in\rr$ such that
\begin{equation}\label{eq:hopf.0}
\lim_{x\to 0}\frac{u(x)}{d(x,\partial\Omega)|x|^{-\am}}=K.
\end{equation}
Moreover, if $u\geq 0$ and $u\not\equiv 0$, then $K>0$.
\end{theorem}
This theorem can be seen as an extension of Hopf's Lemma \cite{gt} in the following sense: when $\gamma=0$ (and therefore $\am=0$), the classical Nash-Moser regularity scheme then yields that $u\in C^1_{loc}$, and when   $u\geq 0$, $u\not\equiv 0$, Hopf's comparison principle yields $\partial_\nu u(0)<0$, which is really a reformulation of \eqref{eq:hopf.0} in the case where  $\am=0$.

The proof of this theorem is quite interesting since, unlike the regular case (i.e., when $L_\gamma=L_0=-\Delta$) or the classical situation when the singularity $0$ is in the interior of the domain $\Omega$ (see Smets \cite{SmetsTAMS}), a direct application of the standard Nash-Moser iterative scheme is not sufficient to obtain the required regularity. Indeed, the scheme only yields the existence of $p_0$, with $1<p_0<\frac{n}{\am-1}$ such that  $u\in L^p$ for all $p<p_0$. Unfortunately, $p_0$ does not reach $\frac{n}{\am-1}$, which is the optimal rate of integration needed to obtain the profile \eqref{eq:hopf.0} for $u$. However, the improved order $p_0$ is enough to allow for the inclusion of the nonlinearity $f(x,u)$ in the linear term of \eqref{regul:eq.0}. We are then reduced to the analysis of the linear equation, that is \eqref{regul:eq.0} with $f(x,u)\equiv 0$,  
in which case we get the conclusion by constructing suitable super- and sub- solutions to the linear equation that have the same profile at $0$ as  \eqref{eq:hopf.0}.  For details, see Section \ref{sec:hopf}.

\noindent As a corollary, one obtains the following description of the profile of variational solutions  of (\ref{one}) on $\rnp$, which improves on a result of Chern-Lin \cite{CL5}, hence allowing us to construct sharper test functions and to prove existence of solutions for (\ref{one}) when $\gamma=\frac{n^2-1}{4}$. 

\begin{theorem} Assume $\gamma < \frac{n^2}{4}$ and let $u\in \dundeuxrnp$, $u\geq 0$, $u\not\equiv 0$ be a weak solution to
\begin{equation}\label{limiting}
-\Delta u-\frac{\gamma}{|x|^2}u=\frac{u^{\crits-1}}{|x|^s}\hbox{ in }\rnp.
\end{equation}
Then, there exist $K_1,K_2>0$ such that
\begin{equation}\nonumber
u(x)\sim_{x\to 0}K_1\frac{x_1}{|x|^{\am}}\quad \hbox{ and } \quad u(x)\sim_{|x|\to +\infty}K_2\frac{x_1}{|x|^{\ap}}.
\end{equation}
\end{theorem}

The above theorem yields in particular, the existence of a solution $u$ for (\ref{limiting}) which satisfies for some $C>0$, the estimates
\begin{equation}\label{bnd:U.0}
u(x)\leq C x_1 |x|^{-\ap} \quad \hbox{ and \quad $|\nabla u(x)|\leq C |x|^{-\ap}$ for all $x\in\rnp$. }
\end{equation}
Noting that 
\begin{eqnarray*}
\gamma < \frac{n^2-1}{4} & \Leftrightarrow & \ap-\am > 1, 
\end{eqnarray*}
it follows from   \eqref{bnd:U.0}, that whenever $\gamma<\frac{n^2-1}{4}$, then  $|x'|^2|\partial_1 u|^2=O(|x'|^{2-2\ap})$ as $|x'|\to +\infty$ on $\partial\rnp=\rr^{n-1}$, from which we could deduce that $x'\mapsto |x'|^2|\partial_1 u(x')|^2$ is in $L^1(\partial\rnp)$. This estimate --which does not hold when $\gamma\geq \frac{n^2-1}{4}$ --  is key for the construction of test functions for $\mu_{\gamma, s}(\Omega)$ based on the solution  $u$ of (\ref{limiting}), in the case when $\gamma\leq \frac{n^2-1}{4}$.  

\noindent In order to deal with the remaining cases for $\gamma$, that is when $\gamma \in (\frac{n^2-1}{4}, \frac{n^2}{4})$, we prove the following result which describes the general profile of any positive solution of $L_\gamma u=a(x)u$, albeit variational or not.  

\begin{theorem}[Classification of singular solutions]
Assume $\gamma < \frac{n^2}{4}$ and let $u\in C^2(B_\delta(0)\cap (\overline{\Omega}\setminus\{0\}))$ be such that
\begin{equation}
\left\{\begin{array}{ll}
-\Delta u-\frac{\gamma+O(|x|^\tau)}{|x|^2}u=0 &\hbox{\rm in } \Omega\cap B_\delta(0)\\
\hfill u > 0&\hbox{\rm in } \Omega\cap B_\delta(0)\\
\hfill u=0&\hbox{\rm on } (\partial\Omega\cap B_\delta(0))\setminus \{0\}.\end{array}\right.
\end{equation}
Then, there exists $K>0$ such that
\begin{equation}\label{profiles}
\hbox{either }u(x)\sim_{x\to 0}K\frac{d(x,\partial\Omega)}{|x|^{\am}}\quad \hbox{ or } \quad u(x)\sim_{x\to 0}K\frac{d(x,\partial\Omega)}{|x|^{\ap}}.
\end{equation}
In the first case, the solution $u$ is variational; in the second case, it is not.
\end{theorem}
This result then allows us to completely classify all positive solutions to $L_\gamma u=0$ on $\rnp$, a fact alluded to in Pinchover-Tintarev (\cite{PT}, Example 1.5).
\begin{proposition} Assume $\gamma < \frac{n^2}{4}$ and let $u\in C^2(\overline{\rnp}\setminus \{0\})$ be such that
\begin{equation}
\left\{\begin{array}{ll}
-\Delta u-\frac{\gamma}{|x|^2}u=0 &\hbox{\rm in } \rnp\\
u>0&\hbox{\rm in }  \rnp\\
u=0&\hbox{\rm on } \partial \rnp.\end{array}\right.
\end{equation}
Then, there exists $\lambda_-,\lambda_+\geq 0$ such that
\begin{equation}
u(x)=\lambda_- x_1|x|^{-\am}+\lambda_+ x_1|x|^{-\ap}\hbox{ for all }x\in \rnp.
\end{equation}
\end{proposition}
\noindent As mentioned above, the case when $\gamma>\frac{n^2-1}{4}$ is more intricate and requires isolating a new notion of  {\it singular  boundary mass} associated to the operator $L_\gamma$ for  domains  of $\rn$ having $0$ on their boundary.
 The following result will be proved in section \ref{sec:mass}. 

\begin{theorem}\label{def:mass:intro} Assume that $\frac{n^2-1}{4}<\gamma<\gamma_H(\Omega)$. Then, up to multiplication by a positive constant, there exists a unique function $H\in C^2(\overline{\Omega}\setminus \{0\})$ such that
\begin{equation}\label{carac:H.0}
-\Delta H-\frac{\gamma}{|x|^2}H=0\hbox{ in }\Omega\; ,\; H>0\hbox{ in }\Omega\; ,\; H=0\hbox{ on }\partial\Omega.
\end{equation}
Moreover, there exists $c_1>0$ and $c_2\in \rr$ such that
\begin{equation*}
\hbox{$H(x)=c_1\frac{d(x,\partial\Omega)}{|x|^{\ap}}+c_2\frac{d(x,\partial\Omega)}{|x|^{\am}}+ o\left(\frac{d(x,\partial\Omega)}{|x|^{\am}}\right)$\qquad as $x\to 0$. }
\end{equation*}
The quantity $m_\gamma(\Omega):=\frac{c_2}{c_1}\in \rr$, which is independent of the choice of $H$ satisfying \eqref{carac:H.0}, will be referred to as the Hardy singular b-mass of $\Omega$. 
\end{theorem}
Indeed, another interpretation of the threshold is the following. The case $\gamma>\frac{n^2-1}{4}$ is the only situation in which one can write a solution $H$ to \eqref{carac:H.0} as the sum of the two profiles given in \eqref{profiles} (plus lower-order terms) for any bounded domain $\Omega$. When $\gamma\leq\frac{n^2-1}{4}$, there might be some intermediate terms between the two profiles.

\medskip\noindent We show that the map $\Omega \to m_\gamma(\Omega)$ is a monotone increasing function on the class of domains having zero on their boundary, once ordered by inclusion. We shall also see below that it is possible to define the mass of some unbounded domains, and that $m_\gamma(\rnp)=0$ for any $\frac{n^2-1}{4}<\gamma<\frac{n^2}{4}$, from which follows that the mass of any one of its smooth subsets having zero on its boundary is non-positive. In particular, $m_\gamma(\Omega)<0$ whenever $\Omega$ is convex bounded and $0\in \partial \Omega$.

\medskip\noindent We shall however exhibit in section \ref{sec:ex:mass} examples of bounded domains $\Omega$ in $\rn$ with $0\in \partial \Omega$ and with positive mass. Among other things, we provide examples of domains with either positive or negative boundary mass, while satisfying any local behavior at $0$ one wishes. In other words, the sign of the Hardy b-mass is totally independent of the local properties of $\partial\Omega$ around $0$.

\medskip\noindent This notion and the preceeding results allow us to establish the following extension of the results of Chern-Lin.

\begin{theorem}  Let $\Omega$ be a bounded smooth domain of $\rr^n$ ($n\geq 3$) such that $0\in \partial \Omega$, hence $\frac{(n-2)^2}{4} <\gamma_H(\Omega) \leq \frac{n^2}{4}$. Let $0\leq s< 2$. 
\begin{enumerate}
\item If $\gamma_H(\Omega)\leq \gamma <\frac{n^2}{4}$, then there are extremals for $\mu_{\gamma,s}(\Omega)$ 
for all $n\geq 3$.
\item If $\gamma <\gamma_H(\Omega)$ and either $s>0$ or $\{s=0$, $n\geq 4$ and $\gamma >0\}$, then there are  extremals for $\mu_{\gamma,s}(\Omega)$, under either one of the following conditions:
\begin{itemize}
\item  $\gamma\leq\frac{n^2-1}{4}$ and the mean curvature of $\partial \Omega$ at $0$ is negative.
\item  $\gamma>\frac{n^2-1}{4}$ and the Hardy b-mass $m_\gamma(\Omega)$ is positive.
\end{itemize}  
 \item  If $\{s=0$ and $\gamma \leq 0\}$, then there are no extremals for $\mu_{\gamma,0}(\Omega)$ for any $n\geq 3$.
\end{enumerate}
\end{theorem}
\noindent We shall also address in section \ref{sec:n:3} the remaining case, i.e.,  $n=3$ and $s=0$ and $\gamma\in (0,\frac{n^2}{4})$. In this situation, there may or may not be extremals for $\mu_{\gamma, 0}(\rnp)$. If they do exist, we can then argue as before --using the same test functions-- to conclude existence of extremals under the same conditions, that is either $\gamma\leq\frac{n^2-1}{4}$ and the mean curvature of $\partial \Omega$ at $0$ is negative, or $\gamma>\frac{n^2-1}{4}$ and the mass $m_\gamma(\Omega)$ is positive. However, if no extremal exist for $\mu_{\gamma, 0}(\rnp)$, then as noted in \eqref{def:K}, we have that
\begin{equation*}
\mu_{\gamma, 0}(\rnp)= \inf_{u\in D^{1,2}(\rn)\setminus\{0\}}\frac{\int_{\rn}|\nabla u|^2\, dx}{\left(\int_{\rn}|u|^{\crit}\, dx\right)^{\frac{2}{\crit}}}=\frac{1}{K(n,2)^2},
\end{equation*}
and we are back to the case of the Yamabe problem with no boundary singularity. This means that 
one needs to resort to a more standard notion of mass $R_\gamma(\Omega, x_0)$ associated to $L_\gamma$ and an interior point $x_0\in \Omega$. 
One can then construct suitable  test-functions in the spirit of Schoen \cite{schoen1}. In order to define the ``internal mass", we show that for a given $\gamma\in (0,\gamma_H(\Omega))$, any solution $G$ of 
$$\left\{\begin{array}{ll}
-\Delta G-\frac{\gamma}{|x|^2} G=0 &\hbox{ in }\Omega\setminus \{x_0\}\\
\hfill G>0 &\hbox{ in }\Omega \setminus \{x_0\}\\
\hfill G=0 &\hbox{ on }\partial\Omega \setminus\{0\},
\end{array}\right.$$
 is unique up to multiplication by a constant, and that for any $x_0\in \Omega$, there exists $R_\gamma(\Omega, x_0)\in \rr$ (independent of $G$) and $c_G>0$ such that
$$G(x)=c_G\left(\frac{1}{|x-x_0|^{n-2}}+R_\gamma(\Omega, x_0)\right)+o(1)\quad \hbox{ as }x\to x_0.$$
The quantity $R_\gamma(\Omega, x_0)$ is well defined, and we prove the following.

\begin{theorem}\label{th:dim3}  Let $\Omega$ be a bounded smooth domain of $\rr^3$ such that $0\in \partial \Omega$. In particular $\frac{1}{4}<\gamma_H(\Omega) \leq \frac{9}{4}$.
\begin{enumerate}
\item If $\gamma_H(\Omega)\leq \gamma <\frac{9}{4}$, then there are extremals for $\mu_{\gamma,0}(\Omega)$.
\item If $0<\gamma <\gamma_H(\Omega)$, 
and if there exists $x_0\in\Omega$ such that $R_\gamma(\Omega, x_0)>0$, then there are  extremals for  $\mu_{\gamma,0}(\Omega)$, under either one of the following conditions:
\begin{enumerate}
\item  $\gamma\leq 2$ and the mean curvature of $\partial \Omega$ at $0$ is negative.
\item  $\gamma>2$ and the mass $m_\gamma(\Omega)$ is positive.
\end{enumerate}  
\end{enumerate}
\end{theorem}
More precisely, if there are extremals for $\mu_{\gamma,0}(\rr^3)$, then conditions (a) and (b) are sufficient to get extremals for $\mu_{\gamma,0}(\Omega)$. If there are no extremals for $\mu_{\gamma,0}(\rr^3)$, then the positivity of the  internal mass $R_\gamma(\Omega, x_0)$ is sufficient to get extremals for $\mu_{\gamma,0}(\Omega)$. We refer to Theorem \ref{th:3} for a precise statement. The following table summarizes our findings.

\begin{table}[ht]
\begin{center}
\caption{Singular Sobolev-Critical  term: $s>0$}

\vspace{2mm}
\begin{tabular}{|c|c|c|c|} \hline
{\bf  Hardy term }& {\bf Dimension}& {\bf Geometric condition} &{\bf  Extremal}\\ \hline\hline
$-\infty < \gamma \leq \frac{n^2-1}{4}$  & $n\geq 3$
  & Negative mean curvature at $0$& Yes\\  \hline
 $\frac{n^2-1}{4} <\gamma <\frac{n^2}{4}$&$n\geq 3$ & Positive boundary-mass & Yes\\
  \hline
\end{tabular}
\end{center}
\end{table}

\begin{table}[ht]
\begin{center}
\caption{ Non-singular Sobolev-Critical term: $s=0$}

\vspace{2mm}
\begin{tabular}{|c|c|c|c|} \hline
{\bf  Hardy term }& {\bf Dim.}& {\bf Geometric condition} &{\bf  Extr.}\\ \hline\hline
$0 < \gamma \leq \frac{n^2-1}{4}$  & $n = 3$
  & Negative mean curvature at $0$ \& Positive internal mass & Yes\\ 
  & $n\geq 4$&Negative mean curvature at $0$&Yes\\ \hline
 $\frac{n^2-1}{4} <\gamma <\frac{n^2}{4}$&$n= 3$ & Positive boundary-mass \& Positive internal mass & Yes\\
 &$n\geq 4$& Positive boundary mass& Yes\\\hline
$\gamma \leq 0$  & $n \geq 3 $& --  &  No\\
  \hline
\end{tabular}
\end{center}
\end{table}

\medskip\noindent{\bf Notations:} in the sequel, $C_i(a,b,...)$ ($i=1,2,...$) will denote constants depending on $a,b,...$. The same notation can be used for different constants, even in the same line.

\section{Old and new inequalities involving singular weights}\label{sec:hardy} The following general form of the Hardy inequality is well known. See for example the thesis of Cowan \cite{craig} or the book of Ghoussoub-Moradifam \cite{gm}.
We include here a proof for completeness. 
\begin{theorem} \label{cowan} Let $\Omega$ be a connected open subset of $\rn$ and consider $\rho\in C^\infty(\Omega)$ such that $\rho>0$ and $-\Delta\rho>0$. Then for any $u\in \dundeux$ we have that  $\sqrt{\rho^{-1}(-\Delta)\rho}u\in L^2(\Omega)$ and 
\begin{equation}
\int_\Omega\frac{-\Delta \rho}{\rho} u^2\, dx\leq \int_\Omega |\nabla u|^2\, dx.\label{ineq:hardy}
\end{equation}
Moreover, the case of equality is achieved exactly on $\rr\rho\cap \dundeux$. In particular, if $\rho\not\in \dundeux$, there are no nontrival extremals for \eqref{ineq:hardy}.
\end{theorem}
\noindent{\it Proof of Theorem \ref{cowan}:} The proof relies of the following integral identity:
\begin{equation}\label{id:ipp}
\int_\Omega |\nabla (\rho v)|^2\, dx-\int_\Omega\frac{-\Delta\rho}{\rho}(\rho v)^2\, dx=\int_\Omega\rho^2|\nabla v|^2\, dx\geq 0
\end{equation}
for all $v\in C^\infty_c(\Omega)$. This identity is a straightforward integration by parts. Since $\rho,-\Delta\rho>0$ in $\Omega$, it follows from density arguments that for any $u\in \dundeux$, then $\sqrt{\rho^{-1}(-\Delta)\rho}u\in L^2(\Omega)$ and \eqref{ineq:hardy} holds. 

\medskip\noindent Assume now that there exists $u_0\in \dundeux\setminus\{0\}$ that is an extremal for \eqref{ineq:hardy}. In other words, we have that
$$\int_\Omega\frac{-\Delta \rho}{\rho} u_0^2\, dx= \int_\Omega |\nabla u_0|^2\, dx.$$
Let $(u_i)_i\in C^\infty_c(\Omega)$ be such that $\lim_{i\to +\infty}u_i=u_0$ in $\dundeux$ and  define $v_i(x):=\frac{u_i(x)}{\rho(x)}$ for all $x\in\Omega$ and all $i$. This is well defined since $u_i$ has compact support in $\Omega$: therefore $v_i\in C^\infty_c(\Omega)$ for all $i$. Since $\dundeux\subset\dundeuxr$, Sobolev's embedding theorem yields convergence of $u_i$ to $u_0$ in $L^{2n/(n-2)}(\Omega)$. Since $\rho>0$ in $\Omega$, we then get that $(v_i)_i$ is uniformly bounded in $H_{1,loc}^2(\Omega)$. It then follows from reflexivity and a diagonal argument that there exists $v\in H_{1,loc}^2(\Omega)$ such that
$$\lim_{i\to +\infty}v_i=v\hbox{ in }H_{1,loc}^2(\Omega).$$
Applying \eqref{id:ipp} to $v_i=\rho^{-1}u_i$ yields $\lim_{i\to +\infty}\int_\omega\rho^2|\nabla v_i|^2\, dx=0$. Therefore, for any $\omega\subset\subset\Omega$, we have that
$$\int_\omega |\nabla v|^2\, dx\leq \liminf_{i\to +\infty}\int_\omega|\nabla v_i|^2\, dx=0.$$
Therefore $\int_\omega |\nabla v|^2\, dx=0$ for all $\omega\subset\subset\Omega$, and then there exists $c\in\rr$ such that $v\equiv c$. Up to extracting additional subsequence, we can assume that $u_i(x)$ and $v_i(x)$ converges respectively to $u_0(x)$ and $v(x)$ respectively when $i\to +\infty$ for a.e. $x\in \Omega$. Therefore, $u_0(x)=c\cdot\rho(x)$ for a.e. $x\in \Omega$. Since $u_0\not\equiv 0$, we have that $c\neq 0$ and then $\rho\in \dundeux$. For dimensional reasons, the equality is then achieved exactly on $\rr\rho\cap \dundeuxr$. This ends the case of equality in case there is a nontrivial extremal. 

\smallskip\noindent Assume now that $\rho\in\dundeux$. We let $(\rho_i)\in C^\infty_c(\Omega)$ such that $\lim_{i\to +\infty}\rho_i=\rho$ in $\dundeux$. Without loss of generality, we can assume that $\rho_i(x)\to \rho(x)$ as $i\to +\infty$ for a.e. $x\in\Omega$. We define $v_i:=\frac{\rho_i}{\rho}\in C^\infty_c(\Omega)$. We have that $v_i(x)\to 1$ as $i\to +\infty$ for a.e. $x\in\Omega$. For any $i, j$, \eqref{id:ipp} yields
\begin{equation*}\label{cauchy:vi}
\int_\Omega |\nabla (\rho_i-\rho_j)|^2\, dx-\int_\Omega\frac{-\Delta\rho}{\rho}(\rho_i-\rho_j)^2\, dx=\int_\Omega\rho^2|\nabla (v_i-v_j)|^2\, dx.
\end{equation*}
Therefore $(\rho\nabla v_i)_i$ is a Cauchy sequence in $L^2(\Omega,\rn)$, and therefore, there exists $\vec{X}\in L^2(\Omega,\rn)$ such that 
\begin{equation}\label{eq:X}
\lim_{i\to +\infty}\rho\nabla v_i=\vec{X}\hbox{ in }L^2(\Omega,\rn).
\end{equation}
Arguing as in the first part of the proof of Theorem \ref{cowan}, we get that there exists $v\in H_{1,loc}^2(\Omega)$ such that $\lim_{i\to +\infty}v_i=v$ in $H_{1,loc}^2(\Omega)$. Since $v_i(x)\to 1$ as $i\to +\infty$ for a.e. $x\in\Omega$, we get that $v\equiv 1$ and therefore $\nabla v=0$, which yields $\vec{X}=0$. It then follows from \eqref{eq:X} that $(\rho\nabla v_i)_i$ goes to $0$ in $L^2(\Omega,\rn)$. Using again \eqref{id:ipp} yields
\begin{equation}
\int_\Omega |\nabla \rho_i|^2\, dx-\int_\Omega\frac{-\Delta\rho}{\rho}\rho_i^2\, dx=\int_\Omega\rho^2|\nabla v_i|^2\, dx.
\end{equation}
Therefore, letting $i\to +\infty$ yields $\int_\Omega |\nabla \rho|^2\, dx=\int_\Omega\frac{-\Delta\rho}{\rho}\rho^2\, dx$, and then $\rho$ is an extremal for \eqref{ineq:hardy}.


\medskip\noindent The general case follows from the case of the existence of an extremal and the case $\rho\in \dundeux$. This ends the proof of Theorem \ref{cowan}.\hfill$\Box$

\medskip\noindent The above theorem applies to various weight functions $\rho$. See for example  the paper of Cowan \cite{craig} or the book \cite{gm}. For this paper, we need it for the following inequality. 
 
\begin{coro}\label{coro:1} Fix $1\leq k\leq n$, we then have the following inequality.
\begin{equation}\nonumber
\left(\frac{n+2k-2}{2}\right)^2=\inf_u \frac{\int_{\rr_+^k\times\rr^{n-k}} |\nabla u|^2\, dx}{\int_{\rr_+^k\times\rr^{n-k}}\frac{u^2}{|x|^2} \, dx},
\end{equation}
where the infimum is taken over all $u$ in $D^{1,2}(\rr_+^k\times\rr^{n-k})\setminus\{0\}$. Moreover, the infimum is never achieved.
\end{coro}
\noindent{\it Proof of Corollary \ref{coro:1}:} Take $\rho(x):=x_1...x_k|x|^{-\alpha}$ for all $x\in\Omega:=\rr_+^k\times\rr^{n-k}\setminus\{0\}$. Then $\frac{-\Delta \rho}{\rho}=\frac{\alpha(n+2k-2-\alpha)}{|x|^2}$. We then maximize the constant by taking $\alpha:=(n+2k-2)/2$. Since $\rho\not\in D^{1,2}(\rr_+^k\times\rr^{n-k})$, Theorem \ref{cowan} applies and we obtain that 
\begin{equation}\label{ineq:hardy:rk}
\left(\frac{n+2k-2}{2}\right)^2\int_{\rr_+^k\times\rr^{n-k}}\frac{u^2}{|x|^2} \, dx\leq \int_{\rr_+^k\times\rr^{n-k}} |\nabla u|^2\, dx
\end{equation}
for all $u\in D^{1,2}(\rr_+^k\times\rr^{n-k})$.
 
\medskip\noindent It remains to prove that the constant in \eqref{ineq:hardy:rk} is optimal. This will be achieved via the following test-function estimates. Construct a sequence $(\rho_\epsilon)_{\epsilon>0}\in D^{1,2}(\rr_+^k\times\rr^{n-k})$ as follows. Starting with  
$\rho(x)=x_1...x_k|x|^{-\alpha}$, 
we fix $\beta>0$ and define
\begin{equation}\label{def:rho:eps}
\rho_\epsilon(x):=\left\{\begin{array}{ll}
\left|\frac{x}{\epsilon}\right|^\beta\rho(x) &\hbox{ if }|x|<\epsilon\\
&\\
\rho(x)&\hbox{ if }\epsilon<|x|<\frac{1}{\epsilon}\\
&\\
|\epsilon\cdot x|^{-\beta} \rho(x) &\hbox{ if }|x|>\frac{1}{\epsilon}
\end{array}\right.\end{equation}
with $\alpha:=(n+2k-2)/2$. As one checks, $\rho_\epsilon\in D^{1,2}(\rr_+^k\times\rr^{n-k})$ for all $\epsilon>0$. The changes of variables $x=\epsilon y$ and $x=\epsilon^{-1}z$ yield
\begin{equation}\label{eq:rho:1}\begin{array}{cc}
\int_{B_\epsilon(0)}\frac{\rho_\epsilon^2}{|x|^2} \, dx= O(1), & \int_{B_\epsilon(0)}|\nabla\rho_\epsilon|^2 \, dx= O(1),\\
\int_{\rn\setminus \overline{B}_{\epsilon^{-1}}(0)}\frac{\rho_\epsilon^2}{|x|^2} \, dx= O(1), & \int_{\rn\setminus \overline{B}_{\epsilon^{-1}}(0)}|\nabla\rho_\epsilon|^2 \, dx= O(1)
\end{array}\end{equation}
when $\epsilon\to 0$. Integrating by parts yields
\begin{eqnarray}
\int_{B_{\epsilon^{-1}}(0)\setminus\overline{B}_\epsilon(0)}|\nabla\rho_\epsilon|^2 \, dx&=&\int_{B_{\epsilon^{-1}}(0)\setminus\overline{B}_\epsilon(0)}\frac{-\Delta \rho}{\rho} \rho^2\, dx+O(1)\nonumber\\
&=&\left(\frac{n+2k-2}{2}\right)^2\int_{B_{\epsilon^{-1}}(0)\setminus\overline{B}_\epsilon(0)}\frac{\rho^2}{|x|^2}\, dx+O(1),\label{eq:rho:2}
\end{eqnarray}
when $\epsilon\to 0$. Using polar coordinates yields
\begin{equation}\label{eq:rho:3}
\int_{B_{\epsilon^{-1}}(0)\setminus\overline{B}_\epsilon(0)}\frac{\rho^2}{|x|^2}\, dx= C(2)\ln\frac{1}{\epsilon}\hbox{ where }C(2):=2\int_{\mathbb{S}^{n-1}}\left|\prod_{i=1}^kx_i\right|^2\, d\sigma.
\end{equation}
Therefore, \eqref{eq:rho:1}, \eqref{eq:rho:2} and \eqref{eq:rho:3} yield
\begin{equation}\nonumber
\frac{\int_{\rr_+^k\times\rr^{n-k}}|\nabla\rho_\epsilon|^2\, dx}{\int_{\rr_+^k\times\rr^{n-k}}\frac{\rho_\epsilon^2}{|x|^2} \, dx}=\left(\frac{n+2k-2}{2}\right)^2+o(1)
\end{equation}
as $\epsilon\to 0$, and we are done. Note that the infimum is never achieved since $\rho\not\in D^{1,2}(\rr_+^k\times\rr^{n-k})$. This ends the proof of Corollary \ref{coro:1}.\hfill$\Box$

\medskip\noindent Another approach to prove Corollary \ref{coro:1} is to see $\rr_+^k\times\rr^{n-k}$ as a cone generated by a domain of the sphere. Then the Hardy constant is given by the Hardy constant of $\rn$ plus the first eigenvalue of the Laplacian of the above domain of the canonical sphere. This point of view is developed in Pinchover-Tintarev \cite{PT} (see also Fall-Musina \cite{FM} and Ghoussoub-Moradifam \cite{gm} for an exposition in book form).

\medskip\noindent We get the following generalized Caffarelli-Kohn-Nirenberg inequality.

\begin{proposition}\label{prop:ckn:gen} Let $\Omega$ be an open subset of $\rn$. Let $\rho,\rho'\in C^\infty(\Omega)$ be such that $\rho,\rho'>0$ and $-\Delta\rho,-\Delta\rho'>0$. Fix $s\in [0,2]$ and assume that there exists $\eps\in (0,1)$ and $\rho_\eps\in C^\infty(\Omega)$ such that
$$\frac{-\Delta \rho}{\rho}\leq (1-\eps)\frac{-\Delta\rho_\eps}{\rho_\eps}\hbox{ in }\Omega\hbox{ with }\rho_\eps,-\Delta\rho_\eps>0.$$
Then we have that
\begin{equation}\label{ineq:ckn:36}
\left(\int_\Omega \left(\frac{-\Delta \rho'}{\rho'}\right)^{s/2}\rho^{\crits}|u|^{\crits}\, dx\right)^{\frac{2}{\crits}}\leq C \int_{\Omega}\rho^{2}|\nabla u|^2\, dx
\end{equation}
for all $u\in C^\infty_c(\Omega)$.
\end{proposition}
\noindent{\it Proof of Proposition \ref{prop:ckn:gen}:} The Sobolev inequality yields the existence of $C(n)>0$ such that
$$\left(\int_\Omega |u|^{\crit}\, dx\right)^{\frac{2}{\crit}}\leq C(n)\int_{\Omega}|\nabla u|^2\, dx$$
for all $u\in C^\infty_c(\Omega)$, where $\crit=\crit(0)=\frac{2n}{n-2}$. A H\"older inequality interpolating between this Sobolev inequality and the Hardy inequality \eqref{ineq:hardy} for $\rho'$ yields the existence of $C> 0$ such that
\begin{equation}\label{ineq:ckn:pf}
\left(\int_\Omega \left(\frac{-\Delta \rho'}{\rho'}\right)^{s/2} |u|^{\crits}\, dx\right)^{\frac{2}{\crits}}\leq C\int_{\Omega}|\nabla u|^2\, dx
\end{equation}
for all $u\in C^\infty_c(\Omega)$. The identity \eqref{id:ipp} for $\rho$ and \eqref{ineq:hardy} for $\rho_\eps$ yield for $v\in C^\infty_c(\Omega)$,  
\begin{eqnarray*}
\int_\Omega \rho^2|\nabla v|^2\, dx &=& \int_\Omega |\nabla (\rho v)|^2\, dx-\int_\Omega \frac{-\Delta\rho}{\rho}(\rho v)^2\, dx\\
&\geq &\int_\Omega |\nabla (\rho v)|^2\, dx-(1-\eps)\int_\Omega \frac{-\Delta\rho_\eps}{\rho_\eps}(\rho v)^2\, dx\\
&\geq & \eps \int_\Omega |\nabla (\rho v)|^2
\end{eqnarray*}
Taking $u:=\rho v$ in \eqref{ineq:ckn:pf} and using this latest inequality yield \eqref{ineq:ckn:36}. This ends the proof of Proposition \ref{prop:ckn:gen}.\hfill$\Box$

\smallskip\noindent Here is an immediate consequence.
\begin{coro}\label{coro:ckn:2}
Fix $k\in\{1,\dots,n-1\}$. There exists then a constant $C:=C(a,b,n)>0$ such that for all $u\in C^\infty_c(\rr_+^k\times\rr^{n-k})$, the following inequality holds:
\begin{equation} \label{CKN:2}
\left(\int_{\rr_+^k\times\rr^{n-k}}|x|^{-bq}\left(\Pi_{i=1}^kx_i\right)^{q}|u|^q \right)^{\frac{2}{q}}\leq C\int_{\rr_+^k\times\rr^{n-k}}\left(\Pi_{i=1}^kx_i\right)^{2}|x|^{-2a}|\nabla u|^2 dx,
\end{equation}
where
\begin{equation}\label{cond2}
-\infty<a<\frac{n-2+2k}{2}, \ \ 0 \leq b-a\leq 1, \ \ q=\frac{2n}{n-2+2(b-a)}.
\end{equation}
\end{coro}

\noindent{\it Proof of Corollary \ref{coro:ckn:2}:} Define $\rho(x)=\rho'(x)=\left(\Pi_{i=1}^kx_i \right)|x|^{-a}$ and $\rho_\eps(x)=\left(\Pi_{i=1}^kx_i \right)|x|^{-\frac{n-2+2k}{2}}$ for all $x\in\rr_+^k\times\rr^{n-k}$. Here, we have that
$$\frac{\Delta\rho'}{\rho'}=\frac{a(n-2+2k-a)}{|x|^2}\hbox{ and }\frac{-\Delta\rho_\eps}{\rho_\eps}=\frac{(n-2+2k)^2}{4|x|^2}.$$
Apply Proposition \ref{prop:ckn:gen} with this data, with suitable $a,b,q$ to get Corollary \ref{coro:ckn:2}.\hfill$\Box$

\medskip\noindent{\bf Remark:} Observe that by taking $k=0$, we recover the classical Caffarelli-Kohn-Nirenberg inequalities \eqref{CKN}. However, one does not see any improvement in the integrability of the weight functions since  
$\left(\Pi_{i=1}^kx_i\right)|x|^{-a}$ is of order $k-a>-(n-2)/2$, hence as close as we wish to $(n-2)/2$. The relevance here appears when one considers the Hardy inequality of Corollary \ref{coro:1}. 

\section{Estimates for the best constant in the Hardy inequality \label{Hardy}}
As mentioned in the introduction, the best constant in the Hardy inequality 
$$\gamma_H(\Omega):=\inf\left\{\frac{\int_{\Omega}|\nabla u|^2\, dx}{\int_{\Omega}\frac{u^2}{|x|^2}\, dx} \, / \; u\in D^{1,2}(\Omega)\setminus\{0\}\right\}$$
does not depend on the domain $\Omega \subset \rn$ if the singularity $0$ belongs to the interior of $\Omega$. It is always equal to $\frac{(n-2)^2}{4}$. We have seen, however, in the last section that the situation changes whenever $0\in \partial \Omega$, since $\gamma_H(\rnp)=\frac{n^2}{4}$. Some properties of the best Hardy constants have been studied by Fall-Musina \cite{FM} and Fall \cite{F}. In this section, we shall collect whatever information we shall need later on about $\gamma_H$. 

\begin{proposition}\label{prop:gamma} $\gamma_H$ satisfies the following properties:
\begin{enumerate}
\item  For any smooth domain $\Omega$ such $0\in \Omega$, we have $\gamma_H(\Omega)=\frac{(n-2)^2}{4}$.
\item If $0\in \partial \Omega$, then $\frac{(n-2)^2}{4}<\gamma_H(\Omega)\leq\frac{n^2}{4}$. 
\item $\gamma_H(\Omega)=\frac{n^2}{4}$ for every $\Omega$ such that $0\in \partial \Omega$ and $\Omega \subset \rnp$.
\item If $\gamma_H(\Omega)<\frac{n^2}{4}$, then it is attained in $D^{1,2}(\Omega)$.  
\item We have $\inf\{\gamma_H(\Omega);\,  0\in\partial\Omega\}=\frac{(n-2)^2}{4}$ for $n\geq 3$.
\item For every $\epsilon>0$, there exists a smooth domain $\rnp\subsetneq \Omega_\epsilon \subsetneq \rn$ such that $0\in \partial\Omega_\epsilon$ and $\frac{n^2}{4}-\epsilon \leq \gamma_H(\Omega_\epsilon) < \frac{n^2}{4}$. 
\end{enumerate}

\end{proposition}

\noindent{\it Proof of Proposition \ref{prop:gamma}:} Properties (1)-(2)-(3)-(4) are well known (See Fall-Musina \cite{FM} and Fall \cite{F}).
However, we sketch proofs since we will make frequent use of the test functions involved. Note first that
Corollary \ref{coro:1} already yields that $\gamma_H(\rnp)=\frac{n^2}{4}$.

\medskip\noindent{\bf Proof of (2):} Since $\Omega\subset\rn$, we have that $\gamma_H(\Omega)\geq \gamma_H(\rn)=\frac{(n-2)^2}{4}$. Assume by contradiction that $\gamma_H(\Omega)=\frac{(n-2)^2}{4}$. It then follows from Theorem \ref{tool} below (applied with $s=2$) that $\gamma_H(\Omega)$ is achieved by a function in $u_0\in \dundeux\setminus \{0\}$ (note that $\mu_{0,\gamma}(\Omega)=\gamma_H(\Omega)-\gamma$). Therefore, $\gamma_H(\rn)$ is achieved in $D^{1,2}(\rn)$. Up to taking $|u_0|$, we can assume that $u_0\geq 0$. Therefore, the Euler-Lagrange equation and the maximum principle yield $u_0>0$ in $\rn$: this is impossible since $u_0\in \dundeux$. Therefore $\gamma_H(\Omega)>\frac{(n-2)^2}{4}$. 

\smallskip\noindent For the other inequality, the standard proof normally uses the fact that the domain contains an interior sphere that is tangent to the boundary at $0$. We choose here to perform another proof based on test-functions, which will be used again to prove Proposition \ref{prop:pptes:inf}. It goes as follows: since $\Omega$ is a smooth bounded domain of $\rn$ such that $0\in\partial\Omega$, there exists $U,V$ open subsets of $\rn$ such that $0\in U$, $0\in V$ and there exists $\varphi\in C^\infty(U,V)$ be a diffeomophism such that $\varphi(0)=0$ and 
$$\varphi(U\cap\{x_1>0\})=\varphi(U)\cap\Omega\hbox{ and }\varphi(U\cap\{x_1=0\})=\varphi(U)\cap\partial\Omega.$$
Moreover, we can and shall assume that $d\varphi_0$ is an isometry. Let $\eta\in C^\infty_c(U)$ such that $\eta(x)=1$ for $x\in B_\delta(0)$ for some $\delta>0$ small enough, and consider  $(\alpha_\epsilon)_{\epsilon>0}\in (0,+\infty)$ such that $\alpha_\epsilon=o(\epsilon)$ as $\epsilon\to 0$. For $\epsilon>0$, define
\begin{equation}\label{def:ue:ext}
\ue(x):=\left\{\begin{array}{ll}
\eta(y)\alpha_\epsilon^{-\frac{n-2}{2}}\rho_\epsilon\left(\frac{y}{\alpha_\epsilon}\right)&\hbox{ for all }x\in \varphi(U)\cap\Omega,\; x=\varphi(y),\\
0&\hbox{elsewhere}.\end{array}\right.
\end{equation}
Here $\rho_\epsilon$ is constructed as in \eqref{def:rho:eps} with $k=1$. Now fix $\sigma\in [0,2]$, and note that only the case $\sigma=2$ is needed for the above proposition. We then have  as $\epsilon\to 0$, 
\begin{eqnarray*}
\int_\Omega \frac{|\ue(y)|^{\crit(\sigma)}}{|y|^\sigma}\, dy&=& \int_{\rnp}\frac{\ue\circ\varphi(x)^{\crit(\sigma)}}{|\varphi(x)|^\sigma}|\hbox{Jac}(\varphi)(x)|\, dx\\
&=& \int_{\rnp}\frac{\ue\circ\varphi(x)^{\crit(\sigma)}}{|x|^\sigma}|(1+O(|x|))\, dx\\
&=& \int_{B_\delta(0)\cap \rnp}\frac{\ue\circ\varphi(x)^{\crit(\sigma)}}{|x|^\sigma}(1+O(|x|))\, dx+O(1).
\end{eqnarray*}
Dividing $B_\delta(0)=\left(B_\delta(0)\setminus B_{\epsilon^{-1}\alpha_\epsilon}(0)\right)\cup\left(B_{\epsilon^{-1}\alpha_\epsilon}(0)\setminus B_{\epsilon\alpha_\epsilon}(0)\right)\cup B_{\epsilon\alpha_\epsilon}(0)$ and arguing as in \eqref{eq:rho:1} to \eqref{eq:rho:3}, we get as $\epsilon\to 0$, 
\begin{eqnarray*}
\int_\Omega \frac{|\ue(y)|^{\crit(\sigma)}}{|y|^\sigma}\, dy&=& \int_{\left[B_{\epsilon^{-1}\alpha_\epsilon}(0)\setminus B_{\epsilon\alpha_\epsilon}(0)\right]\cap \rnp}\frac{\ue\circ\varphi(x)^{\crit(\sigma)}}{|x|^\sigma}(1+O(|x|))\, dx+O(1)\\
&=& \int_{\left[B_{\epsilon^{-1}\alpha_\epsilon}(0)\setminus B_{\epsilon\alpha_\epsilon}(0)\right]\cap \rnp}\frac{\ue\circ\varphi(x)^{\crit(\sigma)}}{|x|^\sigma}\, dx+O(1)\\
&=& \int_{\left[B_{\epsilon^{-1}\alpha_\epsilon}(0)\setminus B_{\epsilon\alpha_\epsilon}(0)\right]\cap \rnp}\frac{\rho(x)^{\crit(\sigma)}}{|x|^\sigma}\, dx+O(1).
\end{eqnarray*}
Passing to polar coordinates yields
\begin{equation}\label{eq:ue:ext:1}
\int_\Omega \frac{|\ue(y)|^{\crit(\sigma)}}{|y|^\sigma}\, dy= C(\sigma)\ln\frac{1}{\epsilon} +O(1)\quad \hbox{as $\epsilon\to 0$,}
\end{equation}
where $C(\sigma):=2\int_{\mathbb{S}^{n-1}}\left|\prod_{i=1}^kx_i\right|^{\crit (\sigma)}\, d\sigma.$

\smallskip\noindent Similar arguments yield
\begin{eqnarray*}
\int_\Omega |\nabla \ue|^2\, dy&= &\int_{B_{\epsilon^{-1}\alpha_\epsilon}(0)\setminus B_{\epsilon\alpha_\epsilon}(0)\cap \rnp}|\nabla \ue\circ\varphi(x)|^2(1+O(|x|)\, dx+O(1)\\
&= &\int_{B_{\epsilon^{-1}\alpha_\epsilon}(0)\setminus B_{\epsilon\alpha_\epsilon}(0)\cap \rnp}|\nabla \ue\circ\varphi(x)|^2\, dx+O(1)\\
&= &\int_{B_{\epsilon^{-1}}(0)\setminus B_{\epsilon}(0)\cap \rnp}|\nabla \rho(x)|^2\, dx+O(1)
\end{eqnarray*}
as $\epsilon\to 0$. Using \eqref{eq:rho:2} and \eqref{eq:rho:3} yield
\begin{equation}\label{eq:ue:ext:2}
\int_\Omega |\nabla \ue|^2\, dy=\frac{n^2}{4} C(2)\ln\frac{1}{\epsilon}+O(1) \quad \hbox{as $\epsilon\to 0$. }
\end{equation}
As a consequence, we get that
\begin{equation}\nonumber
\frac{\int_\Omega|\nabla \ue|^2\, dx}{\int_\Omega\frac{\ue^2}{|x|^2}\, dx}=\frac{n^2}{4}+o(1)\quad \hbox{as $\epsilon\to 0$. }
\end{equation}
In particular, we get that
$\gamma_H(\Omega)\leq \frac{n^2}{4}$,
which proves the large inequality in (2). 

\medskip\noindent{\bf Proof of (3).} Assume that $\Omega\subset\rnp$, then $\dundeux\subset D^{1,2}(\rnp)$, and therefore $\gamma_H(\Omega)\geq \gamma_H(\rnp)=n^2/4$. With the reverse inequality already given by (2), we get that $\gamma_H(\Omega)=n^2/4$ for all $\Omega\subset\rnp$ such that $0\in\partial\Omega$.

\medskip\noindent{\bf Proof of (4).} This will be a particular case of Theorem \ref{tool} when $s=2$. 

\medskip\noindent{\bf Proof of (5).} Let $\Omega_0$ be a bounded domain of $\rn$ such that $0\in\Omega_0$ (i.e., it is not on the boundary). Given $\delta>0$, we chop out a ball of radius $\delta/4$ with $0$ on its boundary to define
$$\Omega_\delta:=\Omega_0\setminus \overline{B}_{\frac{\delta}{4}}\left((\frac{-\delta}{4}, 0,\dots,0)\right)$$
Note that for $\delta>0$ small enough, $\Omega$ is smooth and $0\in\partial\Omega$. We now prove that
\begin{equation}\label{lim:gamma:1}
\lim_{\delta\to 0}\gamma_H(\Omega_\delta)=\frac{(n-2)^2}{4}.
\end{equation}
Define $\eta_1\in C^\infty(\rn)$ such that
$$\eta_1(x)=\left\{\begin{array}{ll}
0& \hbox{ if }|x|<1\\
1& \hbox{ if }|x|>2,
\end{array}\right.$$
and let $\eta_\delta(x):=\eta_1(\delta^{-1}x)$ for all $\delta>0$ and $x\in\rn$. Fix $U\in C^\infty_c(\rn)$ and consider for any $\delta>0$, an $\eps_\delta>0$ such that
$$\lim_{\delta\to 0}\frac{\delta}{\eps_\delta}=\lim_{\delta\to 0}\eps_\delta=0.$$
For $\delta>0$, we define
$$u_{\delta}(x):=\eta_\delta(x)\eps_\delta^{-\frac{n-2}{2}}U(\eps_\delta^{-1}x)\hbox{ for all }x\in \Omega_\delta.$$
For $\delta>0$ small enough, we have that $u_{\delta}\in C^\infty_c(\Omega_\delta)$. A change of variable yields
$$\int_{\Omega_\delta}\frac{u_{\delta}^2}{|x|^2}\, dx=\int_{\rn}\frac{U^2}{|x|^2}\eta_1^2\left(\frac{\eps_\delta x}{\delta}\right)\, dx$$
for all $\delta>0$ small enough. Since $\delta=o(\eps_\delta)$ as $\delta\to 0$,
the dominated convergence theorem yields
\begin{equation}\nonumber
\lim_{\delta\to 0}\int_{\Omega_\delta}\frac{u_{\delta}^2}{|x|^2}\, dx=\int_{\rn}\frac{U^2}{|x|^2}\, dx.
\end{equation}
For $\eps>0$ small enough, we have that
\begin{eqnarray}
\int_{\Omega_\delta}|\nabla u_{\delta}|^2\, dx& =& \int_{\rn}|\nabla u_{\delta}|^2\, dx= \int_{\rn}|\nabla \left(U\cdot \eta_{\frac{\delta}{\eps_\delta}}\right)|^2\, dx\nonumber\\
&=& \int_{\rn}|\nabla U|^2 \eta_{\frac{\delta}{\eps_\delta}}^2\, dx+\int_{\rn} \eta_{\frac{\delta}{\eps_\delta}}\left(-\Delta  \eta_{\frac{\delta}{\eps_\delta}}\right) U^2\, dx.\label{eq:U:1}
\end{eqnarray}
Let $R>0$ be such that $U$ has support in $B_R(0)$. We then have that
\begin{eqnarray*}
\int_{\rn} \eta_{\frac{\delta}{\eps_\delta}}\left(-\Delta\eta_{\frac{\delta}{\eps_\delta}}\right)   U^2\, dx&=&O\left(\left(\frac{\eps_\delta}{\delta}\right)^2 \hbox{Vol}(B_R(0)\cap \hbox{Supp }\left(-\Delta  \eta_{\frac{\delta}{\eps_\delta}}\right))\right)\\
&=&O\left(\left(\frac{\delta}{\eps_\delta}\right)^{n-2}\right)=o(1)
\end{eqnarray*}
as $\delta\to 0$ since $n\geq 3$. This latest identity, \eqref{eq:U:1} and the dominated convergence theorem  yield
\begin{equation}\nonumber
\lim_{\delta\to 0}\int_{\Omega_\delta}|\nabla u_{\delta}|^2\, dx=\int_{\rn}|\nabla U|^2\, dx.
\end{equation}
Therefore, for $U\in C^\infty_c(\rn)$, we have
$$\limsup_{\delta\to 0}\gamma_H(\Omega_\delta)\leq \lim_{\delta\to 0}\frac{\int_{\Omega_\delta}|\nabla u_{\delta}|^2\, dx}{\int_{\Omega_\delta}\frac{u_{\delta}^2}{|x|^2}\, dx}=\frac{\int_{\rn}|\nabla U|^2\, dx}{\int_{\rn}\frac{U^2}{|x|^2}\, dx}.$$
Taking the infimum over all $U\in C^\infty_c(\rn)$, we get that
$$\limsup_{\delta\to 0}\gamma_H(\Omega_\delta)\leq \inf_{U\in D^{1,2}(\rn)\setminus\{0\}}\frac{\int_{\rn}|\nabla U|^2\, dx}{\int_{\rn}\frac{U^2}{|x|^2}\, dx}=\gamma_H(\rn)=\frac{(n-2)^2}{4}.$$
Since $\gamma_H(\Omega_\delta)\geq \frac{(n-2)^2}{4}$ for all $\delta>0$, this completes the proof of \eqref{lim:gamma:1}, yielding (5).


\medskip\noindent{\bf Proof of (6).}  The proof uses the following observation.

\begin{lem}\label{lem:ex:gamma} Let $(\Phi_k)_{k\in\nn}\in C^1(\rn,\rn)$ be such that
\begin{equation}\label{hyp:2}
\lim_{k\to +\infty}\left(\Vert \Phi_k-Id_{\rn}\Vert_\infty+\Vert \nabla(\Phi_k-Id_{\rn})\Vert_\infty\right)=0\hbox{ and }\Phi_k(0)=0.
\end{equation}
Let $D\subset \rn$ be an open domain such that $0\in\partial D$ (the domain is not necessarily bounded nor regular), and set $D_k:=\Phi_k(D)$ for all $k\in\nn$. Then $0\in \partial D_k$ for all $k\in\nn$ and 
\begin{equation}\label{lim:gamma:h}
\lim_{k\to +\infty}\gamma_H(D_k)=\gamma_H(D).
\end{equation} 
\end{lem}
\noindent {\it Proof of Lemma \ref{lem:ex:gamma}:} 
If $u\in C^\infty_c(D_k)$, then $u\circ\Phi_k\in C^\infty_c(D)$ and
\begin{eqnarray*}
\int_{D_k}|\nabla u|^2\, dx &=& \int_{\rnp}|\nabla (u\circ\Phi_k)|^2_{\Phi_k^\star\eucl}|\hbox{Jac}(\Phi_k)|\, dx,\\
\int_{D_k}\frac{u^2}{|x|^2}\, dx &=& \int_{\rnp}\frac{(u\circ\Phi_k(x))^2}{|\Phi_k(x)|^2}|\hbox{Jac}(\Phi_k)|\, dx.
\end{eqnarray*}
Assumption \eqref{hyp:2} yields that
$$\lim_{k\to +\infty}\sup_{x\in D}\left(\left|\frac{|\Phi_k(x)|}{|x|}-1\right|+\sup_{i,j}\left|\left(\partial_i\Phi_k(x),\partial_j\Phi_k(x)\right)-\delta_{ij})\right| +|\hbox{Jac}(\Phi_k)-1|\right)=0,$$
where $\delta_{ij}=1$ if $i=j$ and $0$ otherwise. Therefore, for any $\eps>0$, there exists $k_0$ such that for all $u\in C^\infty_c(D_k)$ and $k\geq k_0$, 
\[
(1+\eps)\int_{D}|\nabla (u\circ\Phi_k)|^2\, dx\geq \int_{D_k}|\nabla u|^2\, dx \geq  (1-\eps)\int_{D}|\nabla (u\circ\Phi_k)|^2\, dx,
\]
and
\[
(1+\eps)\int_{D}\frac{(u\circ\Phi_k(x))^2}{|x|^2}\, dx \geq \int_{D_k}\frac{u^2}{|x|^2}\, dx \geq  (1-\eps)\int_{D}\frac{(u\circ\Phi_k(x))^2}{|x|^2}\, dx.
\]
We can now deduce \eqref{lim:gamma:h} by using a standard density argument. This completes  the proof of Lemma \ref{lem:ex:gamma}.\hfill$\Box$

\medskip\noindent We now prove (6) of Proposition \ref{prop:gamma}. Let $\varphi\in C^\infty(\rr^{n-1})$ such that $0\leq\varphi\leq 1$, $\varphi(0)=0$, and $\varphi(x')=1$ for all $x'\in\rr^{n-1}$ such that $|x'|\geq 1$. For $t\geq 0$, define $\Phi_t(x_1,x'):=(x_1-t\varphi(x'), x')$ for all $(x_1,x')\in\rn$. Set $\tilde{\Omega}_t:=\Phi_t(\rnp)$ and apply Lemma \ref{lem:ex:gamma} to note that $\lim_{\eps\to 0}\gamma_H(\tilde{\Omega}_t)=\gamma_H(\rnp)=\frac{n^2}{4}$. Since $\varphi\geq 0$, $\varphi\not\equiv 0$, we have that $\rnp\subsetneq \tilde{\Omega}_t$ for all $t>0$. To get (6) it suffices to take $\Omega_\eps:=\tilde{\Omega}_t$ for $t>0$ small enough.

\section{Estimates on the best constants in the Hardy-Sobolev inequalities \label{compactness}}
As in the case of  $\gamma_H(\Omega)$, the best Hardy-Sobolev constant
$$\mu_{\gamma, s}(\Omega):=\inf\left\{\frac{\int_{\Omega} |\nabla u|^2\, dx-\gamma \int_{\Omega}\frac{u^2}{|x|^2}dx}{(\int_{\Omega}\frac{u^{\crits}}{|x|^s}dx)^{\frac{2}{\crits}}};\, u\in \dundeux\setminus\{0\}\right\}$$
will depend on the geometry of $\Omega$ whenever $0\in\partial\Omega$. In this section, we collect general facts that will be used throughout the paper.

\begin{proposition}\label{prop:pptes:inf} Let $\Omega$ be a bounded smooth domain such that $0\in \partial \Omega$.
\begin{enumerate}
\item If $\gamma <\frac{n^2}{4}$, then 
$\mu_{\gamma, s}(\Omega)>-\infty.$
\item If $\gamma >\frac{n^2}{4}$, then 
$\mu_{\gamma, s}(\Omega)=-\infty.$\par
\noindent Moreover, 
\item If $\gamma<\gamma_H(\Omega)$, then $\mu_{\gamma, s}(\Omega)>0$. 
\item If $\gamma_H(\Omega)<\gamma<\frac{n^2}{4}$, then $0>\mu_{\gamma, s}(\Omega)>-\infty$.
\item If $\gamma=\gamma_H(\Omega)<\frac{n^2}{4}$, then $\mu_{\gamma, s}(\Omega)=0$.
\end{enumerate}
\end{proposition}

 

\noindent{\it Proof of Proposition \ref{prop:pptes:inf}:} We first assume that $\gamma<\frac{n^2}{4}$. Let $\epsilon>0$ be such that $(1+\epsilon)\gamma\leq \frac{n^2}{4}$. It follows from Proposition \ref{prop:ineq:eps} that there exists $C_\epsilon>0$ such that for all $u\in\dundeux$, 
$$\frac{n^2}{4}\int_\Omega\frac{u^2}{|x|^2}\, dx\leq (1+\epsilon)\int_\Omega|\nabla u|^2\, dx+C_\epsilon\int_\Omega u^2\, dx.$$
For any $u\in\dundeux\setminus \{0\}$, we have 
\begin{eqnarray*}
J^\Omega_{\gamma, s}(u)&\geq& \frac{\left(1-\frac{4\gamma}{n^2}(1+\epsilon)\right)\int_\Omega |\nabla u|^2\, dx-\frac{4\gamma}{n^2}C_\epsilon \int_\Omega u^2\, dx}{\left(\int_\Omega\frac{|u|^{\crits}}{|x|^s}\, dx\right)^{\frac{2}{\crits}}}\\
&\geq& -\frac{4\gamma}{n^2}C_\epsilon\frac{\int_\Omega u^2\, dx}{\left(\int_\Omega\frac{|u|^{\crits}}{|x|^s}\, dx\right)^{\frac{2}{\crits}}}.
\end{eqnarray*} 
It follows from H\"older's inequality that there exists $C>0$ independent of $u$ such that $\int_\Omega u^2\, dx\leq \left(\int_\Omega\frac{|u|^{\crits}}{|x|^s}\, dx\right)^{\frac{2}{\crits}}$. It then follows that $J^\Omega_{\gamma, s}(u)\geq -\frac{4\gamma}{n^2}C_\epsilon C$ for all $u\in\dundeux\setminus \{0\}$. Therefore $\mu_{\gamma, s}(\Omega)>-\infty$ whenever $\gamma<\frac{n^2}{4}$.

\medskip\noindent Assume now that $\gamma>\frac{n^2}{4}$ and define for every $\epsilon>0$ a function $\ue\in \dundeux$ as in \eqref{def:ue:ext}. It then follows from \eqref{eq:ue:ext:1} and \eqref{eq:ue:ext:2} that as $\epsilon\to 0$, 
$$J^\Omega_{\gamma,s}(\ue)=\frac{\left(\frac{n^2}{4}-\gamma\right)C(2)\ln\frac{1}{\epsilon}+O(1)}{\left(C(s)\ln\frac{1}{\epsilon}+O(1)\right)^{\frac{2}{\crits}}}=\left(\left(\frac{n^2}{4}-\gamma\right)\frac{C(2)}{C(s)^{\frac{2}{\crits}}}+o(1)\right)\left(\ln\frac{1}{\epsilon}\right)^{\frac{2-s}{n-s}}.$$
Since $s<2$ and $\gamma>\frac{n^2}{4}$, we then get that $\lim_{\epsilon\to 0}J^\Omega_{\gamma,s}(\ue)=-\infty$, and therefore $\mu_{\gamma, s}(\Omega)=-\infty$.

\medskip\noindent Now assume that $\gamma<\gamma_H(\Omega)$. Sobolev's embedding theorem yields that  $\mu_{0,s}(\Omega)>0$, hence the result is clear for all $\gamma\leq 0$. If now $0\leq \gamma<\gamma_H(\Omega)$, it  follows from the definition of $\gamma_H(\Omega)$ that for all $u\in \dundeux\setminus\{0\}$, 
\begin{eqnarray*}
J^\Omega_{\gamma,s}(u)=\frac{\int_{\Omega} |\nabla u|^2-\gamma \int_{\Omega}\frac{u^2}{|x|^2}dx}{(\int_{\Omega}\frac{u^{\crits}}{|x|^s}dx)^{\frac{2}{\crits}}}
&\geq& \left(1-\frac{\gamma}{\gamma_H(\Omega)}\right)\frac{\int_\Omega|\nabla u|^2\, dx}{\left(\int_\Omega\frac{|u|^{\crits}}{|x|^s}\, dx\right)^{\frac{2}{\crits}}}\\
&\geq& \left(1-\frac{\gamma}{\gamma_H(\Omega)}\right)\mu_{0,s}(\Omega).
\end{eqnarray*}
Therefore $\mu_{\gamma,s}(\Omega)\geq \left(1-\frac{\gamma}{\gamma_H(\Omega)}\right)\mu_{0,s}(\Omega)>0$ when $\gamma<\gamma_H(\Omega)$.\\

\noindent We now assume that $\gamma_H(\Omega)<\gamma<\frac{n^2}{4}$. It follows from Proposition \ref{prop:gamma} (4), that  $\gamma_H(\Omega)$ is attained. We let $u_0$ be such an extremal. In particular 
$J^\Omega_{\gamma_H(\Omega),s}(u)\geq 0=J^\Omega_{\gamma_H(\Omega),s}(u_0)$, and 
 therefore $\mu_{\gamma_H(\Omega),s}(\Omega)=0$. Since $\gamma_H(\Omega)<\gamma<\frac{n^2}{4}$, then $J^\Omega_{\gamma,s}(u_0)<0$, and therefore $\mu_{\gamma,s}(\Omega)<0$ when $\gamma_H(\Omega)<\gamma<\frac{n^2}{4}$.
\noindent This ends the proof of Proposition \ref{prop:pptes:inf}.\hfill$\Box$

\begin{rem} \rm The case $\gamma=\frac{n^2}{4}$ is unclear and anything can happen at that value of $\gamma$. For example, if $\gamma_H(\Omega)<\frac{n^2}{4}$ then $\mu_{\frac{n^2}{4},s}(\Omega)<0$, while if $\gamma_H(\Omega)=\frac{n^2}{4}$ then $\mu_{\frac{n^2}{4},s}(\Omega)\geq 0$. It is our guess that many examples reflecting different regimes can be constructed. 
\end{rem}

\begin{proposition}\label{prop:ineq:eps} Assume $\gamma<\frac{n^2}{4}$ and $s\in [0,2]$. Then, for any $\eps>0$, there exists $C_\eps>0$ such that for all $u\in\dundeux$, 
\begin{equation}\label{ineq:sobo:eps:bis}
\left(\int_\Omega \frac{|u|^{\crits}}{|x|^s}\, dx\right)^{\frac{2}{\crits}}\leq \left(\frac{1}{\mu_{\gamma,s}(\rnp)}+\epsilon\right)\int_\Omega\left(|\nabla u|^2-\gamma\frac{u^2}{|x|^2}\right)\, dx+C_\eps\int_\Omega u^2\, dx.
\end{equation}
\end{proposition}
\noindent{\it Proof of  Proposition \ref{prop:ineq:eps}:} Fix $\epsilon>0$.
We first claim that there exists $\delta_\epsilon>0$ such that for all $u\in C^\infty_c(\Omega\cap B_{\delta_\epsilon}(0))$, 
\begin{equation}\label{eq:step:1}
\left(\int_\Omega \frac{|u|^{\crits}}{|x|^s}\, dx\right)^{\frac{2}{\crits}}\leq (\mu_{\gamma,s}(\rnp)^{-1}+\epsilon)\int_\Omega\left(|\nabla u|^2-\gamma\frac{u^2}{|x|^2}\right)\, dx.
\end{equation}
Indeed, for two open subsets of $\rn$ containing $0$, we may define a diffeomorphism $\varphi: U\to V$ 
such that $\varphi(0)=0$, $\varphi(U\cap\rnp)=\varphi(U)\cap \Omega$ and $\varphi(U\cap\partial\rnp)=\varphi(U)\cap \partial\Omega$. Moreover, we can also assume that $d\varphi_0$ is a linear isometry. In particular 
\begin{equation}\label{eq:pf:1.00}
|\varphi^\star\eucl-\eucl|(x)\leq C|x|\hbox{ and }|\varphi(x)|=|x|\cdot(1+O(|x|))
\end{equation}
for $x\in U$. If now $u\in C^\infty_c(\varphi(B_\delta(0))\cap\Omega)$, then $v:=u\circ\varphi\in C^\infty_c(B_\delta(0)\cap\rnp)$. If $g:=\varphi^{-1\star}\eucl$ denotes the metric induced by $\phi$, then we get  from  \eqref{eq:pf:1.00}, 
 \begin{eqnarray}\label{eq:pf:2.00}
\qquad \left(\int_\Omega \frac{|u|^{\crits}}{|x|^s}\, dx\right)^{\frac{2}{\crits}}&\leq&\left(\int_{B_\delta(0)\cap\rnp} \frac{|v|^{\crits}}{|\varphi(x)|^s}|\hbox{Jac } \varphi(x)|\, dx\right)^{\frac{2}{\crits}}\nonumber\\
&\leq & (1+C\delta)\left(\int_{B_\delta(0)\cap\rnp} \frac{|v|^{\crits}}{|x|^s}\, dx\right)^{\frac{2}{\crits}}\nonumber\\
&\leq & (1+C\delta)\mu_{\gamma,s}(\rnp)^{-1}\int_{B_\delta(0)\cap\rnp}\left(|\nabla v|^2-\gamma\frac{v^2}{|x|^2} \right)\, dx\nonumber\\
&\leq & \frac{1+C\delta}{\mu_{\gamma,s}(\rnp)}\int_{\varphi(B_\delta(0))\cap\Omega}\left(|\nabla u|_{g}^2-\frac{\gamma u^2}{|\varphi^{-1}(x)|^2}\right)|\hbox{Jac } \varphi^{-1}(x)|\, dx\nonumber\\
&\leq&  (1+C_1\delta)\mu_{\gamma,s}(\rnp)^{-1}\int_{\Omega}\left(|\nabla u|^2-\gamma\frac{u^2}{|x|^2}\right)\, dx\nonumber\\
&& +C_2\delta \int_{\Omega}\left(|\nabla u|^2
+\frac{u^2}{|x|^2}\right)\, dx. 
\end{eqnarray}
We also have that
\begin{eqnarray*}
\int_\Omega\frac{u^2}{|x|^2}\, dx&=&\int_{\varphi(B_\delta(0)\cap\rnp)}\frac{u^2}{|x|^2}\, dx=\int_{B_\delta(0)\cap\rnp}\frac{v^2}{|\varphi(x)|^2}|\hbox{Jac}(\varphi)(x)|\, dx\\
&=&\int_{B_\delta(0)\cap\rnp}\frac{v^2}{|x|^2}(1+O(|x|)\, dx\leq (1+C_1\delta) \int_{\rnp}\frac{v^2}{|x|^2}\, dx
\end{eqnarray*}
and
\begin{eqnarray*}
\int_\Omega|\nabla u|^2\, dx&=&\int_{\varphi(B_\delta(0)\cap\rnp)}|\nabla u|^2\, dx=\int_{B_\delta(0)\cap\rnp}|\nabla v|^2_{\varphi^\star\eucl}|\hbox{Jac}(\varphi)(x)|\, dx\\
&=&\int_{B_\delta(0)\cap\rnp}|\nabla v|^2(1+O(|x|)\, dx\geq (1-C_2\delta) \int_{\rnp}|\nabla v|^2\, dx,
\end{eqnarray*}
where $C_1,C_2>0$ are independent of $\delta$ and $v$. Hardy's inequality \eqref{ineq:hardy:rk} then yields for all $u\in C^\infty_c(\varphi(B_\delta(0)\cap\rnp))$,
\begin{equation}\label{ineq:hardy:opt:1}
\frac{n^2}{4}\int_\Omega\frac{u^2}{|x|^2}\, dx\leq \frac{1+C_1\delta}{1-C_2\delta}\int_\Omega|\nabla u|^2\, dx\leq (1+C_3\delta)\int_{\Omega}|\nabla u|^2\, dx.
\end{equation}
Since $\gamma<\frac{n^2}{4}$, there exists then $c>0$ such that for $\delta>0$ small enough,
\begin{equation*}
c^{-1}\int_{\Omega}|\nabla u|^2\, dx\leq \int_{\Omega}\left(|\nabla u|^2-\gamma\frac{u^2}{|x|^2}\right)\, dx\leq c\int_{\Omega}|\nabla u|^2\, dx
\end{equation*}
for all $u\in C^\infty_c(\varphi(B_\delta(0))\cap\Omega)$. Plugging these latest inequalities in \eqref{eq:pf:2.00} yields \eqref{eq:step:1} by taking $\delta_\epsilon$ small enough. 

\medskip\noindent Consider now $\eta\in C^\infty(\rn)$ such that  $\sqrt{\eta},\sqrt{1-\eta}\in C^2(\rn)$, 
such that $\eta(x)=1$ for $x\in B_{\delta_\eps/2}(0)$ and $\eta(x)=0$ for $x\not\in B_{\delta_\eps}(0)$. We shall use the notation 
\[
\|w\|_{{p,|x|^{-s}}}=\left(\int_\Omega \frac{|w|^p}{|x|^s}\, dx\right)^{1/p}.
\]
For $u\in C^\infty_c(\Omega)$, use  H\"older's inequality to write
\begin{eqnarray*}
\left(\int_\Omega \frac{|u|^{\crits}}{|x|^s}\, dx\right)^{\frac{2}{\crits}}&=&\Vert u^2\Vert_{\frac{\crits}{2},|x|^{-s}}=\Vert \eta u^2+(1-\eta)u^2\Vert_{\frac{\crits}{2},|x|^{-s}}\\
&\leq & \Vert \eta u^2\Vert_{\frac{\crits}{2},|x|^{-s}} +\Vert (1-\eta) u^2\Vert_{\frac{\crits}{2},|x|^{-s}}\\
&\leq & \Vert \sqrt{\eta} u\Vert^2_{\crits,|x|^{-s}} +\Vert \sqrt{1-\eta} u\Vert^2_{\crits,|x|^{-s}}.
\end{eqnarray*} 
Since $\sqrt{\eta}u\in C^\infty_c(B_{\delta_\eps}(0)\cap\Omega)$, it follows from inequality \eqref{eq:step:1} 
 that
\begin{eqnarray}
\left(\int_\Omega \frac{|u|^{\crits}}{|x|^s}\, dx\right)^{\frac{2}{\crits}}
&\leq & (\mu_{\gamma,s}(\rnp)^{-1}+\epsilon)\int_\Omega\left(|\nabla (\sqrt{\eta}u)|^2-\gamma\frac{\eta u^2}{|x|^2}\right)\, dx\nonumber\\
&& +\Vert \sqrt{1-\eta} u\Vert^2_{\crits,|x|^{-s}}\nonumber\\
&\leq & (\mu_{\gamma,s}(\rnp)^{-1}+\epsilon)\int_\Omega\eta \left(|\nabla u|^2-\gamma\frac{u^2}{|x|^2}\right)\, dx + C\int_\Omega u^2\, dx\nonumber\\
&&+\Vert \sqrt{1-\eta} u\Vert^2_{\crits,|x|^{-s}}\label{eq:pf:3.00}
\end{eqnarray} 
\medskip\noindent {\bf Case 1: $s=0$.} Then $\crits=\crit$ and it follows from Sobolev's inequality that
\begin{eqnarray}\label{eq:pf:4.00}
\Vert \sqrt{1-\eta} u\Vert^2_{\crits,|x|^{-s}}&\leq& K(n,2)^2\int_\Omega |\nabla(\sqrt{1-\eta}u)|^2\, dx\nonumber \\
&\leq& K(n,2)^2\int_\Omega(1-\eta)|\nabla u|^2\, dx+C\int_\Omega u^2\, dx,
\end{eqnarray}
where $K(n,2)$ is the optimal Sobolev constant defined in \eqref{def:K}. Since $s=0$, it follows from \eqref{bnd:mu} that $K(n,2)^2\leq \mu_{\gamma,s}(\rnp)^{-1}$. It then follows from \eqref{eq:pf:4.00} that
\begin{eqnarray}\label{eq:pf:5.00}
\Vert \sqrt{1-\eta} u\Vert^2_{\crits,|x|^{-s}}&\leq&  (\mu_{\gamma,s}(\rnp)^{-1}+\epsilon)\int_\Omega(1-\eta)\left(|\nabla u|^2-\gamma\frac{u^2}{|x|^2}\right)\, dx\nonumber\\
&&+C\int_\Omega u^2\, dx.
\end{eqnarray}
Plugging together \eqref{eq:pf:3.00} and \eqref{eq:pf:5.00} yields \eqref{ineq:sobo:eps:bis} when $s=0$. \\

\medskip\noindent {\bf Case 2: $0<s<2$.} We let $\nu>0$ be a positive number to be fixed later. Since $2<\crits<\crit$, the interpolation inequality yields the existence of $C_\nu>0$ such that
\begin{eqnarray*}
\Vert \sqrt{1-\eta} u\Vert^2_{\crits,|x|^{-s}}&\leq & C \Vert \sqrt{1-\eta} u\Vert^2_{\crits}\\
&\leq& C\left( \nu \Vert \sqrt{1-\eta} u\Vert^2_{\crit}+C_\nu \Vert \sqrt{1-\eta} u\Vert^2_{2}\right)\\
&\leq& C\left( \nu K(n,2)^2 \Vert \nabla(\sqrt{1-\eta} u)\Vert^2_{2}+C_\nu \Vert \sqrt{1-\eta} u\Vert^2_{2}\right).
\end{eqnarray*}
We choose $\nu>0$ such that  $\nu K(n,2)^2 <\mu_{\gamma,s}(\rnp)^{-1}+\epsilon$. Then we get \eqref{eq:pf:5.00} and we conclude \eqref{ineq:sobo:eps:bis} in the case when $2>s>0$ by combining it with \eqref{eq:pf:3.00}.

\medskip\noindent {\bf Case 3: $s=2$.} This is the easiest case, since then 
\[
\Vert \sqrt{1-\eta} u\Vert^2_{\crits,|x|^{-s}}=\int_\Omega \frac{((1-\eta)u)^2}{|x|^2}\, dx \leq C_\delta \int_\Omega u^2\, dx.
\]
This completes the proof of \eqref{ineq:sobo:eps:bis} for all $s\in [0,2]$, and therefore of Proposition \ref{prop:ineq:eps}. \hfill $\Box$.

\medskip\noindent Now we prove the following result, which will be central for the sequel. The proof is standard.

\begin{theorem} \label{tool} Assume that $\gamma<\frac{n^2}{4}$, $0\leq s\leq  2$ and that $\mu_{\gamma, s}(\Omega)<\mu_{\gamma,s}(\rnp)$. Then there are extremals for  $\mu_{\gamma, s}(\Omega)$.  In particular, there exists a minimizer $u$ in $ \dundeux\setminus\{0\}$ that is a positive solution to the equation 
\begin{eqnarray} \label{pos}
\left\{ \begin{array}{llll}
-\Delta u-\gamma \frac{u}{|x|^2}&=&\mu_{\gamma, s}(\Omega)\frac{u^{\crits-1}}{|x|^s} \ \ &\text{\rm on } \Omega\\
\hfill u&>&0 &\text{\rm on }\partial \Omega\\
\hfill u&=&0 &\text{\rm on }\partial \Omega.
\end{array} \right.
\end{eqnarray} 
\end{theorem}
\noindent {\it Proof of Theorem \ref{tool}:} 
Let $(u_i)\in\dundeux\setminus\{0\}$ be a minimizing sequence for $\mu_{\gamma,s}(\Omega)$, that is $J^\Omega_{\gamma,s}(u_i)=\mu_{\gamma,s}(\Omega)+o(1)$ as $i\to +\infty$. Up to multiplying by a constant, we can assume that
\begin{eqnarray}
&&\int_\Omega \frac{|u_i|^{\crits}}{|x|^s}\, dx=1\quad \hbox{for all $i$,}\label{eq:pf:min:1}\\
&&\int_\Omega\left(|\nabla u_i|^2-\gamma\frac{u_i^2}{|x|^2}\right)\, dx=\mu_{\gamma,s}(\Omega)+o(1)\hbox{ as }i\to +\infty.\label{eq:pf:min:2}
\end{eqnarray}
\medskip\noindent We show that $(u_i)_i$ is bounded in $\dundeux$. 
Indeed, \eqref{eq:pf:min:1} yields that 
\begin{equation}\label{bnd:u:l2}
\hbox{$\int_\Omega u_i^2\, dx\leq C<+\infty$ for all $i$.}
\end{equation}
Fix $\epsilon_0>0$ and use Proposition \ref{prop:ineq:eps} and \eqref{bnd:u:l2} to get that
\begin{equation}
\hbox{$\frac{n^2}{4}\int_\Omega \frac{u_i^2}{|x|^2}\, dx\leq (1+\epsilon_0)\int_\Omega |\nabla u_i|^2\, dx+C$ \quad for all $i$. }
\end{equation}
Since $\gamma<\frac{n^2}{4}$, up to taking $\epsilon_0>0$ small enough, this latest inequality combined with \eqref{eq:pf:min:2} yield the boundedness of $(u_i)_i$ in $\dundeux$. 
It follows that there exists $u\in\dundeux$ such that, up to a subsequence, $(u_i)$ goes to $u$ weakly in $\dundeux$ and strongly in $L^2(\Omega)$ as $i\to +\infty$.

We now show that $\int_\Omega \frac{|u|^{\crits}}{|x|^s}\, dx=1$. For that,  define $\theta_i:=u_i-u\in\dundeux$ for all $i$. In particular, $\theta_i$ goes to $0$ weakly in $\dundeux$ and strongly in $L^2(\Omega)$ as $i\to +\infty$. In particular, we have as $i\to +\infty$, 
\begin{equation}\label{eq:pf:min:3}
1=\int_\Omega \frac{|u_i|^{\crits}}{|x|^s}\, dx=\int_\Omega \frac{|u|^{\crits}}{|x|^s}\, dx+\int_\Omega \frac{|\theta_i|^{\crits}}{|x|^s}\, dx+o(1)
\end{equation}
and 
\begin{equation}\label{eq:pf:min:4}
\mu_{\gamma,s}(\Omega)= \int_\Omega\left(|\nabla u|^2-\gamma\frac{u^2}{|x|^2}\right)\, dx+\int_\Omega\left(|\nabla \theta_i|^2-\gamma\frac{\theta_i^2}{|x|^2}\right)\, dx+o(1),
\end{equation}
For $\epsilon>0$, it follows from the definition of $\mu_{\gamma,s}(\Omega)$ and from \eqref{ineq:sobo:eps:bis} that,  as $i\to +\infty$
\begin{equation}
\mu_{\gamma,s}(\Omega)\left(\int_\Omega \frac{|u|^{\crits}}{|x|^s}\, dx\right)^{\frac{2}{\crits}}\leq \int_\Omega\left(|\nabla u|^2-\gamma\frac{u^2}{|x|^2}\right)\, dx\label{eq:pf:min:5}
\end{equation}
and 
\begin{equation}
(\mu_{\gamma,s}(\rnp)-\epsilon)\left(\int_\Omega \frac{|\theta_i|^{\crits}}{|x|^s}\, dx\right)^{\frac{2}{\crits}}\leq \int_\Omega\left(|\nabla \theta_i|^2-\gamma\frac{\theta_i^2}{|x|^2}\right)\, dx+o(1). \label{eq:pf:min:6}
\end{equation}
Summing these two inequalities and using \eqref{eq:pf:min:3} and \eqref{eq:pf:min:4} and passing to the limit, as $i\to +\infty$, yields
\begin{eqnarray*}
\mu_{\gamma,s}(\Omega)\left(1-\left(\int_\Omega \frac{|u|^{\crits}}{|x|^s}\, dx\right)^{\frac{2}{\crits}}\right)&\geq& (\mu_{\gamma,s}(\rnp)-\epsilon)\left(1-\int_\Omega \frac{|u|^{\crits}}{|x|^s}\, dx\right)^{\frac{2}{\crits}}\\
&\geq& (\mu_{\gamma,s}(\rnp)-\epsilon)\left(1-\left(\int_\Omega \frac{|u|^{\crits}}{|x|^s}\, dx\right)^{\frac{2}{\crits}}\right).
\end{eqnarray*}
Since $\mu_{\gamma,s}(\Omega)<\mu_{\gamma,s}(\rnp)$, then by taking $\epsilon>0$ small enough, we finally conclude that $\int_\Omega \frac{|u|^{\crits}}{|x|^s}\, dx=1$. \\
 It remains to show that $u$ is an extremal for $\mu_{\gamma,s}(\Omega)$. For that, note that since $\int_\Omega \frac{|u|^{\crits}}{|x|^s}\, dx=1$, the definition of $\mu_{\gamma,s}(\Omega)$ 
yields $\int_\Omega\left(|\nabla u|^2-\gamma\frac{u^2}{|x|^2}\right)\, dx\geq \mu_{\gamma,s}(\Omega)$. The second term in the right-hand-side of \eqref{eq:pf:min:4} is nonnegative due to \eqref{eq:pf:min:6}. Therefore, we get that $\int_\Omega\left(|\nabla u|^2-\gamma\frac{u^2}{|x|^2}\right)\, dx= \mu_{\gamma,s}(\Omega)$. This proves the claim and ends the proof of Theorem \ref{tool}.\hfill$\Box$

\section{Sub- and super-solutions for the equation $L_\gamma u=a(x) u$}\label{sec:sub:super}
Here and in the sequel, we shall assume that $0\in\partial\Omega$, where $\Omega$ is a smooth domain. In this section, we shall construct basic sub- and super-solutions for the equation $L_\gamma u=a(x) u$, where $a(x)=O(|x|^{\tau-2})$
for some $\tau >0$. 

\medskip\noindent First recall from the introduction that two solutions for $L_\gamma u=0$, with $u=0$ on $\partial \rnp$ are of the form $u_\alpha(x)=x_1|x|^{-\alpha}$, where  $\alpha\in \{\am,\ap\}$ with 
\begin{equation}
\hbox{$\alpha_{-}(\gamma):=\frac{n}{2}-\sqrt{\frac{n^2}{4}-\gamma}$ \quad and \quad $\alpha_{+}(\gamma):=\frac{n}{2}+ \sqrt{\frac{n^2}{4}-\gamma}.$ }
\end{equation}
These solutions will be the building blocks for sub- and super-solutions of more general linear equations involving $L_\gamma$. 
\begin{proposition}\label{prop:sub:super}
Let $\gamma<\frac{n^2}{4}$ and $\alpha\in \{\am,\ap\}$. Let $0<\tau\leq 1$ and $\beta\in\rr$ such  that $\alpha-\tau<\beta<\alpha$ and $\beta\not\in \{\am,\ap\}$. Then, there exist $u_{\alpha,+},u_{\alpha,-}\in C^\infty(\overline{\Omega}\setminus\{0\})$ such that
\begin{equation}\label{ppty:ua}
\left\{\begin{array}{ll}
u_{\alpha,+},u_{\alpha,-}>0 &\hbox{ in }\Omega\cap B_r(0)\\
u_{\alpha,+},u_{\alpha,-}=0 &\hbox{ in }\partial\Omega \cap B_r(0)\\
-\Delta u_{\alpha,+}-\frac{\gamma+O(|x|^\tau)}{|x|^2}u_{\alpha,+}>0&\hbox{ in }\Omega\cap B_r(0)\\
-\Delta u_{\alpha,-}-\frac{\gamma+O(|x|^\tau)}{|x|^2}u_{\alpha,-}<0&\hbox{ in }\Omega\cap B_r(0).
\end{array}
\right.
\end{equation}
Moreover, we have as $x\to 0$, $x\in \Omega$,  that
\begin{equation}\label{asymp:ua:plus}
u_{\alpha,+}(x)=\frac{d(x,\partial\Omega)}{|x|^{\alpha}}(1+O(|x|^{\alpha-\beta}))\hbox{\rm \, and \, }u_{\alpha,-}(x)=\frac{d(x,\partial\Omega)}{|x|^{\alpha}}(1+O(|x|^{\alpha-\beta})).
\end{equation}

\end{proposition}

\medskip\noindent{\it Proof of Proposition \ref{prop:sub:super}}: We first choose an adapted chart to lift the basic solutions from $\rnp$. 
Since $0\in\partial\Omega$ and $\Omega$ is smooth, there exists $\tilde{U},\tilde{V}$ two bounded domains of $\rn$ such that $0\in \tilde{U}$, $0\in \tilde{V}$, and there exists $c\in C^\infty(\tilde{U},\tilde{V})$ a $C^\infty-$diffeomorphism such that $c(0)=0$,
$$c(\tilde{U}\cap\{x_1>0\})= c(\tilde{U})\cap\Omega\hbox{ and }c(\tilde{U}\cap\{x_1=0\})= c(\tilde{U})\cap\partial\Omega.$$
The orientation of $\partial\Omega$ is chosen in such a way  that for any $x'\in\tilde{U}\cap\{x_1=0\}$, 
$$\left\{\partial_1c(0,x'),\partial_2 c(0,x'),\dots,\partial_n c(0,x') \right\}$$
is a direct basis of $\rn$ (canonically oriented). For $x'\in\tilde{U}\cap\{x_1=0\}$, we define $\nu(x')$ as the unique orthonormal inner vector at the tangent space $T_{c(0,x')}\partial\Omega$ (it is chosen such that $\left\{\nu(x'),\partial_2 c(0,x'),\dots,\partial_n c(0,x') \right\}$ is a direct basis of $\rn$). In particular, on $\rnp:=\{x_1>0\}$, $\nu(x'):=(1,0,\dots,0)$.

\medskip\noindent  Here and in the sequel, we write for any $r>0$ 
\begin{equation}\label{def:ball}
\tilde{B}_r:=(-r,r)\times B_r^{(n-1)}(0)
\end{equation}
where $B_r^{(n-1)}(0)$ denotes the ball of center $0$ and radius $r$ in $\rr^{n-1}$. It is standard that there exists $\delta>0$ such that
\begin{equation}\label{def:chart:1}
\begin{array}{cccc}
\varphi: & \tilde{B}_{2\delta}& \to &\rn\\
& (x_1,x')\in \rr\times\rr^{n-1}& \mapsto & c(0,x')+x_1\nu(x')
\end{array}
\end{equation}
is a $C^\infty-$diffeomorphism onto its open image $\varphi(\tilde{B}_{2\delta})$, and
\begin{equation}\label{def:chart:2}
\varphi(\tilde{B}_{2\delta}\cap\{x_1>0\})= \varphi(\tilde{B}_{2\delta})\cap\Omega\hbox{ and }\varphi(\tilde{B}_{2\delta}\cap\{x_1=0\})= \varphi(\tilde{B}_{2\delta})\cap\partial\Omega.
\end{equation}
We also have for all $x'\in B_\delta(0)^{(n-1)}$, 
\begin{equation}\label{def:chart:3}
\nu(x')\hbox{ is the outer orthonormal vector at the tangent space }T_{\varphi(0,x')}\partial\Omega.
\end{equation}
An important remark is that
\begin{equation}\label{def:chart:4}
d(\varphi(x_1,x'),\partial\Omega)=|x_1|\hbox{ for all }(x_1,x')\in \tilde{B}_{2\delta}\hbox{ close to }0.
\end{equation}

\medskip\noindent Consider the metric $g:=\varphi^\star\eucl$ on $\tilde{B}_{2\delta}$, that  is the pull-back of the Euclidean metric $\eucl$ via the diffeomorphism $\varphi$. Following classical notations, we define 
\begin{equation}\label{def:chart:5}
g_{ij}(x):=\left(\partial_i\varphi(x),\partial_j\varphi(x)\right)_{\eucl}\hbox{ for all }x\in \tilde{B}_{2\delta}\hbox{ and }i,j=1,...,n.
\end{equation}

\medskip\noindent Up to a change of coordinates, we can assume that $(\partial_2\varphi(0),...,\partial_n\varphi(0))$ is an orthogonal basis of $T_0\partial\Omega.$ In other words, we then have that
\begin{equation}\label{def:chart:6}
g_{ij}(0)=\delta_{ij}\hbox{ for all }i,j=1,...,n.
\end{equation}

\medskip\noindent We claim that 
\begin{equation}\label{eq:g:ortho}
g_{i1}(x)=\delta_{i1}\hbox{ for all }x\in \tilde{B}_{2\delta}\hbox{ and }i=1,...,n.
\end{equation}
Indeed, for any $x=(x_1,x')\in \tilde{B}_{2\delta}$, we have that $\partial_1\varphi(x)=\nu(x')$, which is a unitary vector. Therefore $g_{11}(x)=1$. For $i\geq 2$, we have
$$g_{1i}(x)=\left(\nu(x'), \partial_i\varphi(0,x')+x_1\partial_i\nu(x')\right)_{\eucl}=\left(\nu(x'), \partial_i\varphi(0,x')\right)_{\eucl}+x_1\partial_i\left(|\nu(x')|^2\right)/2.$$
Since $\nu(x')$ is orthogonal to the tangent space spanned by $(\partial_2\varphi(0,x'),...,\partial_n\varphi(0,x'))$ and $|\nu(x')|=1$, we get that $g_{1i}(x)=0$, which proves (\ref{eq:g:ortho}).

\medskip\noindent Fix now $\alpha\in\rr$ and consider $\Theta\in C^\infty(\tilde{B}_{2\delta})$ such that $\Theta(0)=0$ and which will be constructed later (independently of $\alpha$) with additional needed properties. Fix $\eta\in C^\infty_c(\tilde{B}_{2\delta})$ such that $\eta(x)=1$ for all $x\in \tilde{B}_{\delta}$. Define $u_\alpha\in C^\infty(\overline{\Omega}\setminus\{0\})$ as
\begin{equation}\label{def:ua}
u_\alpha\circ\varphi(x_1,x'):=\eta(x)x_1|x|^{-\alpha}(1+\Theta(x))\hbox{ for all }(x_1,x')\in \tilde{B}_{2\delta}\setminus \{0\}.
\end{equation}
In particular, $u_\alpha(x)>0$ for all $x\in \varphi(\tilde{B}_{2\delta})\cap \Omega$ and $u_\alpha(x)=0$ on $\Omega\setminus  \varphi(\tilde{B}_{2\delta})$. 

\medskip\noindent We claim that with a good choice of $\Theta$, we have that
\begin{equation}\label{eq:ua}
\hbox{$-\Delta u_\alpha=\frac{\alpha(n-\alpha)}{|x|^2}u_\alpha+O\left(\frac{u_\alpha(x)}{|x|}\right)$ \quad as $x\to 0$.}
\end{equation}
Indeed, using the chart $\varphi$, we have that
$$\left(-\Delta u_\alpha\right)\circ\varphi(x_1,x')=-\Delta_g(u_\alpha\circ\varphi)(x_1,x')$$
for all $(x_1,x')\in \tilde{B}_{\delta}\setminus\{0\}$. Here, $-\Delta_g$ is the Laplace operator associated to the metric $g$, that is
$$-\Delta_g:=-g^{ij}\left(\partial_{ij}-\Gamma_{ij}^k\partial_k\right),$$
where
$$\Gamma_{ij}^k:=\frac{1}{2}g^{km}\left(\partial_ig_{jm}+\partial_j g_{im}-\partial_m g_{ij}\right),$$
and $(g^{ij})$ is the inverse of the matrix $(g_{ij})$. Here and in the sequel, we have adopted Einstein's convention of summation. It follows from \eqref{eq:g:ortho} that
\begin{eqnarray}
\left(-\Delta u_\alpha\right)\circ\varphi&=& -\Delta_{\eucl}(u_\alpha\circ\varphi)-\sum_{i,j\geq 2}\left(g^{ij}-\delta^{ij}\right)\partial_{ij}(u_\alpha\circ\varphi)\nonumber\\
&&+ g^{ij}\Gamma_{ij}^1\partial_1(u_\alpha\circ\varphi)+\sum_{k\geq 2}g^{ij}\Gamma_{ij}^k\partial_k(u_\alpha\circ\varphi). \label{eq:ua:2}
\end{eqnarray}
It follows from the definition \eqref{def:ua} that there exists $C>0$ such that for any $i,j,k\geq 2$, we have that
$$|\partial_{ij}(u_\alpha\circ\varphi)(x_1,x')|\leq C|x_1|\cdot|x|^{-\alpha-2}\hbox{ and }|\partial_{k}(u_\alpha\circ\varphi)(x_1,x')|\leq C|x_1|\cdot|x|^{-\alpha-1},$$
for all $(x_1,x')\in \tilde{B}_{\delta}\setminus\{0\}$. It follows from \eqref{def:chart:6} that $g^{ij}-\delta^{ij}=O(|x|)$ as $x\to 0$. Therefore, \eqref{eq:ua:2} yields that as $x\to 0$, 
\begin{eqnarray}\label{eq:ua:3}
\left(-\Delta u_\alpha\right)\circ\varphi&=& -\Delta_{\eucl}(u_\alpha\circ\varphi)+g^{ij}\Gamma_{ij}^1\partial_1(u_\alpha\circ\varphi)+O(x_1|x|^{-\alpha-1})
\end{eqnarray}
 The definition of $g_{ij}$ and the expression of $\varphi(x_1,x')$ then yield that as $x\to 0$, 
\begin{eqnarray*}
g^{ij}\Gamma_{ij}^1&=&-\frac{1}{2}\sum_{i,j\geq 2}g^{ij}\partial_1 g_{ij}\\
&=& -\sum_{i,j\geq 2}g^{ij}(x_1,x')\left((\partial_i\varphi(0,x'),\partial_j\nu(x'))+x_1(\partial_i(x'),\partial_j\nu(x'))\right)\\
&=&-\sum_{i,j\geq 2}g^{ij}(0,x')\left(\partial_i\varphi(0,x'),\partial_j\nu(x')\right)+O(|x_1|)\\
&=&H(x')+O(|x_1|),
\end{eqnarray*}
where $H(x')$ is the mean curvature of the $(n-1)-$manifold $\partial\Omega$ at $\varphi(0,x')$ oriented by the outer normal vector $-\nu(x')$. Using the expression \eqref{def:ua} and using the smoothness of $\Theta$, \eqref{eq:ua:3} yields
\begin{eqnarray*}
\left(-\Delta u_\alpha\right)\circ\varphi&=& \left(-\Delta_{\eucl}(x_1|x|^{-\alpha})\right)\cdot(1+\Theta)+|x|^{-\alpha}\left(H(x')(1+\Theta)-2\partial_1\Theta\right)\nonumber\\
&&+O(x_1|x|^{-\alpha-1}) \qquad \hbox{as $x\to 0$}.
\end{eqnarray*}
We now define
$$\hbox{$\Theta(x_1,x'):=e^{-\frac{1}{2}x_1 H(x')}-1$ for all $x=(x_1,x')\in \tilde{B}_{2\delta}$. }$$
Clearly $\Theta(0)=0$ and $\Theta\in C^\infty(\tilde{B}_{2\delta})$. Noting that 
$$-\Delta_{\eucl}\left(x_1|x|^{-\alpha}\right)=\frac{\alpha(n-\alpha)}{|x|^2}x_1|x|^{-\alpha},$$
we then get that as $x\to 0$, 
\begin{eqnarray}\label{eq:delta:ua:1}
\left(-\Delta u_\alpha\right)\circ\varphi&=& \frac{\alpha(n-\alpha)}{|x|^2}x_1|x|^{-\alpha}\cdot(1+\Theta)+O(x_1|x|^{-\alpha-1})
\end{eqnarray}
With the choice that $g_{ij}(0)=\delta_{ij}$, we have that $(\partial_i\varphi(0))_{i=1,...,n}$ is an orthonormal basis of $\rn$, and therefore $|\varphi(x)|=|x|(1+O(|x|))$ as $x\to 0$. It then follows from \eqref{eq:delta:ua:1} and \eqref{def:ua} that
\begin{equation}\label{est:ua:1}
-\Delta u_\alpha=\frac{\alpha(n-\alpha)}{|x|^2}u_\alpha+O(|x|^{-1}u_\alpha) \qquad \hbox{as $x\to 0$. }
\end{equation}
This proves \eqref{eq:ua}. We now proceed with the construction of the sub- and super-solutions. Let $\alpha\in \{\am,\ap\}$ in such a way that $\alpha(n-\alpha)=\gamma$ and consider $\beta,\lambda\in \rr$ to be chosen later. It follows from \eqref{eq:ua} that
\begin{eqnarray*}
\left(-\Delta-\frac{\gamma+O(|x|^\tau)}{|x|^2}\right)(u_\alpha+\lambda u_\beta)&=& \frac{\lambda(\beta(n-\beta)-\gamma)}{|x|^2}u_\beta\\
&&+\frac{O(|x|^\tau)}{|x|^2}u_\alpha+O(|x|^{-1}u_\alpha)+O(|x|^{\tau-2}u_\beta)\\
&=& \frac{u_\beta}{|x|^2}\big(\lambda(\beta(n-\beta)-\gamma)\\
&& +O(|x|^\tau)+ O(|x|^{\tau+\beta-\alpha})+ O(|x|^{1+\beta-\alpha})\big)
\end{eqnarray*}
as $x\to 0$. Choose $\beta$ such that $\alpha-\tau<\beta<\alpha$ in such a way that $\beta\neq\am$ and $\beta\neq \ap$. In particular, $\beta>\alpha-1$ and $\beta(n-\beta)-\gamma\neq 0$. We then have 
\begin{equation}\label{eq:super:sub:sol}
\left(-\Delta-\frac{\gamma+O(|x|^\tau)}{|x|^2}\right)(u_\alpha+\lambda u_\beta)=\frac{u_\beta}{|x|^2}\left(\lambda(\beta(n-\beta)-\gamma) + O(|x|^{\tau+\beta-\alpha}))\right)
\end{equation}
as $x\to 0$. Choose $\lambda\in\rr$ such that $\lambda(\beta(n-\beta)-\gamma)>0$. Finally, let $u_{\alpha,+}:=u_\alpha+\lambda u_\beta$ and $u_{\alpha,-}:=u_\alpha-\lambda u_\beta$. They clearly verify \eqref{ppty:ua} and \eqref{asymp:ua:plus}, which completes the proof of Proposition \ref{prop:sub:super}.\hfill$\Box$

\section{Regularity and Hopf-type result for the operator $L_\gamma$}\label{sec:hopf}

This section is devoted to the proof of the following key result. 

\begin{theorem}[Optimal regularity and Generalized Hopf's Lemma]\label{th:hopf} Fix $\gamma<\frac{n^2}{4}$ and let $f: \Omega\times \rr\to \rr$ be a Caratheodory function such that
$$|f(x, v)|\leq C|v| \left(1+\frac{|v|^{\crits-2}}{|x|^s}\right)\hbox{\rm  for all }x\in \Omega\hbox{ \rm and }v\in \rr.$$
Let $u\in \dundeux_{loc, 0}$ be a weak solution of 
\begin{equation}\label{regul:eq}
-\Delta u-\frac{\gamma+O(|x|^\tau)}{|x|^2}u=f(x,u)\,\, \hbox{\rm  in }\dundeux_{loc, 0}
\end{equation}
for some $\tau>0$. Then there exists $K\in\rr$ such that
\begin{equation}\label{eq:hopf}
\lim_{x\to 0}\frac{u(x)}{d(x,\partial\Omega)|x|^{-\am}}=K.
\end{equation}
Moreover, if $u\geq 0$ and $u\not\equiv 0$, we have that $K>0$.
\end{theorem}
\noindent As mentioned in the introduction, this can be viewed as a generalization of Hopf's Lemma in the following sense: when $\gamma=0$ (and then $\am=0$), the classical Nash-Moser  regularity scheme yields $u\in C^1_{loc}$, and when $u\geq 0$, $u\not\equiv 0$, Hopf's comparison principle yields $\partial_\nu u(0)<0$, which is a reformulation of \eqref{eq:hopf} when $\am=0$.

\smallskip\noindent The remaining of this section is devoted to the proof of Theorem \ref{th:hopf}. We shall need the following two lemmae, which will be used frequently throughout the paper.

\begin{lem}\label{lem:funda}{\rm (Rigidity of solutions)} Let $u\in C^2(\overline{\rnp}\setminus \{0\})$ be a nonnegative function such that
\begin{equation}\label{eq:entire}
-\Delta u-\frac{\gamma}{|x|^2}u=0\hbox{ in }\rnp\; ;\; u=0\hbox{ on }\partial \rnp.
\end{equation}
Suppose there exists $\alpha\in \{\am,\ap\}$ such that $u(x)\leq C|x|^{1-\alpha}$, then there exists $\lambda\geq 0$ such that
$$u(x)=\lambda x_1|x|^{-\alpha}\hbox{ for all }x\in \rnp.$$
\end{lem}
\noindent We note that this lemma is only a first step in proving rigidity for solutions of $L_\gamma u=0$ on $\rnp$. Indeed, the pointwise assumption above  is not necessary as it will be removed in Proposition \ref{prop:funda:2}, which will be a consequence of the classification Theorem \ref{th:classif}. 
\par
\smallskip\noindent{\it Proof of Lemma \ref{lem:funda}:} 
We first assume that $\alpha:=\am$ and prove that
\begin{equation}\label{ineq:1}
\hbox{either }u\equiv 0\hbox{ or }\liminf_{|x|\to +\infty}\frac{u(x)}{x_1|x|^{-\am}}>0.
\end{equation}
Indeed suppose that the second situation does not hold, that is there exists $(x_i)_i\in\rnp$ such that 
$$\lim_{i\to +\infty}|x_i|=+\infty\hbox{ and }\lim_{i\to +\infty}\frac{u(x_i)}{(x_i)_1|x_i|^{-\am}}=0.$$
Define $r_i:=|x_i|$, $\theta_i:=\frac{x_i}{|x_i|}$ and $u_i(x):=r_i^{\am-1}u(r_i x)$ for all $i$ and all $x\in\rnp$. It follows from the hypothesis of the theorem that for all $i$, 
$$-\Delta u_i-\frac{\gamma}{|x|^2}u_i=0\hbox{ in }\rnp\; ;\; u_i=0\hbox{ on }\partial \rnp,$$
and $0\leq u_i(x)\leq C|x|^{1-\am}$ for all $x\in\rnp$. It follows from elliptic theory that there exist $\hat{u}\in C^2(\rnp)$ such that $u_i\to \hat{u}$ in $C^2_{loc}(\overline{\rnp}\setminus\{0\})$. In particular, we have that
$$-\Delta \hat{u}-\frac{\gamma}{|x|^2}\hat{u}=0\hbox{ in }\rnp \; ;\; \hat{u}\geq 0\hbox{ in }\rnp,\; ;\; \hat{u}=0\hbox{ on }\partial \rnp.$$
Let $\theta:=\lim_{i\to +\infty}\theta_i$. It follows from the convergence that $\hat{u}(\theta)=0$ if $\theta\in\rnp$, and $\partial_1\hat{u}(\theta)=0$ if $\theta\in\partial\rnp$. Hopf's maximum principle yields that $\hat{u}\equiv 0$. In particular, we get that
\begin{equation}\nonumber
\lim_{i\to +\infty}\sup_{x\in \partial B_{r_i}(0)}\frac{u(x)}{x_1|x|^{-\am}}=0.
\end{equation}
For $\epsilon>0$, there exists $i_0$ such that $u(x)\leq \epsilon x_1|x|^{-\am}$ for all $x\in \partial B_{r_i}(0)$ and $i\geq i_0$. Since $u-\epsilon x_1|x|^{-\am}$ is locally in $D^{1,2}$ and $(-\Delta-\frac{\gamma}{|x|^2})(u-\epsilon x_1|x|^{-\am})=0$, it follows from the maximum principle (and coercivity) that $u(x)\leq \epsilon x_1|x|^{-\am}$ for $x\in B_{r_i}(0)$ with $i\geq i_0$. Letting $i\to +\infty$ yields 
$u(x)\leq \epsilon x_1|x|^{-\am}$ for all $x\in\rnp$ and all $\epsilon>0$. Letting $\epsilon\to 0$ yields $u\leq 0$, and therefore $u\equiv 0$. This proves the claim \eqref{ineq:1}.

We now assume that $u\not\equiv 0$. It then follows from \eqref{ineq:1} that there exists $\epsilon_0>0$ and $R_0>0$ such that $u(x)\geq \epsilon_0 x_1|x|^{-\am}$ for all $|x|\geq B_{R_0}(0)$. Applying again the maximum principle on $B_{R_0}(0)$, we get that $u(x)\geq \epsilon_0 x_1|x|^{-\am}$ for all $x\in\rnp$. We have so far proved that
\begin{equation}\label{eq:dicho}
u\equiv 0\hbox{ or there exists }\epsilon_0>0\hbox{ such that }u(x)\geq \epsilon_0 x_1|x|^{-\am}\hbox{ for all }x\in\rnp.
\end{equation}
Let now $\lambda:=\max\{k\geq 0\hbox{ such that }u(x)\geq k x_1|x|^{-\am}\hbox{ for all }x\in\rnp\}$. Then $\bar{u}(x):=u(x)-\lambda x_1|x|^{-\am}\geq 0$ satisfies \eqref{eq:entire}. It then follows from \eqref{eq:dicho} that $\bar{u}\equiv 0$ or $\bar{u}(x)\geq \epsilon_0 x_1|x|^{-\am}$ for all $x\in\rnp$ for some $\epsilon_0>0$. This second case cannot happens since it would imply that $u\geq (\lambda+\epsilon_0)x_1|x|^{-\am}$, which is a contradiction. Therefore $\bar{u}\equiv 0$ and the Proposition is proved when $\alpha=\am$.

To finish, it remains to consider the case where $\alpha=\ap$. Here we define $\tilde{u}(x):=|x|^{2-n}u(x/|x|^2)$ to be the Kelvin transform of $u$. The function $\tilde{u}$ then satisfies \eqref{eq:entire} with $\am$. It then follows from the first part of this proof that $\tilde{u}=\lambda x_1|x|^{-\am}$. Coming back to the initial function $u$ yields $u=\lambda x_1|x|^{-\ap}$. This completes the proof of Lemma \ref{lem:funda}.\hfill$\Box$\\

\begin{lem}\label{lem:deriv:1} Let $f: \Omega\times \rr\to \rr$ be as in the statement of Theorem \ref{th:hopf}, and consider $u\in \dundeux_{loc, 0}$ such that \eqref{regul:eq} holds. Assume there exists $C>0$ and $\alpha\in\{\ap,\am\}$ such that
\begin{equation}\label{bnd:u:lem}
|u(x)|\leq C |x|^{1-\alpha}\hbox{ for }x\to 0, \, x\in\Omega.
\end{equation}
If $f\not\equiv 0$, we assume that $\alpha=\am$.
\begin{enumerate}
\item Then, there exists $C_1>0$ such that
\begin{equation}\label{ineq:der}
|\nabla u(x)|\leq C_1 |x|^{-\alpha}\hbox{ as }x\to 0, \, x\in\Omega.
\end{equation}
\item If $\lim_{x\to 0}|x|^{\alpha-1}u(x)=0$, then $\lim_{x\to 0}|x|^{\alpha}|\nabla u(x)|=0$. Moreover, if $u>0$, then there exists $l\geq 0$ such that
\begin{equation}\label{lim:u:uni}
\lim_{x\to 0}\frac{|x|^\alpha u(x)}{d(x,\partial\Omega)}=l\hbox{ and }\lim_{x\to 0,\, x\in\partial\Omega}|x|^\alpha |\nabla u(x)|=l
\end{equation}
\end{enumerate}
\end{lem}
\noindent{\it Proof of Lemma \ref{lem:deriv:1}:} Assume that \eqref{bnd:u:lem} holds. As a first remark, we claim that we can assume that for some $\tau>0$, 
\begin{equation}\label{regul:eq:BIS}
-\Delta u-\frac{\gamma+O(|x|^\tau)}{|x|^2}u=0\,\, \hbox{\rm  in }\dundeux_{loc, 0}.
\end{equation}
Indeed, this is clear if $f\equiv 0$. If $f\not\equiv 0$, since $\alpha=\am$, we have as $x\to 0$, 
\begin{eqnarray*}
|f(x,u)|&\leq &C |u| \left(1+|x|^{-s}|x|^{-(\crits-2)(\am-1)}\right)\\
&&\leq C\frac{|u|}{|x|^2}\left(|x|^2+|x|^{(\crits-2)(\frac{n}{2}-\am)}=O\left(|x|^{\tau'}\frac{u}{|x|^2}\right)\leq C\right)
\end{eqnarray*}
for some $\tau'>0$. Plugging this inequality into \eqref{regul:eq} and replacing $\tau$ by $\min\{\tau, \tau'\} $ yields \eqref{regul:eq:BIS}. 

\medskip\noindent In the sequel, we shall write 
$\omega(x):=\frac{|x|^\alpha u(x)}{d(x,\partial\Omega)}$
for $x\in\Omega$. Let $(x_i)_i\in\Omega$ be such that 
\begin{equation}\label{lim:l}
\lim_{i\to +\infty}x_i=0
\end{equation}
Choose a chart $\varphi$ as in \eqref{def:chart:1} such that $d\varphi_0=Id_{\rn}$. For any $i$, define $X_i\in\rnp$ such that $x_i=\varphi(X_i)$, $r_i:=|X_i|$ and $\theta_i:=\frac{X_i}{|X_i|}$. In particular, $\lim_{i\to +\infty}r_i=0$ and $|\theta_i|$=1 for all $i$. Set 
$$\tilde{u}_i(x):=r_i^{\alpha-1}u(\varphi(r_i x))\hbox{ for all }i\hbox{ and }x\in B_R(0)\cap\rnp\, ;\, x\neq 0.$$
Equation \eqref{regul:eq:BIS} then rewrites
\begin{equation}\label{eq:ui}
\left\{\begin{array}{ll}
-\Delta_{g_i} \tilde{u}_i-\frac{\gamma+o(1)}{|x|^2}\tilde{u}_i=0&\hbox{ in }B_R(0)\cap\rnp\\
\hfill \tilde{u}_i=0&\hbox{ in }B_R(0)\cap\partial\rnp,
\end{array}\right.
\end{equation}
where $g_i(x):=(\varphi^\star\eucl)(r_ix)$ is a metric that goes to $\eucl$ on every compact subset of $\rn$ as $i\to \infty$. Here, $o(1)\to 0$ in $C^0_{loc}(\overline{\rnp}\setminus \{0\})$. It follows from \eqref{bnd:u:lem} and \eqref{lim:l} that
\begin{equation}\label{lim:l:2}
|u_i(x)|\leq C |x|^{1-\alpha}\hbox{ for all }i\hbox{ and all }x\in B_R(0)\cap\rnp, 
\end{equation}
It follows from elliptic theory, that there exists $\tilde{u}\in C^2(\overline{\rnp}\setminus\{0\})$ such that $\tilde{u}_i\to \tilde{u}$ in $C^1_{loc}(\overline{\rnp}\setminus\{0\})$. By letting $\theta:=\lim_{i\to +\infty}\theta_i$ ($|\theta|=1$), we then have that for any $j=1,...,n$, $\partial_j\tilde{u}_i(\theta_i)\to \partial_j \tilde{u}(\theta)$ as $i\to +\infty$, which rewrites
\begin{equation}\label{lim:der}
\lim_{i\to +\infty}|x_i|^\alpha\partial_j u (x_i)=\partial_j \tilde{u}(\theta)\hbox{ for all } j=1,...,n.
\end{equation}

\smallskip\noindent We now prove \eqref{ineq:der}. For that, we argue by contradiction and  assume that there exists a sequence $(x_i)_i\in\Omega$ that goes to $0$ as $i\to +\infty$ and such that $|x_i|^\alpha|\nabla u(x_i)|\to +\infty$ as $i\to +\infty$. It then follows from \eqref{lim:der}  that $|x_i|^\alpha|\nabla u(x_i)|=O(1)$ as $i\to +\infty$. A contradiction with our assumption, which proves \eqref{ineq:der}. The case when $|x|^\alpha u(x)\to 0$ as $x\to 0$ goes similarly.


\medskip\noindent Now we consider the case when $u>0$, which implies 
that $\tilde{u}_i\geq 0$ and $\tilde{u}\geq 0$. We let $l\in [0,+\infty]$ and $(x_i)_i\in\Omega$ be such that 
\begin{equation}\label{lim:l:2:bis}
\lim_{i\to +\infty}x_i=0\hbox{ and }\lim_{i\to +\infty}\omega(x_i)=l.
\end{equation}
We claim that 
\begin{equation}\label{eq:lim:u}
0\leq l<+\infty\hbox{ and }\lim_{x\to 0}\omega(x)=l\in [0,+\infty).
\end{equation}
Indeed, using the notations above, we get that
$$\lim_{i\to +\infty}\frac{\tilde{u}_i(\theta_i)}{(\theta_i)_1}=l.$$
The convergence of $\tilde{u}_i$ in $C^1_{loc}(\overline{\rnp}\setminus\{0\})$ then yields $l<+\infty$. Passing to the limit as $i\to +\infty$ in \eqref{eq:ui}, we get
$$\left\{\begin{array}{ll}
-\Delta_{\eucl} \tilde{u}-\frac{\gamma}{|x|^2}\tilde{u}=0&\hbox{ in }\rnp\\
\hfill \tilde{u}\geq 0&\hbox{ in }\rnp\\
\hfill \tilde{u}=0&\hbox{ in }\partial\rnp.
\end{array}\right.$$
The limit \eqref{lim:l:2:bis} can be rewriten as $\tilde{u}(\theta)=l\theta_1$ if $\theta\in\rnp$ and $\partial_1\tilde{u}(\theta)=l$ if $\theta\in\partial\rnp$. The rigidity Lemma \ref{lem:funda} then yields
$$\tilde{u}(x)=l x_1|x|^{-\alpha}\hbox{ for all }x\in\rnp.$$
In particular, since the differential of $\varphi$ at $0$ is the identity map, it follows from the convergence of $\tilde{u}_i$ to $\tilde{u}$ locally in $C^1$ that
\begin{equation}\label{est:u:bnd:sup}
\lim_{i\to +\infty}\sup_{x\in\Omega\cap \partial B_{r_i}(0)}\frac{u(x)}{d(x,\partial\Omega)|x|^{-\alpha}}=\sup_{x\in \rnp\cap \partial B_{1}(0)}\frac{\tilde{u}(x)}{x_1|x|^{-\alpha}}=l
\end{equation}
and
\begin{equation}\label{est:u:bnd:inf}
\lim_{i\to +\infty}\inf_{x\in\Omega\cap \partial B_{r_i}(0)}\frac{u(x)}{d(x,\partial\Omega)|x|^{-\alpha}}=\inf_{x\in\rnp\cap \partial B_{1}(0)}\frac{\tilde{u}(x)}{x_1|x|^{-\alpha}}=l.
\end{equation}
We distinguish two cases:

\smallskip\noindent{\bf Case 1: $\alpha=\ap$.} Let $W$ and $u_{\ap, -}^{(d)}$ be as in Proposition \ref{prop:super:dirichlet}, and  fix $\epsilon>0$. It follows from \eqref{est:u:bnd:inf} that there exists $i_0$ such that for $i\geq i_0$, we have that
\begin{equation}\nonumber
u(x)\geq (l-\epsilon)u_{\ap, -}^{(d)}(x)\hbox{ for all }x\in W\cap \partial B_{r_i}(0).
\end{equation}
Since $(-\Delta-(\gamma+O(|x|^\tau))|x|^{-2})(u- (l-\epsilon)u_{\ap, -}^{(d)})\geq 0$ in $W\setminus B_{r_i}(0)$ and since $u_{\ap,-}$ vanishes on $\partial W\setminus\{0\}$, it follows from the comparison principle that
\begin{equation}\nonumber
u(x)\geq (l-\epsilon)u_{\ap, -}^{(d)}(x)\hbox{ for all }x\in W\setminus \partial B_{r_i}(0).
\end{equation}
Letting $i\to +\infty$ yields 
\begin{equation}\nonumber
u(x)\geq (l-\epsilon)u_{\ap, -}^{(d)}(x)\hbox{ for all }x\in W\setminus \{0\}.
\end{equation}
It then follows from this inequality and the asymptotics for $u_{\ap,-}^{(d)}$ that
$$\liminf_{x\to 0}\omega(x)\geq l.$$
Note that this is valid for any $l\in\rr$ satisfying \eqref{lim:l:2:bis}. By taking $l:=\limsup_{x\to 0}\omega(x)$, we then get that $\lim_{x\to 0}\omega(x)=l$. 

\smallskip\noindent{\bf Case 2: $\alpha=\am$.} Consider the super- and sub-solutions $u_{\am,+}, u_{\am,-}$ constructed in Proposition \ref{prop:sub:super}.  It follows from \eqref{est:u:bnd:sup} and \eqref{est:u:bnd:inf} that for $\epsilon>0$, there exists $i_0$ such that for $i\geq i_0$, we have  
$$(l-\epsilon) u_{\am,-}(x)\leq u(x)\leq (l+\epsilon)u_{\am,+}(x)\hbox{ for all }x\in \Omega\cap \partial B_{r_i}(0).$$
Since the operator $-\Delta-(\gamma+O(|x|^\tau))|x|^{-2}$ is coercive on $\Omega\cap B_{r_i}(0)$ and that the functions we consider are in $D^{1,2}$ (i.e., they are variational), it follows from the maximum principle that
$$(l-\epsilon) u_{\am,-}(x)\leq u(x)\leq (l+\epsilon)u_{\am,+}(x)\hbox{ for all }x\in \Omega\cap B_{r_i}(0).$$
Using the asymptotics \eqref{asymp:ua:plus} of the sub- and super-solution, we get that
$$(l-\epsilon) \leq \liminf_{x\to 0}\frac{u(x)}{d(x,\partial\Omega)|x|^{-\am}}\leq \limsup_{x\to 0}\frac{u(x)}{d(x,\partial\Omega)|x|^{-\am}}\leq (l+\epsilon).$$
Letting $\epsilon\to 0$ yields $\lim_{x\to 0}\omega(x)=l\geq0$. This ends Case 2 and completes the proof of \eqref{eq:lim:u}.

\smallskip\noindent The case $u>0$ is a consequence of \eqref{eq:lim:u} and \eqref{lim:der} (note that for the second limit, $x_i\in\partial\Omega$ rewrites  $\theta_i\in\partial\rnp$ and therefore $(\theta_i)_1=0$).
This ends the proof of Lemma \ref{lem:deriv:1}.\hfill$\Box$

\medskip\noindent The following lemma is essentially Theorem \ref{th:hopf} in the case of linear equations of the form $L_\gamma u=a(x)u$. 

\begin{lem}\label{Hopf.linear} Assume that  $u\in \dundeux_{loc, 0}$ is a weak solution of 
\begin{equation}\label{hop.lin}
\left\{\begin{array}{ll}
-\Delta u-\frac{\gamma+O(|x|^\tau)}{|x|^2}u=0 &\hbox{\rm in }\dundeux_{loc, 0}\\
\hfill u=0  & \hbox{\rm on }B_{2\delta}(0)\cap \partial \Omega,
\end{array}\right.
\end{equation}
for some $\tau>0$. Then, there exists $\ell\in \rr$ such that
$$\lim_{x\to 0}\frac{u(x)}{d(x,\partial\Omega)|x|^{-\am}}=\ell.$$
\end{lem}
\smallskip\noindent{\it Proof of Lemma \ref{Hopf.linear}: } First, we assume that  $u\in \dundeux_{loc, 0}$, satisfies (\ref{hop.lin}) and $u>0$ on $B_\delta(0)\cap\Omega$. 
\smallskip\noindent We claim that there exists $C_0>0$ such that 
\begin{equation}\label{est:u:up:low}
\frac{1}{C_0} \frac{d(x,\partial\Omega)}{|x|^{\am}}\leq u(x)\leq C_0 \frac{d(x,\partial\Omega)}{|x|^{\am}}\hbox{ for all }x\in \Omega\cap B_\delta(0).
\end{equation}
Indeed, since $u$ is smooth outside $0$, it follows from Hopf's Maximum principle that there exists $C_1,C_2>0$ such that
\begin{equation}\label{est:u}
C_1 d(x,\partial\Omega)\leq u(x)\leq C_2d(x,\partial\Omega)\hbox{ for all }x\in \Omega\cap\partial B_\delta(0).
\end{equation}
Let $u_{\am,+}$ be the super-solution constructed in Proposition \ref{prop:sub:super}. It follows from \eqref{est:u} and the asymptotics \eqref{asymp:ua:plus} of $u_{\am,+}$ that there exists $C_3>0$ such that
\begin{equation}\nonumber
u(x)\leq C_3 u_{\am,+}(x)\hbox{ for all }x\in \partial (B_\delta(0)\cap\Omega).
\end{equation}
Since $u$ is a solution and $u_{\am,+}$ is a supersolution, both being in $\dundeux_{loc,0}$, it follows from the maximum principle (by choosing $\delta>0$ small enough so that $-\Delta-(\gamma+O(|x|^\tau))|x|^{-2}$ is coercive on $B_\delta(0)\cap\Omega$) that $u(x)\leq C_3 u_{\am,+}(x)$ for all $x\in B_\delta(0)\cap\Omega$. In particular, it follows from the asymptotics \eqref{asymp:ua:plus} of $u_{\am,+}$ that there exists $C_4>0$ such that
\begin{equation}\nonumber
u(x)\leq C_4 d(x,\partial\Omega)|x|^{-\am}\hbox{ for all }x\in \Omega\cap B_\delta(0).
\end{equation}
Arguing similarly with the lower-bound in \eqref{est:u} and the subsolution $u_{\am,-}$, we get the existence of $C_0>0$ such that \eqref{est:u:up:low} holds. 

\smallskip\noindent Since $u>0$, \eqref{est:u:up:low} coupled with Lemma \ref{lem:deriv:1} yields Lemma \ref{Hopf.linear} for $u>0$. 
 

\medskip\noindent Now we deal with the case when $u$ is a sign-changing solution for (\ref{hop.lin}). We then define $u_1,u_2: B_\delta(0)\cap\Omega\to \rr$ be such that
$$\left\{\begin{array}{ll}
-\Delta u_1-\frac{\gamma+O(|x|^\tau)}{|x|^2}u_1=0 & \hbox{ in }B_\delta(0)\cap\Omega\\
\hfill u_1(x)=\max\{u(x), 0\} & \hbox{ on }\partial (B_\delta(0)\cap\Omega).
\end{array}\right.$$
$$\left\{\begin{array}{ll}
-\Delta u_2-\frac{\gamma+O(|x|^\tau)}{|x|^2}u_2=0 & \hbox{ in }B_\delta(0)\cap\Omega\\
\hfill u_2(x)=\max\{-u(x), 0\} & \hbox{ on }\partial (B_\delta(0)\cap\Omega).
\end{array}\right.$$
The existence of such solutions is ensured by choosing $\delta>0$ small enough so that the operator $-\Delta-(\gamma+O(|x|^\tau))|x|^{-2}$ is coercive on $B_\delta(0)\cap\Omega$. In particular, $u_1,u_2\in D^{1,2}(\Omega)_{loc,0}$, $u_1,u_2\geq 0$ and $u=u_1-u_2$. It follows from the maximum principle that for all $i$, either $u_i\equiv 0$ or $u_i>0$. The first part of the proof yields that there exists $l_1,l_2\geq  0$ such that
$$\lim_{x\to 0}\frac{u_1(x)}{d(x,\partial\Omega)|x|^{-\am}}=l_1\hbox{ and }\lim_{x\to 0}\frac{u_2(x)}{d(x,\partial\Omega)|x|^{-\am}}=l_2.$$
Therefore, we get that
$$\lim_{x\to 0}\frac{u(x)}{d(x,\partial\Omega)|x|^{-\am}}=l_1-l_2\in\rr.$$
This completes the proof of Lemma  \ref{Hopf.linear}.\hfill$\Box$\\

\noindent{\it Proof of Theorem \ref{th:hopf}:}
We let here $u\in \dundeux_{loc,0}$ be a solution to \eqref{regul:eq}, that is
\begin{equation}\label{regul:eq:2}
-\Delta u-\frac{\gamma+O(|x|^\tau)}{|x|^2}u=f(x,u)\hbox{ weakly in }\dundeux_{loc, 0}
\end{equation}
for some $\tau>0$. We shall first use the classical DeGiorgi-Nash-Moser iterative scheme (see Gilbarg-Trudinger \cite{gt}, and Hebey \cite{hebeybook1} for expositions in book form).  We skip most of the computations and refer to Ghoussoub-Robert (Proposition A.1 of \cite{gr2}) for the details. The proof goes through four steps.

\medskip\noindent{\bf Step 1:} Let $\beta\geq 1$ be such that $\frac{4\beta}{(\beta+1)^2}>\frac{4}{n^2}\gamma$. Assume that $u\in L^{\beta+1}(\Omega\cap B_{\delta_0}(0))$. We claim that
 \begin{equation}\label{eq:iter:4}
 u\in L^{\frac{n}{n-2}(\beta+1)}(\Omega\cap B_{\delta_0}(0)).
 \end{equation}
Indeed,  fix $\beta\geq 1$, $L>0$, and define
 \begin{equation}\label{G:1}
 G_L(t)=\left\{\begin{array}{ll}
|t|^{\beta-1}t& \hbox{ if }|t|\leq L\\
\beta L^{\beta-1}(t-L)+L^\beta & \hbox{ if }t\geq L\\
\beta L^{\beta-1}(t+L)-L^\beta & \hbox{ if }t\leq -L
\end{array}\right.
\end{equation}
and
 \begin{equation}\label{G:2}
H_L(t)=\left\{\begin{array}{ll}
|t|^{\frac{\beta-1}{2}}t& \hbox{ if }|t|\leq L\\
\frac{\beta+1}{2}L^{\frac{\beta-1}{2}}(t-L)+L^{\frac{\beta+1}{2}} & \hbox{
if }t\geq L\\
\frac{\beta+1}{2}L^{\frac{\beta-1}{2}}(t+L)-L^{\frac{\beta+1}{2}} & \hbox{
if }t\leq -L
\end{array}\right.
\end{equation}
As easily checked,
\begin{equation}\label{G:3}
0\leq t G_L(t)\leq H_L(t)^2\hbox{ and
}G_L'(t)=\frac{4\beta}{(\beta+1)^2}(H_L'(t))^2
\end{equation}
for all $t\in\rr$ and all $L>0$. We fix $\delta>0$ small that will be chosen later. We let $\eta\in C^\infty_c(\rn)$ be such that $\eta(x)=1$ for $x\in B_{\delta/2}(0)$ and $\eta(x)=0$ for $x\in \rn\setminus B_{\delta}(0)$. Multiplying equation \eqref{regul:eq:2} with $\eta^2 G_L(u)\in \dundeux$, we get that
\begin{eqnarray}
\int_\Omega (\nabla u,\nabla (\eta^2 G_L(u)))\, dx&-&\int_\Omega\frac{\gamma+O(|x|^\tau)}{|x|^2}\eta^2 u G_L(u)\, dx\nonumber\\
&=&\int_\Omega f(x,u)\eta^2 G_L(u)\, dx.\label{eq:iter:1}
\end{eqnarray}
Integrating by parts, and using formulae \eqref{G:1} to \eqref{G:3} above (see \cite{gr2} for details) yield
\begin{eqnarray}\nonumber
\int_\Omega (\nabla u,\nabla (\eta^2 G_L(u)))\, dx&=&\frac{4\beta}{(\beta+1)^2}\int_\Omega\left(\left|\nabla (\eta H_L(u))\right|^2-\eta(-\Delta)\eta H_L(u)^2\right)\, dx\\
&&+\int_\Omega -\Delta(\eta^2) J_L(u)\, dx
\end{eqnarray}
where $J_L(t):=\int_0^t G_L(\tau)\, d\tau$. This identity and \eqref{eq:iter:1} yield
\begin{eqnarray}
\frac{4\beta}{(\beta+1)^2}\int_\Omega\left|\nabla (\eta H_L(u))\right|^2\, dx&-&\int_\Omega\frac{\gamma+O(|x|^\tau)}{|x|^2}\eta^2 uG_L(u)\, dx\nonumber\\
& \leq& \int_\Omega |-\Delta(\eta^2)|\cdot | J_L(u)|\, dx\nonumber \\
&&
+ C(\beta,\delta)\int_{\Omega\cap B_\delta(0)}|H_L(u)|^2\, dx\nonumber\\
&& + C\int_\Omega\frac{|u|^{\crits-2}}{|x|^s}(\eta H_L(u))^2\, dx.\label{eq:iter:2}
\end{eqnarray}
H\"older's inequality and the Sobolev inequality \eqref{ineq:sobo:s:g} yield
\begin{eqnarray*}
&&\int_\Omega\frac{|u|^{\crits-2}}{|x|^s}(\eta H_L(u))^2\, dx\\
&&\leq \left(\int_{\Omega\cap B_\delta(0)}\frac{|u|^{\crits}}{|x|^s}\, dx\right)^{\frac{\crits-2}{\crits}}\left(\int_\Omega \frac{|\eta H_L(u)|^{\crits}}{|x|^s}\, dx\right)^{\frac{2}{\crits}}\\
&&\leq \left(\int_{\Omega\cap B_\delta(0)}\frac{|u|^{\crits}}{|x|^s}\, dx\right)^{\frac{\crits-2}{\crits}}\cdot \frac{1}{\mu_{0,s}(\Omega)}\int_\Omega |\nabla (\eta H_L(u))|^2\, dx.
\end{eqnarray*}
Plugging this estimate into \eqref{eq:iter:2} and noting $\gamma_+:=\max\{\gamma,0\}$ yields
\begin{eqnarray*}
\frac{4\beta}{(\beta+1)^2}\int_\Omega\left|\nabla (\eta H_L(u))\right|^2\, dx&-&\left(\gamma_++C\delta^\tau\right)\int_\Omega\frac{(\eta H_L(u))^2}{|x|^2}\, dx\\
&&\leq  C(\beta,\delta)\int_{\Omega\cap B_\delta(0)}\left(|H_L(u)|^2+| J_L(u)|\right)\, dx\\
&&+ \alpha(\delta)\int_\Omega |\nabla (\eta H_L(u))|^2\, dx,
\end{eqnarray*}
where
$$\alpha(\delta):=C\left(\int_{\Omega\cap B_\delta(0)}\frac{|u|^{\crits}}{|x|^s}\, dx\right)^{\frac{\crits-2}{\crits}}\cdot \frac{1}{\mu_{0,s}(\Omega)},$$
so that 
$$\lim_{\delta\to 0}\alpha(\delta)=0.$$
It follows from \eqref{ineq:hardy:opt:1} that 
\begin{equation}\nonumber
\frac{n^2}{4}\int_\Omega\frac{(\eta H_L(u))^2}{|x|^2}\, dx\leq (1+\epsilon(\delta))\int_\Omega |\nabla (\eta H_L(u))|^2\, dx,
\end{equation}
where $\lim_{\delta\to 0}\epsilon(\delta)=0$. Therefore, we get that
\begin{eqnarray*}
&&\left(\frac{4\beta}{(\beta+1)^2}-\alpha(\delta)-\left(\gamma_++C\delta^\tau\right)\frac{4}{n^2}(1+\epsilon(\delta))\right)\int_\Omega\left|\nabla (\eta H_L(u))\right|^2\, dx\\
&&\leq  C(\beta,\delta)\int_{\Omega\cap B_\delta(0)}\left(|H_L(u)|^2+| J_L(u)|\right)\, dx\leq C(\beta,\delta)\int_{B_{\delta}(0)\cap \Omega}|u|^{\beta+1}\, dx.
\end{eqnarray*}
Let $\delta\in (0,\delta_0)$ be such that 
$$\frac{4\beta}{(\beta+1)^2}-\alpha(\delta)-\left(\gamma_++C\delta^\tau\right)\frac{4}{n^2}(1+\epsilon(\delta))>0.$$
This is possible since $\frac{4\beta}{(\beta+1)^2}> \frac{4}{n^2}\gamma$. Using Sobolev's embedding, we then get that
\begin{eqnarray*}
\left(\int_{B_{\delta/2}(0)\cap\Omega}|H_L(u)|^{\crit}\, dx\right)^{\frac{2}{\crit}}&\leq& \left(\int_{\rn}|\eta H_L(u)|^{\crit}\, dx\right)^{\frac{2}{\crit}}\\
&\leq& \mu_{0,0}(\Omega)^{-1}\int_\Omega\left|\nabla (\eta H_L(u))\right|^2\, dx\\
&\leq & C(\beta,\delta,\gamma)\int_{B_{\delta}(0)\cap \Omega}|u|^{\beta+1}\, dx.
\end{eqnarray*}
Since $u\in L^{\beta+1}(B_{\delta_0}(0)\cap \Omega)$, let $L\to +\infty$ and use Fatou's Lemma  to obtain that $u\in L^{\frac{\crit}{2}(\beta+1)}(B_{\delta/2}(0)\cap\Omega)$. The standard iterative scheme then yields that $u\in C^1(\overline{\Omega}\cap B_{\delta_0}(0)\setminus\{0\})$. Therefore $u\in L^{\frac{\crit}{2}(\beta+1)}(B_{\delta_0}(0)\cap\Omega)$, which proves  claim \eqref{eq:iter:4}.

\medskip\noindent{\bf Step 2:} We now show that
\begin{eqnarray}
&&\hbox{if }\gamma\leq 0,\hbox{ then }u\in L^p(\Omega\cap B_{\delta}(0))\hbox{ for all }p\geq 1,\\
&&\hbox{if }\gamma>0,\hbox{ then }u\in L^p(\Omega\cap B_{\delta}(0))\hbox{ for all }p\in \left(1,\frac{n}{n-2}\frac{n}{\am}\right).\label{integ:max}
\end{eqnarray}
The case $\gamma\leq 0$ is standard, so we only consider the case where $\gamma>0$. Fix $p\geq 2$ and set $\beta:=p-1$. we have 
$$\frac{4\beta}{(\beta+1)^2}>\frac{4}{n^2}\gamma\; \Leftrightarrow \; \frac{n}{\ap}<p<\frac{n}{\am}.$$
Since $\ap>n/2$ and $p\geq 2$, then 
$$\frac{4\beta}{(\beta+1)^2}>\frac{4}{n^2}\gamma\; \Leftrightarrow \; p<\frac{n}{\am}.$$
Therefore, it follows from Claim 1 that if $u\in L^p(\Omega\cap B_{\delta_0})$, with $p<n/\am$, then $u\in L^{\frac{n}{n-2}p}(\Omega\cap B_{\delta_0})$. Since $u\in L^2(\Omega\cap B_{\delta_0})$, then \eqref{integ:max} follows.

\medskip\noindent{\bf Step 3:} We claim that for any $\lambda>0$, then 
\begin{equation}\label{eq:iter:6}
\hbox{$|x|^{\frac{n-2}{2}}|u(x)|=O(|x|^{\frac{n-2}{n}\left(\frac{n}{2}-\max\{\am,0\}-\lambda\right)}$ \quad as $x\to 0$.}
\end{equation}
Take $p\in \left(\crit,\frac{n^2}{(n-2)\am}\right)$ if $\gamma>0$, and $p>\crit$ if $\gamma\leq 0$. This is possible since $\crit=2n/(n-2)$ and $\am<n/2$. We fix a sequence $(\eps_i)_i\in (0,+\infty)$ such that $\lim_{i\to +\infty}\eps_i=0$ and we fix a chart $\varphi$ as in \eqref{def:chart:1} to \eqref{def:chart:6}. For any $i\in\nn$, we define
$$u_i(x):=\eps_i^{\frac{n}{p}}u(\varphi(\eps_i x))\hbox{ for all }x\in \tilde{B}_{\delta/\eps_i}.$$
Equation \eqref{regul:eq:2} then rewrites
\begin{equation}\label{regul:eq:3}
-\Delta_{g_i}u_i-\frac{\epsilon_i^2(\gamma+O(\epsilon_i^\tau|x|^\tau))}{|\varphi(\epsilon_i x)|^2}u_i=f_i(x, u_i)\; ;\; u_i=0\hbox{ on }\partial\rnp\cap \tilde{B}_{\delta/\eps_i}
\end{equation}
where $g_i(x):=\varphi^\star\eucl(\epsilon_i x)$ and 
$$|f_i(x,u_i)|\leq C\epsilon_i^2|u_i|+ C \eps_i^{(\crits-2)\left(\frac{n-2}{2}-\frac{n}{p}\right)}|x|^{-s}|u_i|^{\crits-1}$$
in $\tilde{B}_{\delta/\eps_i}$. We fix $R>0$ and we define $\omega_R:=\left(\tilde{B}_R\setminus\tilde{B}_{R^{-1}}\right)\cap\rnp$. With our choice of $p$ above and using \eqref{integ:max}, we get that
\begin{equation}\label{init:bnd:p}
\Vert u_i\Vert_{L^p(\omega_R)}\leq C,
\end{equation}
and
\begin{equation}
|f_i(x,u_i)|\leq C_R|u_i|+ C_R |u_i|^{\crits-1}\quad \hbox{for all $x\in \omega_R$. }
\end{equation}
Fix $q\geq p>\crit$. It follows from elliptic regularity that 
$$\Vert u_i\Vert_{L^q(\omega_{R})}\leq C\; \Rightarrow\; \left\{\begin{array}{ll}
\Vert u_i\Vert_{L^{q'}(\omega_{R/2})}\leq C' &\hbox{ if }q<\frac{n}{2}(\crits-1)\\
\Vert u_i\Vert_{L^{r}(\omega_{R/2})}\leq C' &\hbox{ for all }r\geq 1\hbox{ if }q=\frac{n}{2}(\crits-1)\\
\Vert u_i\Vert_{L^\infty(\omega_{R/2})}\leq C' &\hbox{ if }q>\frac{n}{2}(\crits-1)
\end{array}\right.$$
where $\frac{1}{q'}=\frac{\crits-1}{q}-\frac{2}{n}$ and the constants $C,C'$ are uniform with respect to $i$. It then follows from the standard bootstrap iterative argument and the initial bound \eqref{init:bnd:p} that $\Vert u_i\Vert_{L^\infty(\omega_{R/4})}\leq C'$. Taking $R>0$ large enough and going back to the definition of $u_i$, we get that for all $i\in\nn$, 
$$|x|^{\frac{n}{p}}|u(x)|\leq C\hbox{ for all }x\in\Omega\cap B_{2\eps_i}(0)\setminus B_{\eps_i/2}(0)$$
Since this holds for any sequence $(\eps_i)_i$, we get that $|x|^{\frac{n}{p}}|u(x)|\leq C$ around $0$ for any $\crit<p<\frac{n^2}{(n-2)\am}$ when $\gamma>0$. Letting $p$ go to $\frac{n^2}{(n-2)\am}$ yields \eqref{eq:iter:6} when $\gamma>0$. For $\gamma\leq 0$, we let $p\to +\infty$. This ends Step 3.

\medskip\noindent To finish the proof of Theorem \ref{th:hopf}, we rewrite equation \eqref{regul:eq:2} as
$$-\Delta u-\frac{a(x)}{|x|^2}u=0$$
where
\begin{eqnarray}\nonumber
a(x)&=& \gamma +O(|x|^\tau)+O(|x|^2)+O\left(|x|^{2-s}|u|^{\crits-2}\right)\\
&=&\gamma +O(|x|^\tau)+O(|x|^2)+O\left(|x|^{\frac{n-2}{2}}|u(x)|\right)^{\crits-2}\nonumber
\end{eqnarray}
for all $x\in\Omega$. Since $\am<\frac{n}{2}$, it then follows from \eqref{eq:iter:6} that there exists $\tau'>0$ such that $a(x)=\gamma+O(|x|^{\tau'})$ as $x\to 0$. Therefore we are back to the linear case in Lemma \ref{Hopf.linear} and we are done. 
\hfill$\Box$

\medskip\noindent Here are a few consequences of Theorem \ref{th:hopf}.

\begin{coro}\label{coro:3} Suppose $\gamma <\gamma_H(\Omega)$ and consider the first eigenvalue of the operator $L_\gamma$, that is
$$\lambda_1(\Omega,\gamma):=\inf_{u\in \dundeux\setminus\{0\}}\frac{\int_\Omega\left(|\nabla u|^2-\frac{\gamma}{|x|^2}u^2\right)\, dx}{\int_\Omega u^2\, dx}>0, $$ 
and let $u_0\in \dundeux\setminus\{0\}$ be a minimizer. Then, there exists $A\neq 0$ such that
$$u_0(x)\sim_{x\to 0}A\frac{d(x,\partial \Omega)}{|x|^{\am}}.$$
\end{coro}
\noindent{\it Proof of Corollary \ref{coro:3}:} The existence of a $u_0$ that doesn't change sign is standard. 
The Euler-Lagrange equation is $-\Delta u-\frac{\gamma}{|x|^2}u=k u$ for some $k\in\rr$. We then apply Theorem \ref{th:hopf}.\hfill$\Box$

\begin{coro}\label{coro:2} Suppose $u\in D^{1,2}(\rnp)$, $u\geq 0$, $u\not\equiv 0$ is a weak solution of
$$-\Delta u-\frac{\gamma}{|x|^2}u=\frac{u^{\crit-1}}{|x|^s}\hbox{ in }\rnp.$$
Then, there exist $K_1,K_2>0$ such that
\begin{equation}
u(x)\sim_{x\to 0}K_1\frac{x_1}{|x|^{\am}}\hbox{ and }u(x)\sim_{|x|\to +\infty}K_2\frac{x_1}{|x|^{\ap}}.
\end{equation}
\end{coro}
\noindent{\it Proof of Corollary \ref{coro:2}:} Theorem \ref{th:hopf} yields the behavior when $x\to 0$. The Kelvin transform $\hat{u}(x):=|x|^{2-n}u(x/|x|^2)$ is a solution to the same equation in $D^{1,2}(\rnp)$, and its behavior at $0$ is given by Theorem \ref{th:hopf}. Going back to $u$ yields the behavior at $\infty$. 
\hfill$\Box$

\section{A classification of singular solutions of $L_\gamma u=a(x)u$}\label{sec:sing}

In this section we describe the profile of any positive solution --variational or not-- of linear equations involving $L_\gamma$. Here is the main result of this section. 
 
\begin{theorem}\label{th:classif} Let $u\in C^2(B_\delta(0)\cap (\overline{\Omega}\setminus\{0\}))$ be such that
\begin{equation}\label{eq:u:c:1}
\left\{\begin{array}{ll}
-\Delta u-\frac{\gamma+O(|x|^\tau)}{|x|^2}u=0 &\hbox{\rm in } \Omega\cap B_\delta(0)\\
\hfill u>0&\hbox{\rm  in } \Omega\cap B_\delta(0)\\
\hfill u=0&\hbox{\rm  on } (\partial\Omega\cap B_\delta(0))\setminus \{0\}.\end{array}\right.
\end{equation}
Then, there exists $K>0$ such that
$$\hbox{either }u(x)\sim_{x\to 0}K\frac{d(x,\partial\Omega)}{|x|^{\am}}\quad \hbox{ or }\quad u(x)\sim_{x\to 0}K\frac{d(x,\partial\Omega)}{|x|^{\ap}}.$$
In the first case, the solution $u\in D^{1,2}(\Omega)_{loc,0}$ is a variational solution to (\ref{eq:u:c:1}).
\end{theorem}
\noindent It is worth noting that Pinchover \cite{PinchoverAIHP} proved that the quotient of any two positive solutions to \eqref{eq:u:c:1} has a limit at $0$. \\
The proof will require the following two lemmas. The first gives a Harnack-type inequality. 

\begin{proposition}\label{prop:harnack}
Let $\Omega$ be a smooth bounded domain of $\rn$, and let $a\in L^\infty(\Omega)$ be such that $\Vert a\Vert_\infty\leq M$ for some $M>0$. Assume $U$ is an open subset of $\rn$ and consider $u\in C^2(U\cap\overline{\Omega})$ to be a solution of 
$$\left\{\begin{array}{ll}
-\Delta_g u+au=0&\hbox{ in }U\cap\Omega\\
\hfill u\geq 0&\hbox{ in }U\cap\Omega\\
\hfill u=0&\hbox{ on }U\cap\partial\Omega.
\end{array}\right.$$
Here $g$ is a smooth metric on $U$. If $U'\subset\subset U$ is such that $U'\cap\Omega$ is connected, then there exists $C>0$ depending only on $\Omega, U',M$ and $g$ such that
\begin{equation}\label{ineq:harnack}
\frac{u(x)}{d(x,\partial\Omega)}\leq C\frac{u(y)}{d(y,\partial\Omega)}\hbox{ for all }x,y\in U'\cap\Omega.
\end{equation}
\end{proposition}

\medskip\noindent{\it Proof of Proposition \ref{prop:harnack}:} We first prove a local result. The global result will be the consequence of a covering of $U'$. Fix $x_0\in \partial\Omega$. For $\delta>0$ small enough, there exists a smooth open domain $W$ such that
\begin{equation}\label{ppte:W}
B_\delta(x_0)\cap\Omega\subset W\subset B_{2\delta}(x_0)\cap\Omega\hbox{ and } B_\delta(x_0)\cap\partial W=B_\delta(x_0)\cap\partial \Omega.
\end{equation}
Let $G$ is the Green's function of $-\Delta_g+a$ with Dirichlet boundary condition on $W$, then its representation formula reads as  
 \begin{equation}\label{eq:u:harnack}
 u(x)=\int_{\partial W}u(\sigma)\left(-\partial_{\nu,\sigma}G(x,\sigma)\right)\, d\sigma=\int_{\partial W\setminus\partial\Omega}u(\sigma)\left(-\partial_{\nu,\sigma}G(x,\sigma)\right)\, d\sigma
 \end{equation}
for all $x\in W$, where $\partial_{\nu,\sigma}G(x,\sigma)$ is the normal derivative of $y\mapsto G(x,y)$ at $\sigma\in \partial W$. Estimates of the Green's function (See Robert \cite{r1} and Ghoussoub-Robert \cite{gr2}) yield the existence of $C>0$ such that for all $x\in W$ and $\sigma\in \partial W$, 
$$\frac{1}{C}\frac{d(x,\partial W)}{|x-\sigma|^n}\leq -\partial_{\nu,\sigma}G(x,\sigma)\leq C\frac{d(x,\partial W)}{|x-\sigma|^n}$$
It follows from \eqref{ppte:W} that there exists $C(\delta)>0$ such that for all $x\in B_{\delta/2}(x_0)\cap\Omega\subset W$ and $\sigma\in \partial W\setminus\partial\Omega$, 
$$\frac{1}{C(\delta)}d(x,\partial W)\leq -\partial_{\nu,\sigma}G(x,\sigma)\leq C(\delta)d(x,\partial W)$$
Since $u$ vanishes on $\partial\Omega$, it then follows from \eqref{eq:u:harnack} that for all $x\in B_{\delta/2}(x_0)\cap\Omega$, 
$$\frac{1}{C(\delta)}d(x,\partial W)\int_{\partial W}u(\sigma)\, d\sigma \leq u(x)\leq C(\delta)d(x,\partial W)\int_{\partial W}u(\sigma)\, d\sigma.$$
It is easy to check, that under the assumption \eqref{ppte:W}, we have that $d(x,\partial\Omega)=d(x,\partial W)$. Therefore, we have for all $x\in B_{\delta/2}(x_0)\cap\Omega$, 
$$\frac{1}{C(\delta)}\int_{\partial W}u(\sigma)\, d\sigma \leq \frac{u(x)}{d(x,\partial\Omega)}\leq C(\delta)\int_{\partial W}u(\sigma)\, d\sigma$$
These lower and upper bounds being independent of $x$, we get inequality \eqref{ineq:harnack} for any $x,y\in B_{\delta/2}(x_0)\cap\Omega$. 

\medskip\noindent The general case is a consequence of a covering of $U'\cap \Omega$ by finitely many balls. Note that for balls intersecting $\partial\Omega$, we apply the preceding result, while for balls not intersecting $\partial\Omega$, we apply the classical Harnack inequality. This completes the proof of Proposition \ref{prop:harnack}.\hfill$\Box$\\

In the next proposition we construct  sub- and super solutions with Dirichlet boundary conditions on any smooth domain $W$ of $\rn$ satisfying property \eqref{ppte:W}.

 \begin{proposition}\label{prop:super:dirichlet} Let $\Omega$ be a smooth bounded domain of $\rn$, and let $W$ be a smooth domain of $\rn$ such that for some $r>0$ small enough, we have
 \begin{equation}\label{ppte:W.bis}
B_r(0)\cap\Omega\subset W\subset B_{2r}(0)\cap\Omega\hbox{ and } B_r(0)\cap\partial W=B_r(0)\cap\partial \Omega.
\end{equation}
Fix $\gamma<\frac{n^2}{4}$, $0<\tau\leq 1$ and $\beta\in\rr$ such that $\ap-\tau<\beta<\ap$ and $\beta\neq \am$. Then, up to choosing $r$ small enough, there exists $u_{\ap,+}^{(d)},u_{\ap,-}^{(d)}\in C^\infty(\overline{W}\setminus\{0\})$ such that
\begin{equation}\label{lalala}
\left\{\begin{array}{ll}
\hfill u_{\ap,+}^{(d)},u_{\ap,+}^{(d)}=0 &\hbox{\rm in }\partial W\setminus \{0\}\\
-\Delta u_{\ap,+}^{(d)}-\frac{\gamma+O(|x|^\tau)}{|x|^2}u_{\ap,+}^{(d)}>0&\hbox{\rm in }W\\
-\Delta u_{\ap,-}^{(d)}-\frac{\gamma+O(|x|^\tau)}{|x|^2}u_{\ap,-}^{(d)}<0&\hbox{\rm in }W.
\end{array}
\right.
\end{equation}
Moreover, we have as $x\to 0$, $x\in \Omega$ that
\begin{equation}\label{asymp:u:d1}
u_{\ap,+}^{(d)}(x)=\frac{d(x,\partial\Omega)}{|x|^{\ap}}(1+O(|x|^{\alpha-\beta})) \end{equation}
and 
\begin{equation}\label{asymp:u:d2}
 u_{\ap,-}^{(d)}(x)=\frac{d(x,\partial\Omega)}{|x|^{\ap}}(1+O(|x|^{\alpha-\beta})).
\end{equation}
 \end{proposition}

\medskip\noindent{\it Proof of Proposition \ref{prop:super:dirichlet}:} Take $\eta\in C^\infty(\rn)$ such that $\eta(x)=0$ for $x\in B_{\delta/4}(0)$ and $\eta(x)=1$ for $x\in \rn\setminus B_{\delta/3}(0)$. Define on $W$ the function 
$$f(x):=\left(-\Delta-\frac{\gamma+O(|x|^\tau)}{|x|^2}\right)(\eta u_{\ap,+}), $$
where $u_{\ap,+}$ is given by Proposition \ref{prop:sub:super}.
Note that $f$ vanishes around $0$ and that it is in $C^\infty(\overline{W})$. Let $v\in D^{1,2}(W)$ be such that
$$\left\{\begin{array}{ll}
-\Delta v-\frac{\gamma+O(|x|^\tau)}{|x|^2}v=f&\hbox{ in }W\\
\hfill v=0&\hbox{ on }\partial W.
\end{array}\right.$$
Note that for $r>0$ small enough, $-\Delta-\gamma|x|^{-2}$ is coercive on $W$, and therefore, the existence of $v$ is ensured for small $r$. Define
$$u_{\ap,+}^{(d)}:=u_{\ap,+}-\eta u_{\ap,+}+v.$$
The properties of $W$ and the definition of $\eta$ and $v$ yield
$$\left\{\begin{array}{ll}
\hfill u_{\ap,+}^{(d)}=0 &\hbox{ in }\partial W\setminus \{0\}\\
-\Delta u_{\ap,+}^{(d)}-\frac{\gamma+O(|x|^\tau)}{|x|^2}u_{\ap,+}^{(d)}>0&\hbox{ in }W.
\end{array}
\right.$$
Moreover, since $-\Delta v-(\gamma+O(|x|^\tau))|x|^{-2}v=0$ around $0$ and $v\in D^{1,2}(W)$, it follows from Theorem \ref{th:hopf} that there exists $C>0$ such that $|v(x)|\leq C d(x, W)|x|^{-\am}$ for all $x\in W$. Then \eqref{asymp:u:d1} follows from the asymptotics \eqref{asymp:ua:plus}
of $u_{\ap,+}$ and the fact that $\am<\ap$. We argue similarly for $u_{\ap,-}^{(d)}$. This proves Proposition \ref{prop:super:dirichlet}.\hfill$\Box$\\

\noindent{\it Proof of Theorem \ref{th:classif}:} Let $u$ be a solution of \eqref{eq:u:c:1} as in the statement of Theorem \ref{th:classif}. We claim that
\begin{equation}\label{bnd:fct:positive}
u(x)=O(d(x,\partial\Omega)|x|^{-\ap})\hbox{ for } x\to 0,\; x\in\Omega.
\end{equation}
Indeed, otherwise we can assume that 
\begin{equation}\label{hyp:bnd:absurd}
\limsup_{x\to 0}\frac{u(x)}{d(x,\partial\Omega)|x|^{-\ap}}=+\infty.
\end{equation}
In particular, there exists $(x_k)_{k}\in\Omega$ such that for all $k\in\nn$, 
\begin{equation}\label{eq:u:c:4}
\lim_{k\to +\infty}x_k=0\hbox{ and }\frac{u(x_k)}{d(x_k,\partial\Omega)|x_k|^{-\ap}}\geq k,
 \end{equation}
We claim that there exists then $C>0$ such that
\begin{equation}\label{ineq:lower:harnack}
\hbox{$\frac{u(x)}{d(x,\partial\Omega)|x|^{-\ap}}\geq C k$ for all $x\in \Omega\cap\partial B_{r_k}(0)$, with $r_k:=|x_k|\to 0$.}
\end{equation}
 We prove the claim by using the Harnack inequality \eqref{ineq:harnack}: First take the chart $\varphi$ at $0$ as in \eqref{def:chart:1}, and define
$$u_k(x):=u\circ\varphi(r_k x)\hbox{ for }x\in \rnp\cap B_3(0)\setminus\{0\}.$$
Equation \eqref{eq:u:c:1} rewrites
\begin{equation}\label{eq:u:c:3}
-\Delta_{g_k}u_k+a_k u_k=0\hbox{ in }\rnp\cap B_3(0)\setminus\{0\},
\end{equation}
with $a_k(x):=-r_k^2\frac{\gamma+O(r_k^\tau|x|^\tau)}{|\varphi(r_k x)|^2}$. In particular, there exists $M>0$ such that $|a_k(x)|\leq M$ for all $x\in \rnp\cap B_3(0)\setminus \overline{B}_{1/3}(0)$. Since $u_k\geq 0$, the Harnack inequality \eqref{ineq:harnack} yields the existence of $C>0$ such that
\begin{equation}\label{eq:u:c:2}
\frac{u_k(y)}{y_1}\geq C \frac{u_k(x)}{x_1}\hbox{\quad for all $x,y\in \rnp\cap B_2(0)\setminus \overline{B}_{1/2}(0)$.}
\end{equation}
 Let $\tilde{x}_k\in\rnp$ be such that $x_k=\varphi(r_k \tilde{x}_k)$. In particular, $|\tilde{x}_k|=1+o(1)$ as $k\to +\infty$. It then follows from \eqref{eq:u:c:4}, \eqref{eq:u:c:3} and \eqref{eq:u:c:2} that
$$\frac{u\circ\varphi(r_k y)}{d(\varphi(r_k y),\partial\Omega)}\geq C\cdot k \quad \hbox{for all $y\in \rnp\cap B_2(0)\setminus \overline{B}_{1/2}(0)$.}$$
 In particular, \eqref{ineq:lower:harnack} holds. 

\medskip\noindent We let now $W$ be a smooth domain such that \eqref{ppte:W.bis} holds for $r>0$ small enough. Take the super-solution $u_{\ap,-}^{(d)}$ defined in Proposition \ref{prop:super:dirichlet}. We have that
$$u(x)\geq \frac{C\cdot k}{2}u_{\ap,-}^{(d)}(x)\hbox{ for all }x\in W\cap\partial B_{r_k}(0).$$
Since $u_{\ap,-}^{(d)}$ vanishes on $\partial W$, we have
\begin{equation}\nonumber
u(x)\geq \frac{C\cdot k}{2}u_{\ap,-}^{(d)}(x)\hbox{ for all }x\in \partial(W\cap\partial B_{r_k}(0)).
\end{equation}
Moreover, we have that
$$\hbox{$-\Delta u_{\ap,-}^{(d)}-\frac{\gamma+O(|x|^\tau)}{|x|^2}u_{\ap,-}^{(d)}<0=-\Delta u -\frac{\gamma+O(|x|^\tau)}{|x|^2}u$\quad  on $W$.}
$$
Up to taking $r$ even smaller, it follows from the coercivity of the operator and the maximum principle that
\begin{equation}\label{ineq:u:up:bnd}
u(x)\geq \frac{C\cdot k}{2}u_{\ap,-}^{(d)}(x)\hbox{ for all }x\in W\cap B_{r_k}(0).
\end{equation}
For any $x\in W$, we let $k_0\in\nn$ such that $r_k<|x|$ for all $k\geq k_0$. It then follows from \eqref{ineq:u:up:bnd} that $u(x)\geq \frac{C\cdot k}{2}u_{\ap,-}^{(d)}(x)$ for all $k\geq k_0$. Letting $k\to +\infty$ yields that $u_{\ap,-}^{(d)}(x)$ goes to zero for all $x\in W$. This is a contradiction with \eqref{asymp:u:d2}. Hence \eqref{hyp:bnd:absurd} does not hold, and therefore \eqref{bnd:fct:positive} holds.  

A straightforward consequence of  \eqref{bnd:fct:positive} and Lemma \ref{lem:deriv:1} is that 
there exists $l\in\rr$ such that
\begin{equation}\label{lim:l:3}
\lim_{x\to 0}\frac{u(x)}{d(x,\partial\Omega)|x|^{-\ap}}=l.
\end{equation}

\medskip\noindent We now show the following lemma:

\begin{lem}\label{lem:inter} If $\lim_{x\to 0}\frac{u(x)}{d(x,\partial\Omega)|x|^{-\ap}}=0$, then  $u\in \dundeux_{loc, 0}$ and there exists $K>0$ such that  $u(x)\sim_{x\to 0}K\frac{d(x,\partial \Omega)}{|x|^{\am}}$.
\end{lem}

\smallskip\noindent {\it Proof of Lemma \ref{lem:inter}:} We shall use Theorem \ref{th:hopf}. Take $W$ as in \eqref{ppte:W.bis} and let $\eta\in C^\infty(\rn)$ such that $\eta(x)=0$ for $x\in B_{\delta/4}(0)$ and $\eta(x)=1$ for $x\in \rn\setminus B_{\delta/3}(0)$. Define
$$
\hbox{$f(x):=\left(-\Delta-\frac{\gamma+O(|x|^\tau)}{|x|^2}\right)(\eta u)$ for $x\in W$.}
 $$
The function $f\in C^\infty(\overline{W})$ and vanishes around $0$. Let $v\in D^{1,2}(\Omega)$ be such that
$$\left\{\begin{array}{ll}
-\Delta v-\frac{\gamma+O(|x|^\tau)}{|x|^2}v=f&\hbox{ in }W\\
\hfill v=0&\hbox{ on }\partial W.
\end{array}\right.$$
Note again that for $r>0$ small enough, $-\Delta-\gamma|x|^{-2}$ is coercive on $W$, and therefore, the existence of $v$ is ensured for small $r$. Define
$$\tilde{u}:=u-\eta u+v.$$
The properties of $W$ and the definition of $\eta$ and $v$ yield
$$\left\{\begin{array}{ll}
-\Delta \tilde{u}-\frac{\gamma+O(|x|^\tau)}{|x|^2}\tilde{u}=0&\hbox{ in }W\\
\hfill \tilde{u}=0 &\hbox{ in }\partial W\setminus \{0\}.
\end{array}
\right.$$
Moreover, since $-\Delta v-(\gamma+O(|x|^\tau))|x|^{-2}v=0$ around $0$ and $v\in D^{1,2}(W)$, it follows from Theorem \ref{th:hopf} that there exists $C>0$ such that $|v(x)|\leq C d(x, W)|x|^{-\am}$ for all $x\in W$. Therefore, we have that
\begin{equation}\label{lim:u:0}
\lim_{x\to 0}\frac{\tilde{u}(x)}{d(x,\partial\Omega)|x|^{-\ap}}=0.
\end{equation}
It then follows from Lemma \ref{lem:deriv:1} that
\begin{equation}\label{lim:nablau:0}
\lim_{x\to 0}|x|^{\ap}|\nabla\tilde{u}(x)|=0.
\end{equation}
Let $\psi\in C^\infty_c(W)$ and $w\in D^{1,2}(W)$ be such that
$$\left\{\begin{array}{ll}
-\Delta w-\frac{\gamma+O(|x|^\tau)}{|x|^2}w=\psi&\hbox{ in }W\\
\hfill w=0&\hbox{ on }\partial W.
\end{array}\right.$$
Since $\psi$ vanishes around $0$, it follows from Theorem \ref{th:hopf} and Lemma \ref{lem:deriv:1} that 
\begin{equation}\label{bnd:w}
w(x)=O(d(x,\partial W)|x|^{-\am})\hbox{\quad and \quad}|\nabla w(x)|=O(|x|^{-\am})\quad \hbox{as $x\to 0$.}
\end{equation}
 Fix $\epsilon>0$ small and integrate by parts using that both $\tilde{u}$ and $w$ vanish on $\partial W$, to get 
\begin{eqnarray*}
0&=& \int_{W\setminus B_\epsilon(0)}\left(-\Delta \tilde{u}-\frac{\gamma+O(|x|^\tau)}{|x|^2}\tilde{u}\right)w\, dx\\
&=&  \int_{W\setminus B_\epsilon(0)}\left(-\Delta w-\frac{\gamma+O(|x|^\tau)}{|x|^2}w\right)\tilde{u}\, dx+\int_{\partial (W\setminus B_\epsilon(0))}\left(-w\partial_\nu \tilde{u} +\tilde{u}\partial_\nu w\right)\, d\sigma\\
&=& \int_{W\setminus B_\epsilon(0)}\psi\tilde{u}\, dx-\int_{\Omega\cap \partial B_\epsilon(0)}\left(-w\partial_\nu \tilde{u} +\tilde{u}\partial_\nu w\right)\, d\sigma.
\end{eqnarray*}
Using the upper-bounds \eqref{lim:u:0}, \eqref{lim:nablau:0} and \eqref{bnd:w}, and that $\psi$ vanishes around $0$, we get
\begin{eqnarray*}
0&=& \int_{W\setminus B_\epsilon(0)}f\tilde{u}\, dx+o\left(\epsilon^{n-1}(\epsilon^{1-\am}\epsilon^{-\ap}+\epsilon^{1-\ap}\epsilon^{-\am})\right)\\
&=& \int_{W\setminus B_\epsilon(0)}f\tilde{u}\, dx+o(1), \quad \hbox{as $\epsilon\to 0$.}
\end{eqnarray*}
Therefore, we have
$\int_W f\tilde{u}\, dx=0$
for all $f\in C_c^\infty(W)$. Since $\tilde{u}\in L^p$ is smooth outside $0$, we then get that $\tilde{u}\equiv 0$, and therefore $u=\eta u+v$. In particular, $u\in D^{1,2}(\Omega)_{loc,0}$ is a distributional positive solution to $-\Delta u-\frac{\gamma+O(|x|^\tau)}{|x|^2}u=0$ on $W$. It then follows from Theorem \ref{th:hopf} that there exists $K>0$ such that 
$u(x)\sim_{x\to 0}K\frac{d(x,\partial \Omega)}{|x|^{\am}}$. This prove Lemma \ref{lem:inter}.\hfill$\Box$

\medskip\noindent Combining Lemma \ref{lem:inter} with (\ref{lim:l:3}) completes the proof of Theorem \ref{th:classif}.\hfill$\Box$


\medskip\noindent As a consequence of Theorem \ref{th:classif}, we 
improve Lemma \ref{lem:funda} as follows.
\begin{proposition}\label{prop:funda:2} Let $u\in C^2(\overline{\rnp}\setminus \{0\})$ be a nonnegative function such that
\begin{equation}\label{eq:u:36}
-\Delta u-\frac{\gamma}{|x|^2}u=0\hbox{ in }\rnp\; ;\; u=0\hbox{ on }\partial \rnp.
\end{equation}
Then there exist $\lambda_-,\lambda_+\geq 0$ such that
$$u(x)=\lambda_- x_1|x|^{-\am}+\lambda_+ x_1|x|^{-\ap}\hbox{ for all }x\in \rnp.$$
\end{proposition}
\noindent{\it Proof of Proposition \ref{prop:funda:2}:} Without loss of generality, we assume that $u\not\equiv 0$, so that $u>0$. We consider the Kelvin transform of $u$ defined by $\hat{u}(x):=|x|^{2-n}u(x/|x|^2)$  for all $x\in\rnp$. Both $u$ and $\hat{u}$ are then nonnegative solutions of \eqref{eq:u:36}. It follows from Theorem \ref{th:classif} that, after performing back the Kelvin transform, there exist $\alpha_1,\alpha_2\in \{\ap,\am\}$ such that
$$\lim_{x\to 0}\frac{u(x)}{x_1|x|^{-\alpha_1}}=l_1>0\hbox{ and }\lim_{|x|\to \infty}\frac{u(x)}{x_1|x|^{-\alpha_2}}=l_2>0.$$
If $\alpha_1\leq \alpha_2$, then $u(x)\leq C x_1|x|^{-\alpha_1}$ for all $x\in\rnp$. The result then follows from Lemma \ref{lem:funda}. If $\alpha_1>\alpha_2$, then $\alpha_1=\ap$ and $\alpha_2=\am$. We define 
$$\tilde{u}(x):=u(x)-l_1 x_1|x|^{-\ap}\hbox{ for all }x\in\rnp.$$
to obtain that 
\begin{equation*}
-\Delta \tilde{u}-\frac{\gamma}{|x|^2}\tilde{u}=0\hbox{ in }\rnp\; ;\; \tilde{u}=0\hbox{ on }\partial \rnp,
\end{equation*}
and 
$\tilde{u}(x)=o(x_1|x|^{-\ap})\hbox{ as }x\to 0.$
Arguing as in the proof of Lemma \ref{lem:inter}, we get that $\tilde{u}\in D^{1,2}(\rnp)_{loc,0}$ and $\tilde{u}(x)=O(x_1|x|^{-\am})$ as $x\to 0$. Moreover, we have that $\tilde{u}(x)=(l_2+o(1))x_1|x|^{-\am}$ as $|x|\to +\infty$, therefore $\tilde{u}(x)>0$ for $|x|>>1$. Since $\tilde{u}\in D^{1,2}(\rnp)_{loc,0}$, the comparison principle then yields $\tilde{u}>0$ everywhere. We also have that $\tilde{u}(x)\leq C x_1|x|^{-\am}$ for all $x\in\rnp$. It then follows from Lemma \ref{lem:funda} that there exists $\lambda_-\geq 0$ such that $\tilde{u}(x)=\lambda_- x_1|x|^{-\am}$ for all $x\in\rnp$. We then get the conclusion of Proposition \ref{prop:funda:2}. \hfill$\Box$

\section{The Hardy singular b-mass of a domain in the case $\gamma>\frac{n^2-1}{4}$}\label{sec:mass}

We shall proceed in the following theorem to define the mass of a smooth bounded domain $\Omega$ of $\rn$ such as $0\in\partial\Omega$. It will involve the expansion of positive singular solutions of the Dirichlet boundary problem $L_\gamma u=0$. 

\begin{theorem}\label{def:mass} Let $\Omega$ be a smooth bounded domain $\Omega$ of $\rn$ such as $0\in\partial\Omega$, and assume that $\frac{n^2-1}{4}<\gamma<\gamma_H(\Omega)$. Then, up to multiplication by a positive constant, there exists a unique function $H\in C^2(\overline{\Omega}\setminus \{0\})$ such that
\begin{equation}\label{carac:H}
-\Delta H-\frac{\gamma}{|x|^2}H=0\hbox{ in }\Omega\; ,\; H>0\hbox{ in }\Omega\; ,\; H=0\hbox{ on }\partial\Omega.
\end{equation}
Moreover, there exists $c_1>0$ and $c_2\in \rr$ such that 
\begin{equation}\label{def:c1:c2}
\hbox{$H(x)=c_1\frac{d(x,\partial\Omega)}{|x|^{\ap}}+c_2\frac{d(x,\partial\Omega)}{|x|^{\am}}+ o\left(\frac{d(x,\partial\Omega)}{|x|^{\am}}\right)$ as $x\to 0$.}
\end{equation}
The quantity $m_\gamma(\Omega):=\frac{c_2}{c_1}\in \rr$, which is independent of the choice of $H$ satisfying \eqref{carac:H}, will be called the Hardy b-mass of $\Omega$ associated to $L_\gamma$.  
\end{theorem}

\medskip\noindent{\it Proof of Theorem \ref{def:mass}.}  First, we start by constructing a singular solution $H_0$ for \eqref{carac:H}. For that, consider $u_{\ap}$ as in \eqref{def:ua} and let $\eta\in C^\infty_c(\rn)$ be such that $\eta(x)=1$ for $x\in B_{\delta/2}(0)$ and $\eta(x)=0$ for $x\in \rn\setminus B_\delta(0)$.  Set
$$f:=-\Delta (\eta u_{\ap})-\frac{\gamma}{|x|^2}(\eta u_{\ap})\hbox{ in }\overline{\Omega}\setminus\{0\}.$$
It follows from \eqref{est:ua:1} that $f$ is smooth outside $0$ and that
$$f(x)=O\left(d(x,\partial\Omega)|x|^{-\ap-1}\right)=O\left(|x|^{-\ap}\right)\hbox{ in }\Omega\cap B_{\delta/2}(0).$$
Since $\gamma>\frac{n^2-1}{4}$, we have that $\ap<\frac{n+1}{2}$, and therefore $f\in L^{\frac{2n}{n+2}}(\Omega)=\left(L^{\crit}(\Omega)\right)^\prime\subset \left(\dundeux\right)^\prime$. It then follows from the coercivity assumption $\gamma<\gamma_H(\Omega)$ that there exists $v\in \dundeux$ such that
$$
-\Delta v-\frac{\gamma}{|x|^2}v=f\hbox{ in }\left(\dundeux\right)^\prime.
$$
Let $v_1,v_2\in \dundeux$ be such that
\begin{equation}\label{eq:fplus:fmoins}
-\Delta v_1-\frac{\gamma}{|x|^2}v_1=f_+\hbox{ and }-\Delta v_2-\frac{\gamma}{|x|^2}v_2=f_-\hbox{ in }\left(\dundeux\right)^\prime.
\end{equation}
In particular, $v=v_1-v_2$ and $v_1,v_2\in C^1(\overline{\Omega}\setminus\{0\})$, and they vanish on $\partial\Omega\setminus\{0\}$. 
\medskip\noindent Assume that $f_+\not\equiv 0$. Since $f_+\geq 0$, the comparison principle yields $v_1>0$ on $\Omega\setminus\{0\}$ and $\partial_\nu v_1<0$ on $\partial\Omega\setminus\{0\}$. Therefore, for any $\delta>0$ small enough, there exists $C(\delta)>0$ such that 
\begin{equation}\nonumber
v_1(x)\geq C(\delta) d(x,\partial\Omega)\hbox{ for all }x\in \partial B_\delta(0)\cap\Omega.
\end{equation}
Let $u_{\am,-}$ be the sub-solution defined in \eqref{ppty:ua}. It follows from the asymptotic \eqref{asymp:ua:plus} that there exists $C'(\delta)>0$ such that
$$v_1\geq C'(\delta)u_{\am,-}\hbox{ in } \partial B_\delta(0)\cap\Omega.$$
Since this inequality also holds on $\partial(B_\delta(0)\cap\Omega)$ and that 
$$
\hbox{$(-\Delta-\frac{\gamma}{|x|^2})(v_1- C'(\delta)u_{\am,-})\geq 0$\quad  in $B_\delta(0)\cap\Omega$,}
$$
 coercivity and the maximum principle yield $v_1\geq C'(\delta)u_{\am,-}$ in $B_\delta(0)\cap\Omega$. It then follows from \eqref{asymp:ua:plus} that there exists $c>0$ such that
$$v_1(x)\geq c\cdot d(x,\partial\Omega)|x|^{-\am}\hbox{ in }B_\delta(0)\cap\Omega.$$
Therefore, we have that
$$f_+(x)\leq C d(x,\partial\Omega)|x|^{-\ap-1}\leq \frac{C}{c}|x|^{\am-\ap-1}v_1(x)\leq \frac{C}{c}|x|^{\am-\ap+1}\frac{v_1(x)}{|x|^2}$$
in $B_\delta(0)\cap\Omega$. Therefore, \eqref{eq:fplus:fmoins} yields
$$-\Delta v_1+\frac{\gamma+O(|x|^{\am-\ap+1})}{|x|^2}v_1=0 \hbox{ in }B_\delta(0)\cap\Omega.$$
Since $\gamma>\frac{n^2-1}{4}$, we have that $\am-\ap+1>0$. Since $v_1\in \dundeux$, $v_1\geq 0$ and $v_1\not\equiv 0$, it follows from Theorem \ref{th:hopf} that there exists $K_1>0$ such that 
\begin{equation}\label{eq:v1}
\hbox{$v_1(x)=K_1\frac{d(x,\partial\Omega)}{|x|^{\am}}+o\left(\frac{d(x,\partial\Omega)}{|x|^{\am}}\right)$\quad as $x\to 0$. }
\end{equation}
If $f_+\equiv 0$, then $v_1\equiv 0$ and \eqref{eq:v1} holds with $K_1=0$. Arguing similarly for $f_-$, we then get that there exists $K_1,K_2\geq 0$ such that for any $i=1,2$, we have that
\begin{equation}\nonumber
\hbox{$v_i(x)=K_i\frac{d(x,\partial\Omega)}{|x|^{\am}}+o\left(\frac{d(x,\partial\Omega)}{|x|^{\am}}\right)$\quad when $x\to 0$. }
\end{equation}
Since $v=v_1-v_2$, we then get that there exists $K\in\rr$ such that
\begin{equation}\label{asymp:v}
\hbox{$v(x)=-K\frac{d(x,\partial\Omega)}{|x|^{\am}}+o\left(\frac{d(x,\partial\Omega)}{|x|^{\am}}\right)$\quad as $x\to 0$. }
\end{equation}
\medskip\noindent Set
\begin{equation}\nonumber
H_0(x):=\eta(x)u_{\ap}(x)-v(x)\hbox{ for all }x\in\overline{\Omega}\setminus\{0\}.
\end{equation}
It follows from the definition of $v$ and the regularity outside $0$ that
\begin{equation}\nonumber
-\Delta H_0-\frac{\gamma}{|x|^2}H_0=0\hbox{ in }\Omega\, ;\, H_0(x)=0\hbox{ in }\partial\Omega\setminus\{0\}.
\end{equation}
Moreover, the asymptotics \eqref{asymp:ua:plus} and \eqref{asymp:v} yield $H_0(x)>0$ on $\Omega\cap B_{\delta'}(0)$ for some $\delta'>0$ small enough. It  follows from the comparison principle that $H_0>0$ in $\Omega$. 

\noindent We now perform an expansion of $H_0$. First note that from the definition \eqref{def:ua} of $u_{\ap}$, the asymptotic \eqref{asymp:v} of $v$ and the fact that $\ap-\am<1$, we have
\begin{eqnarray}
H_0(x)&=&\frac{d(x,\partial\Omega)}{|x|^{\ap}}(1+O(|x|))+K\frac{d(x,\partial\Omega)}{|x|^{\am}}+o\left(\frac{d(x,\partial\Omega)}{|x|^{\am}}\right)\nonumber\\
&=&\frac{d(x,\partial\Omega)}{|x|^{\ap}}+K\frac{d(x,\partial\Omega)}{|x|^{\am}}+o\left(\frac{d(x,\partial\Omega)}{|x|^{\am}}\right)\nonumber
\end{eqnarray}
as $x\to 0$. In particular, since in addition $H_0>0$ in $\Omega$, there exists $c>1$ such that
\begin{equation}\label{est:H:3}
\hbox{$\frac{1}{c}\frac{d(x,\partial\Omega)}{|x|^{\ap}}\leq H_0(x)\leq c\frac{d(x,\partial\Omega)}{|x|^{\ap}}$ \quad for all $x\in \Omega$.}
\end{equation}

\noindent Finally, we establish the uniqueness. For that, we let $H\in C^2(\overline{\Omega}\setminus\{0\})$ be as in \eqref{carac:H} and set
$$\lambda_0:=\max\{\lambda\geq 0/\; H\geq \lambda H_0\}.$$
The number $\lambda_0$ is clearly defined, and so we set $\tilde{H}:=H-\lambda_0 H_0\geq 0$. Assume that $\tilde{H}\not\equiv 0$. Since $-\Delta\tilde{H}-\gamma|x|^{-2}\tilde{H}=0$, it follows from Theorem \ref{th:classif} that there exists $\alpha\in \{\ap,\am\}$ and $K>0$ such that 
\begin{equation}\label{est:H}
H(x)\sim_{x\to 0}K\frac{d(x,\partial\Omega)}{|x|^{\alpha}}.
\end{equation}
If $\alpha=\am$, then $\tilde{H}\in \dundeux$ is a variational solution to $-\Delta\tilde{H}-\frac{\gamma}{|x|^2}\tilde{H}=0$ in $\Omega$. Then coercivity then yields that $\tilde{H}\equiv 0$, contradicting the initial hypothesis.

\smallskip\noindent Therefore $\alpha=\ap$. Since $\tilde{H}>0$ vanishes on $\partial\Omega\setminus\{0\}$, there exists for any $\delta>0$,  $c(\delta)>0$ such that
\begin{equation}\label{est:H:2}
\tilde{H}(x)\geq c(\delta)d(x,\partial\Omega)\hbox{ for }x\in \Omega\setminus B_\delta(0).
\end{equation}
Therefore, \eqref{est:H}, \eqref{est:H:2} and \eqref{est:H:3} yield the existence of $c>0$ such that $\tilde{H}\geq c H_0$, and then $H\geq (\lambda_0+c)H_0$, contradicting the definition of $\lambda_0$. It follows that $\tilde{H}\equiv 0$, which means that  $H=\lambda_0 H_0$ for some $\lambda_0>0$. This proves uniqueness and completes the proof of Theorem \ref{def:mass}.\hfill$\Box$

Now we establish the monotonicity of the mass with respect to set inclusion.

\begin{proposition}\label{prop:mono} The Hardy b-mass  is strictly increasing in the following sense: Assume $\Omega_1,\Omega_2$ are two smooth bounded domains such that $0\in \partial\Omega_1\cap\partial\Omega_2$, and $\frac{n^2-1}{4}<\gamma<\min\{\gamma_H(\Omega_1),\gamma_H(\Omega_2)\}$, then
\begin{equation}\label{ineq:mono}
\Omega_1\subsetneq \Omega_2\; \Rightarrow \; m_\gamma(\Omega_2)<m_\gamma(\Omega_1).
\end{equation}
Moreover, if $\Omega\subsetneq \rnp$ and  $\frac{n^2-1}{4}<\gamma <\frac{n^2}{4}$, then $m_\gamma(\Omega)<0$.
\end{proposition}
\noindent{\it Proof of Proposition \ref{prop:mono}:} It follows from the definition of the mass that for $i=1,2$, there exists $H_i\in C^2(\overline{\Omega_i}\setminus\{0\})$ such that
\begin{equation}\label{carac:Hi}
-\Delta H_i-\frac{\gamma}{|x|^2}H_i=0\hbox{ in }\Omega_i\; ,\; H_i>0\hbox{ in }\Omega_i\; ,\; H_i=0\hbox{ on }\partial\Omega_i,
\end{equation}
with 
\begin{equation}\label{exp:Hi}
H_i(x)=\frac{d(x,\partial\Omega_i)}{|x|^{\ap}}+m_{\gamma}(\Omega_i)\frac{d(x,\partial\Omega_i)}{|x|^{\am}}+o\left(\frac{d(x,\partial\Omega_i)}{|x|^{\am}}\right)
\end{equation}
as $x\to 0$, $x\in \Omega_i$. Set $h:=H_2-H_1$ on $\Omega_1$. Since $\Omega_1\subsetneq \Omega_2$, we have that
\begin{equation}\label{syst:h}
\left\{\begin{array}{ll}
-\Delta h-\frac{\gamma}{|x|^2}h=0 &\hbox{ in }\Omega_1\\
\hfill h\geq 0&\hbox{ on }\partial\Omega_1\\
\hfill h\not\equiv 0&\hbox{ in }\partial\Omega_1.
\end{array}\right.\end{equation}
First, we  claim that $h\in H^{1,2}(\Omega_1)$. Indeed, it follows from the construction of the singular function in the previous theorem, that there exists $w\in H^{1,2}(\Omega_1)$ such that
\begin{equation}\label{eq:h}
h(x)=\frac{d(x,\partial\Omega_2)-d(x,\partial\Omega_1)}{|x|^{\ap}}+w(x)\hbox{ for all }x\in \Omega_1.
\end{equation}
Since $\Omega_1\subset\Omega_2$ and $0$ is on the boundary of both domains, then the tangent spaces at $0$ of $\Omega_1$ and $\Omega_2$ are equal, and one gets that 
$$\hbox{$d(x,\partial\Omega_1)-d(x,\partial\Omega_2)=O(|x|^2)$\quad as $x\to 0$. }$$
Since $\ap-\am<1$, we then get that
\begin{equation}\nonumber
\tilde{h}(x):=\frac{d(x,\partial\Omega_2)-d(x,\partial\Omega_1)}{|x|^{\ap}}=O(|x|^{1-\am})\hbox{ as }x\to 0.
\end{equation}
Similarly, $|\nabla \tilde{h}(x)|=O(|x|^{-\am})$ as $x\to 0$, we deduce that $\tilde{h}\in H^{1,2}(\Omega_1)$. It then follows from \eqref{eq:h} that $h\in H^{1,2}(\Omega_1)$.\\  
\noindent To prove the monotonicity, note first that since $\gamma<\gamma_H(\Omega_1)$ and $h\in H^{1,2}(\Omega_1)$, it follows from \eqref{syst:h} and the comparison principle that $h\geq 0$ in $\Omega_1$. Since $h\not\equiv 0$, it follows from Hopf's maximum principle that for any $\delta>0$ small, there exists $C(\delta)>0$ such that
$$h(x)\geq C(\delta)d(x,\partial\Omega_1)\hbox{ for all }x\in \partial B_\delta(0)\cap\Omega_1.$$
We define the sub-solution $u_{\am,-}$ as in Proposition \ref{prop:sub:super}. It then follows from the inequality above and the asymptotics in \eqref{asymp:ua:plus} that there exists $\epsilon_0>0$ such that
$$h(x)\geq 2\epsilon_0 u_{\am,-}(x)\hbox{ for all }x\in \partial B_\delta(0)\cap\Omega_1.$$
This inequality also holds on $B_\delta(0)\cap\partial\Omega_1$ since $u_{\am,-}$ vanishes on $\partial\Omega_1$. It the follows from the maximum principle that $h(x)\geq 2\epsilon_0 u_{\am,-}(x)$ for all $x\in B_\delta(0)\cap\Omega_1$. With the definition of $h$ and the asymptotic \eqref{asymp:ua:plus}, we then that for $\delta'>0$ small enough, we have that
\begin{equation}\label{lower:bnd}
H_2(x)-H_1(x)\geq \epsilon_0 \frac{d(x,\partial\Omega_1)}{|x|^{\am}}\hbox{ for all }x\in B_{\delta'}(0)\cap\Omega_1.
\end{equation}
We let $\vec{\nu}$ be the inward orthonormal vector of $\partial\Omega_1$ at $0$. This is also the inward orthonormal vector of $\partial\Omega_2$ at $0$. Therefore, for any $t>0$ small enough, we have that $d(t\vec{\nu},\partial\Omega_i)=t$ for $i=1,2$. It then follows from the expressions \eqref{exp:Hi} and \eqref{lower:bnd} that
$$
\hbox{$(m_\gamma(\Omega_2)-m_\gamma(\Omega_1))\frac{t}{t^{\am}}+o\left(\frac{t}{t^{\am}}\right)\geq \epsilon_0\frac{t}{t^{\am}}$\quad as $t\downarrow 0$.}$$
 We then get that $m_\gamma(\Omega_2)-m_\gamma(\Omega_1)\geq \epsilon_0$, and therefore $m_\gamma(\Omega_2)>m_\gamma(\Omega_1)$. This proves \eqref{ineq:mono} and ends the first part of Proposition \ref{prop:mono}.
 
 \medskip\noindent Finally, now that $H_2(x):=\frac{x_1}{|x|^{\ap}}$ is a solution for $\rnp$. If now  $\Omega\subset\rnp$, then the above yields that $0=m_\gamma(\rnp)>m_\gamma(\Omega)$, which completes the proof of Proposition \ref{prop:mono}.\hfill$\Box$

\medskip\noindent In Section \ref{sec:ex:mass}, we will prove that one can define the mass  $m_\gamma(\rnp)$ of $\rnp$, and that $m_\gamma(\rnp)=0$.

\section{Test functions and the existence of extremals in the Hardy-Sobolev inequalities}\label{sec:best:cst}
\noindent Let $\Omega$ be a domain of $\rn$ such that $0\in \partial\Omega$. For $\gamma\in\rr$ and $s\in [0,2)$, recall that 
\begin{equation}\label{def:mu}
\mu_{\gamma,s}(\Omega):=\inf_{u\in D^{1,2}(\Omega)\setminus\{0\}}J^\Omega_{\gamma,s}(u),
\end{equation}
where
\begin{equation}\nonumber
J^\Omega_{\gamma,s}(u):=\frac{\int_\Omega\left(|\nabla u|^2-\frac{\gamma}{|x|^2}u^2\right)\, dx}{\left(\int_\Omega\frac{|u|^{\crit}}{|x|^s}\, dx\right)^{\frac{2}{\crit}}}.
\end{equation}
Note that critical points $u\in D^{1,2}(\Omega)$ of $J^\Omega_{\gamma,s}$ are weak solutions to the pde
\begin{equation}\label{equation}
-\Delta u-\frac{\gamma}{|x|^2}=\lambda \frac{|u|^{\crit-2}u}{|x|^s},
\end{equation}
for some $\lambda \in \rr$, which can be rescaled to be equal to $1$ if $\lambda >0$ and to be $-1$ is $\lambda <0$. In this section, we investigate the existence of minimizers for $J^\Omega_{\gamma,s}$.  We start with the following easy case, where we don't have extremals.

\begin{proposition}\label{prop:s0:gamma:neg} Let $\Omega\subset \rn$ be a smooth domain such that $0\in\partial\Omega$ (No boundedness is assumed).
 When $s=0$ and $\gamma\leq 0$, we have that $\mu_{\gamma,0}(\Omega)=\frac{1}{K(n,2)^2}$ (where $K(n,2)=\mu_{0,0}(\rn)$ is the best constant in the Sobolev inequality  \eqref{def:K}) and there is no extremal.
\end{proposition}
\smallskip\noindent{\it Proof of Proposition \ref{prop:s0:gamma:neg}:}  Note that $\crits=\crit(0)=\crit$. Since $\gamma\leq 0$, we have for any $u\in C^\infty_c(\Omega)\setminus\{0\}$, 
\begin{equation}\label{ineq:ext:s:0}
\frac{\int_{\Omega}\left(|\nabla u|^2-\gamma\frac{u^2}{|x|^2}\right)\, dx}{\left(\int_\Omega |u|^{\crit}\, dx \right)^{\frac{2}{\crit}}}\geq \frac{\int_{\Omega}|\nabla u|^2\, dx}{\left(\int_\Omega |u|^{\crit}\, dx \right)^{\frac{2}{\crit}}}\geq \frac{1}{K(n,2)^2},
\end{equation}
and therefore $\mu_{\gamma,0}(\Omega)\geq \frac{1}{K(n,2)^2}$. Fix now $x_0\in\Omega$ and let $\eta\in C^\infty_c(\Omega)$ be such that $\eta(x)=1$ around $x_0$. Set $\ue(x):=\eta(x)\left(\frac{\eps}{\eps^2+|x-x_0|^2}\right)^{\frac{n-2}{2}}$ for all $x\in\Omega$ and $\eps>0$. Since $x_0\neq 0$, it is easy to check that $\lim_{\eps\to 0}\int_\Omega \frac{\ue^2}{|x|^2}\, dx=0$. It is also classical (see for example Aubin \cite{aubin}) that $$\lim_{\eps\to 0}\frac{\int_{\Omega}|\nabla \ue|^2\, dx}{\left(\int_\Omega |\ue|^{\crit}\, dx \right)^{\frac{2}{\crit}}}= \frac{1}{K(n,2)^2}.$$ 
It follows that $\mu_{\gamma,0}(\Omega)\leq \frac{1}{K(n,2)^2}$. This proves that $\mu_{\gamma,0}(\Omega)= \frac{1}{K(n,2)^2}$.

\medskip\noindent Assume now that there exists an extremal $u_0$ for $\mu_{\gamma,0}(\Omega)$ in $\dundeux\setminus\{0\}$. The inequalities in \eqref{ineq:ext:s:0} and the fact that 
$$\frac{\int_{\Omega}|\nabla u_0|^2\, dx}{\left(\int_\Omega |u_0|^{\crit}\, dx \right)^{\frac{2}{\crit}}}= \frac{1}{K(n,2)^2},$$
means that $u_0\in \dundeux\subset D^{1,2}(\rn)$ is an extremal for the classical Sobolev inequality on $\rn$. But these extremals are known (see Aubin \cite{aubin} or Talenti \cite{Tal}) and their support is the whole of $\rn$, which is a contradiction since $u_0$ has support in $\Omega\neq \rn$. It follows that there is no extremal for $\mu_{\gamma,0}(\Omega)$. This proves Proposition \ref{prop:s0:gamma:neg}.\hfill$\Box$

The remainder of the section is devoted to the proof of  the following.

\begin{theorem} \label{main} Let $\Omega$ be a smooth bounded domain in $\rn$ ($n\geq 3$) such that $0\in \partial \Omega$ and let $0\leq s< 2$ and $ \gamma <\frac{n^2}{4}$. Assume that either  $s>0$, or that $\{s=0$, $n\geq 4$ and $\gamma>0\}$. There are then extremals for $\mu_{\gamma,s}(\Omega)$ under one of the following two conditions:
\begin{enumerate} 
\item  $\gamma\leq\frac{n^2-1}{4}$ and the mean curvature of $\partial \Omega$ at $0$ is negative. 
\item $\gamma>\frac{n^2-1}{4}$ and the mass $m_\gamma (\Omega)$ of $\Omega$ is positive. 
\end{enumerate} 
Moreover, if $\gamma <\gamma_H(\Omega)$ (resp., $\gamma \geq \gamma_H(\Omega)$), then such extremals are positive solutions for (\ref{equation}) with $\lambda >0$ (resp., $\lambda \leq 0$).

\end{theorem}
\noindent The remaining case $n=3$, $s=0$ and $\gamma >0$ will be dealt with in section 11. 

\medskip\noindent According to Theorem \ref{tool}, in order to establish existence of extremals,  it suffices to show that  $\mu_{\gamma,s}(\Omega)<\mu_{\gamma,s}(\rnp)$. The rest of the section consists of showing that the above mentioned geometric conditions lead to such gap.

\smallskip\noindent  In the sequel, $h_\Omega(0)$ will denote the mean curvature of $\partial \Omega$ at $0$. The orientation is chosen such that the mean curvature of the canonical sphere (as the boundary of the ball) is positive.  Since $\{s>0\}$, or that $\{s=0$, $n\geq 4$ and $\gamma>0\}$, it follows from Proposition \ref{prop:ext:rnp} in Section \ref{sec:app:1} of the appendix (see also Bartsch-Peng-Zhang \cite{BPZ} and Chern-Lin \cite{CL5}) that there are extremals for $\mu_{\gamma,s}(\rnp)$. The following proposition combined with Theorem \ref{tool} clearly yield the  claims in Theorem \ref{main}. 

\begin{proposition}\label{prop:test:fct} We fix $\gamma<\frac{n^2}{4}$. Assume that there are extremals for $\mu_{\gamma,s}(\rnp)$. There exist then two families $(\ue^1)_{\eps>0}$ and $(\ue^2)_{\eps>0}$ in $\dundeux$, and two positive constants $c^1_{\gamma,s}$ and $c^2_{\gamma,s}$ such that:
\begin{enumerate}

\item For $\gamma<\frac{n^2-1}{4}$, we have that
\begin{equation}\label{asymp:J:gamma:small}
J(u^1_\epsilon)=\mu_{\gamma,s}(\rnp)\left(1+c^1_{\gamma,s}\cdot h_\Omega(0)\cdot \eps +o(\eps)\right)\hbox{ when }\eps\to 0.
\end{equation}
\item For $\gamma=\frac{n^2-1}{4}$, we have that
\begin{equation}\label{asymp:J:gamma:crit}
J(u^1_\epsilon)=\mu_{\gamma,s}(\rnp)\left(1+c^1_{\gamma,s}\cdot h_\Omega(0)\cdot \eps\ln\frac{1}{\eps} +o\left(\eps\ln\frac{1}{\eps}\right)\right)\hbox{ when }\eps\to 0.
\end{equation}
\item For $\gamma>\frac{n^2-1}{4}$, we have as $\epsilon\to 0$, 
that
\begin{equation}\label{asymp:J:gamma:large}
J(u^2_\epsilon)=\mu_{\gamma,s}(\rnp)\left(1-c^2_{\gamma,s}\cdot m_\gamma(\Omega)\cdot \eps^{\ap-\am}+o(\eps^{\ap-\am})\right).
\end{equation}
\end{enumerate}
\end{proposition}
\noindent{\bf Remark:} When $\gamma<\frac{n^2-1}{4}$, this result is due to Chern-Lin \cite{CL5}. Actually, they stated the result for $\gamma<\frac{(n-2)^2}{4}$, but their proof works for $\gamma<\frac{n^2-1}{4}$. However, when $\gamma\geq\frac{n^2-1}{4}$, we need the exact asymptotic profile of $U$ that was described by Corollary \ref{coro:2}.

\medskip\noindent{\it Proof of Proposition \ref{prop:test:fct}:} By assumption, there exists $U\in D^{1,2}(\rnp)\setminus \{0\}$, $U\geq 0$, that is a minimizer for $\mu_{\gamma,s}(\rnp)$. In other words, 
\begin{equation}\nonumber
J^{\rnp}_{\gamma,s}(U)=\frac{\int_{\rnp}\left(|\nabla U|^2-\frac{\gamma}{|x|^2}U^2\right)\, dx}{\left(\int_{\rnp}\frac{|U|^{\crits}}{|x|^s}\, dx\right)^{\frac{2}{\crits}}}=\mu_{\gamma,s}(\rnp).
\end{equation}
Therefore (see Corollary \ref{coro:2}), there exists $\lambda>0$ such that
\begin{equation}\left\{\begin{array}{ll}\label{eq:U}
-\Delta U-\frac{\gamma}{|x|^2}U=\lambda\frac{U^{\crits-1}}{|x|^s}&\hbox{ in }\rnp\\
\hfill U>0 &\hbox{ in }\rnp\\
\hfill U=0&\hbox{ in }\partial\rnp
\end{array}\right\}\end{equation}
and  there exists $K_1,K_2>0$ such that
\begin{equation}\label{asymp:U}
U(x)\sim_{x\to 0}K_1\frac{x_1}{|x|^{\ams}}\hbox{ and }U(x)\sim_{|x|\to +\infty}K_2\frac{x_1}{|x|^{\aps}},
\end{equation}
where here and in the sequel, we write for convenience
$$\aps:=\ap\hbox{ and }\ams:=\am.$$
In particular, it follows from Lemma \ref{lem:deriv:1} (after reducing all limits to happen at $0$ via the Kelvin transform) that there exists $C>0$ such that
\begin{equation}\label{bnd:U}
\hbox{$U(x)\leq C x_1 |x|^{-\aps}\hbox{ and }|\nabla U(x)|\leq C |x|^{-\aps}$ for all $x\in\rnp$.}
\end{equation}
We shall now construct suitable test-function for each range of $\gamma$. First note that 

\begin{eqnarray}
\gamma<\frac{n^2-1}{4} & \Leftrightarrow & \aps-\ams>1\nonumber\\
\gamma=\frac{n^2-1}{4} & \Leftrightarrow & \aps-\ams=1.\nonumber
\end{eqnarray}
\noindent Concerning terminology, here and in the sequel, we define as in \eqref{def:ball}
$$\tilde{B}_{r}:=(-r,r)\times B_{r}^{(n-1)}(0)\subset \rr\times\rr^{n-1},$$
for all $r>0$ and
$$V_+:=V\cap\rnp$$
for all $V\subset\rn$. Since $\Omega$ is smooth, up to a rotation, there exists $\delta>0$ and $\varphi_0: B_\delta^{(n-1)}(0)\to \rr$ such that $\varphi_0(0)=|\nabla\varphi_0(0)|=0$ and 
\begin{equation}\label{def:phi}
\left\{\begin{array}{llll}
\varphi: & \tilde{B}_{3\delta} & \to &\rn\\
& (x_1,x') & \mapsto & (x_1+\varphi_0(x'), x'), 
\end{array}\right.
\end{equation}
that realizes a diffeomorphism onto its image and such that
$$\varphi(\tilde{B}_{3\delta}\cap\rnp)=\varphi(\tilde{B}_{3\delta})\cap\Omega\hbox{ and }\varphi(\tilde{B}_{3\delta}\cap\partial\rnp)=\varphi(\tilde{B}_{3\delta})\cap\partial\Omega.$$
Let $\eta\in C^\infty_c(\rn)$ be such that $\eta(x)=1$ for all $x\in \tilde{B}_\delta$, $\eta(x)=0$ for all $x\not\in \tilde{B}_{2\delta}$. 

\medskip\noindent{\bf Case 1: $\gamma\leq\frac{n^2-1}{4}$.} As in Chern-Lin \cite{CL5}, for any $\epsilon>0$, we define 
$$u_\epsilon(x):=\left(\eta\epsilon^{-\frac{n-2}{2}}U(\epsilon^{-1}x)\right)\circ \varphi^{-1}(x)\hbox{ for }x\in\varphi(\tilde{B}_{2\delta})\cap\Omega\hbox{ and }0\hbox{ elsewhere.}$$
This subsection is devoting to give a Taylor expansion of $J^\Omega_{\gamma,s}(\ue)$ as $\epsilon\to 0$. In the sequel, we adopt the following notation: given $(a_\epsilon)_{\epsilon>0}\in\rr$, $\Theta_\gamma(a_\epsilon)$ denotes a quantity such that, as $\epsilon\to 0$.
$$\Theta_\gamma(a_\epsilon):=\left\{\begin{array}{ll}
o(a_\epsilon) & \hbox{ if }\gamma<\frac{n^2-1}{4}\\
O(a_\epsilon) & \hbox{ if }\gamma=\frac{n^2-1}{4}\\
\end{array}\right. 
$$

\medskip\noindent{\bf Estimate of $\int_\Omega |\nabla \ue|^2\, dx$:} \\

\noindent It follows from \eqref{bnd:U} that 
\begin{equation}\label{bnd:ue}
\hbox{$|\nabla\ue(x)|\leq C\eps^{\aps-\frac{n}{2}}|x|^{-\aps}$ for all $x\in \Omega$ and $\eps>0$. }
\end{equation}
Therefore,
\begin{equation}\nonumber
\hbox{$\int_{\varphi(\left(\tilde{B}_{3\delta}\setminus \tilde{B}_{\delta}\right)\cap\rnp)}|\nabla\ue|^2\, dx=\Theta_\gamma(\epsilon)$ as $\eps\to 0$. }
\end{equation}
It follows  that
\begin{equation}\nonumber
\int_\Omega |\nabla \ue|^2\, dx=\int_{\tilde{B}_{\delta,+}}|\nabla (\ue\circ\varphi)|^2_{\varphi^\star\eucl}|\hbox{Jac}(\varphi)|\, dx+\Theta_\gamma(\epsilon)\quad \hbox{as $\eps\to 0$, }
\end{equation}
where $\tilde{B}_{\delta,+}:=\tilde{B}_{\delta}\cap\rnp$. The definition \eqref{def:phi} of $\varphi$ yields $\hbox{Jac}(\varphi)=1$. Moreover, for any $\theta\in (0,1)$, we have as $x\to 0$, 
$$\varphi^\star\eucl:=\left(\begin{array}{cc}
1 & \partial_j\varphi_0\\ \partial_i\varphi_0 & \delta_{ij}+\partial_i\varphi_0\partial_j\varphi_0\end{array}\right)= Id+ H+O(|x|^{1+\theta})$$
where
$$H:=\left(\begin{array}{cc}
0 & \partial_j\varphi_0\\ \partial_i\varphi_0 & 0\end{array}\right).$$
It follows that
\begin{eqnarray}
\int_\Omega |\nabla \ue|^2\, dx&=&\int_{\tilde{B}_{\delta,+}}|\nabla (\ue\circ\varphi)|^2_{\eucl}\, dx-\int_{\tilde{B}_{\delta,+}}H^{ij}\partial_i(\ue\circ\varphi)\partial_j (\ue\circ\varphi)\, dx\nonumber\\
&& +O\left(\int_{\tilde{B}_{\delta,+}}|x|^{1+\theta}|\nabla (\ue\circ\varphi)|^2\, dx\right)+\Theta_\gamma(\epsilon) \quad \hbox{as $\eps\to 0$. }\label{eq:pf:2}
\end{eqnarray}
We have that
\begin{eqnarray}
&&\int_{\tilde{B}_{\delta,+}}H^{ij}\partial_i(\ue\circ\varphi)\partial_j (\ue\circ\varphi)\, dx\nonumber\\
&&= 2\sum_{i\geq 2}\int_{\tilde{B}_{\delta,+}}H^{1i}\partial_1(\ue\circ\varphi)\partial_i (\ue\circ\varphi)\, dx\nonumber\\
&&=2\sum_{i\geq 2}\int_{\tilde{B}_{\delta,+}}\partial_i\varphi_0(x')\partial_1(\ue\circ\varphi)\partial_i (\ue\circ\varphi)\, dx\nonumber\\
&&= 2\sum_{i,j\geq 2}\int_{\tilde{B}_{\delta,+}}\partial_{ij}\varphi_0(0) (x')^j\partial_1(\ue\circ\varphi)\partial_i (\ue\circ\varphi)\, dx \nonumber\\
&&+O\left(\int_{\tilde{B}_{\delta,+}}|x|^{2}|\nabla (\ue\circ\varphi)|^2\, dx\right)\quad \hbox{as $\eps\to 0$. }\label{eq:pf:1}
\end{eqnarray}
We let $II$ be the second fundamental form at $0$ of the oriented boundary $\partial\Omega$. By definition, for any $X,Y\in T_0\partial\Omega$, we have that
$$II(X,Y):=\left(d\vec{\nu}_0(X),Y\right)_{\eucl}$$
where $\vec{\nu}: \partial\Omega\to \rn$ is the outer orthonormal vector of $\partial\Omega$. In particular, we have that $\vec{\nu}(0)=(-1,0,\cdot, 0)$. For any $i,j\geq 2$, we have that
$$II_{ij}:=II(\partial_i\varphi (0),\partial_j\varphi(0))= (\partial_i(\vec{\nu}\circ\varphi)(0),\partial_j\varphi(0))=-(\vec{\nu}(0),\partial_{ij}\varphi(0))=\partial_{ij}\varphi_0(0).$$
Plugging \eqref{eq:pf:1} in \eqref{eq:pf:2}, and using a change of variables, we get that
\begin{eqnarray}
\int_\Omega |\nabla \ue|^2\, dx&=&\int_{\tilde{B}_{\eps^{-1}\delta,+}}|\nabla U|^2\, dx-2II_{ij}\sum_{i,j\geq 2}\int_{\tilde{B}_{\eps^{-1}\delta,+}} (x')^j\partial_1U\partial_i U\, dx\nonumber\\
&& +O\left(\int_{\tilde{B}_{\delta,+}}|x|^{1+\theta}|\nabla (\ue\circ\varphi)|^2\, dx\right)+\Theta_\gamma(\epsilon)\label{eq:pf:3} \quad {\rm as}\, \,  \eps\to 0. 
\end{eqnarray}
We now choose $\theta$:\\
(i)\,\, If $\gamma<\frac{n^2-1}{4}$, then choose $0<\theta<\aps-\ams-1$;\\
 (ii) \,\, If $\gamma=\frac{n^2-1}{4}$, we take any $\theta\in (0,1)$. \\
  In both cases, we get by using \eqref{bnd:ue}, that
\begin{equation}\label{eq:pf:4}
\int_{\tilde{B}_{\delta,+}}|x|^{1+\theta}|\nabla (\ue\circ\varphi)|^2\, dx=\Theta_\gamma(\epsilon)\quad \hbox{as $\eps\to 0$. }
\end{equation}
Moreover, using \eqref{bnd:U}, we have that
\begin{equation}\label{eq:pf:5}
\hbox{$\int_{\tilde{B}_{\eps^{-1}\delta,+}}|\nabla U|^2\, dx=\int_{\rnp}|\nabla U|^2\, dx+\Theta_\gamma(\eps)$ as $\eps\to 0$. }
\end{equation}
Plugging together \eqref{eq:pf:3}, \eqref{eq:pf:4}, \eqref{eq:pf:5} yields
\begin{eqnarray}\label{id:nabla}
\int_\Omega |\nabla \ue|^2\, dx&=&\int_{\rnp}|\nabla U|^2\, dx\nonumber\\
&&-2II_{ij}\sum_{i,j\geq 2}\int_{\tilde{B}_{\eps^{-1}\delta,+}} (x')^j\partial_1U\partial_i U\, dx +\Theta_\gamma(\epsilon)
\end{eqnarray}

\medskip\noindent{\bf Estimate of $\int_\Omega \frac{|\ue|^{\crits}}{|x|^s}\, dx$:}\\

\noindent Fix $\sigma\in [0,2]$. We will apply the estimates below to $\sigma=s\in [0,2)$ or to $\sigma:=2$. The first estimate in \eqref{bnd:U} yields
\begin{equation}\label{est:ue:co}
|\ue(x)|\leq C\eps^{\aps-\frac{n}{2}}d(x,\partial\Omega)|x|^{-\aps}\leq C\eps^{\aps-\frac{n}{2}}|x|^{1-\aps}
\end{equation}
for all $\eps>0$ and all $x\in\Omega$. Since $\hbox{Jac }\varphi=1$, this estimate then yields
\begin{eqnarray}
\int_\Omega \frac{|\ue|^{\crit(\sigma)}}{|x|^\sigma}\, dx&=&\int_{\varphi(\tilde{B}_{\delta,+})} \frac{|\ue|^{\crit(\sigma)}}{|x|^\sigma}\, dx +\Theta_\gamma(\eps)\nonumber\\
&=& \hbox{$\int_{\tilde{B}_{\delta,+}} \frac{|\ue\circ\varphi|^{\crit(\sigma)}}{|\varphi(x)|^\sigma}\, dx +\Theta_\gamma(\eps)$ \quad as $\eps\to 0$.} \label{eq:pf:6}
\end{eqnarray}
If $\gamma<\frac{n^2-1}{4}$ or if $\gamma=\frac{n^2-1}{4}$ and $\sigma<2$, we choose $\theta\in (0,(\aps-\ams)\frac{\crit(\sigma)}{2}-1)\cap (0,1)$. If $\gamma=\frac{n^2-1}{4}$ and $\sigma=2$, we choose any $\theta\in (0,1)$. Using the expression of $\varphi(x_1,x')$, a Taylor expansion yields
\begin{equation}
|\varphi(x)|^{-\sigma}=|x|^{-\sigma}\left(1-\frac{\sigma}{2}\frac{x_1}{|x|^2}\sum_{i,j\geq 2}\partial_{ij}\varphi_0(0)(x')^i(x')^j+O(|x|^{1+\theta})\right)\quad \hbox{as $\eps\to 0$. }\label{eq:pf:7}
\end{equation}
The choice of $\theta$ yields
\begin{equation}\label{eq:pf:8}
\hbox{$\int_{\tilde{B}_{\delta,+}} \frac{|\ue\circ\varphi|^{\crit(\sigma)}}{|\varphi(x)|^\sigma}|x|^{1+\theta}\, dx=\Theta_\gamma(\eps)$ \quad as $\eps\to 0$. }
\end{equation}
Plugging together \eqref{eq:pf:6}, \eqref{eq:pf:7}, \eqref{eq:pf:8}, and using a change of variable, we get as $\eps\to 0$ that
\begin{eqnarray}
\int_\Omega \frac{|\ue|^{\crit(\sigma)}}{|x|^\sigma}\, dx&=& \int_{\tilde{B}_{\eps^{-1}\delta,+}} \frac{|U|^{\crit(\sigma)}}{|x|^\sigma}\, dx\nonumber\\
&& -\frac{\sigma}{2}\sum_{i,j\geq 2}\eps II_{ij}\int_{\tilde{B}_{\eps^{-1}\delta,+}} \frac{|U|^{\crit(\sigma)}}{|x|^\sigma}\frac{x_1}{|x|^2}(x')^i(x')^j\, dx+\Theta_\gamma(\eps).\nonumber
\end{eqnarray}
Moreover, \eqref{bnd:U} yields
$$\int_{\tilde{B}_{\eps^{-1}\delta,+}} \frac{|U|^{\crit(\sigma)}}{|x|^\sigma}\, dx=\int_{\rnp} \frac{|U|^{\crit(\sigma)}}{|x|^\sigma}\, dx+\Theta_\gamma(\eps) \quad \hbox{as $\eps\to 0$. }$$
Therefore, we get that
\begin{eqnarray}\label{id:s}
\int_\Omega \frac{|\ue|^{\crit(\sigma)}}{|x|^\sigma}\, dx&=& \int_{\rnp} \frac{|U|^{\crit(\sigma)}}{|x|^\sigma}\, dx\nonumber\\
&& -\frac{\sigma}{2}\sum_{i,j\geq 2}\eps II_{ij}\int_{\tilde{B}_{\eps^{-1}\delta,+}} \frac{|U|^{\crit(\sigma)}}{|x|^\sigma}\frac{x_1}{|x|^2}(x')^i(x')^j\, dx+\Theta_\gamma(\eps).
\end{eqnarray}
We now compute the terms in $U$ by using its symmetry property established in Chern-Lin \cite{CL5} (see also Theorem \ref{th:sym} in the Appendix). Indeed, it follows from \eqref{eq:U} that there exists $\tilde{U}: (0,+\infty)\times \rr$ such that $U(x_1,x')=\tilde{U}(x_1, |x'|)$ for all $(x_1,x')\in\rnp$. Therefore, for any $i,j\geq 2$, we get that
\begin{equation}\nonumber
\int_{\tilde{B}_{\eps^{-1}\delta,+}} \frac{|U|^{\crit(\sigma)}}{|x|^\sigma}\frac{x_1}{|x|^2}(x')^i(x')^j\, dx= \frac{\delta_{ij}}{n-1}\int_{\tilde{B}_{\eps^{-1}\delta,+}} \frac{|U|^{\crit(\sigma)}}{|x|^\sigma}\frac{x_1}{|x|^2}|x'|^2\, dx
\end{equation}
and that
\begin{equation}\nonumber
\int_{\tilde{B}_{\eps^{-1}\delta,+}} (x')^j\partial_1U\partial_i U\, dx=\frac{\delta_{ij}}{n-1}\int_{\tilde{B}_{\eps^{-1}\delta,+}} \partial_1U (x',\nabla U)\, dx
\end{equation}
where $x=(x_1,x')\in\rnp$. Therefore, the identities \eqref{id:nabla} and \eqref{id:s} rewrite
\begin{eqnarray}\label{id:nabla:bis}
\qquad \int_\Omega |\nabla \ue|^2\, dx&=&\int_{\rnp}|\nabla U|^2\, dx-\frac{2h_\Omega(0)}{n-1}\eps\int_{\tilde{B}_{\eps^{-1}\delta,+}} \partial_1U(x',\nabla U)\, dx +\Theta_\gamma(\epsilon)
\end{eqnarray}
and 
\begin{eqnarray}\label{id:s:bis}
\int_\Omega \frac{|\ue|^{\crit(\sigma)}}{|x|^\sigma}\, dx&=& \int_{\rnp} \frac{|U|^{\crit(\sigma)}}{|x|^\sigma}\, dx\\
&& -\frac{\sigma h_\Omega(0)}{2(n-1)}\eps\int_{\tilde{B}_{\eps^{-1}\delta,+}} \frac{|U|^{\crit(\sigma)}}{|x|^\sigma}\frac{x_1}{|x|^2}|x'|^2\, dx+\Theta_\gamma(\eps)\nonumber
\end{eqnarray}
as $\eps\to 0$, where $h_\Omega(0)=\sum_i II_{ii}$ is the mean curvature at $0$.\\

\medskip\noindent{\bf An intermediate identity.} We now claim that as $\eps\to 0$, 
\begin{eqnarray}
\int_{\tilde{B}_{\eps^{-1}\delta,+}} \partial_1U(x',\nabla U)\, dx &=& \int_{\tilde{B}_{\eps^{-1}\delta,+}} \frac{|x'|^2x_1}{2|x|^{2}}\left(\lambda \frac{s}{\crits}\frac{U^{\crits}}{|x|^s}+\gamma \frac{U^{2}}{|x|^2}\right)\, dx\nonumber\\
&&-\int_{\partial\rnp\cap \tilde{B}_{\eps^{-1}\delta}} \frac{|x'|^2(\partial_1 U)^2}{4}\, dx + \Theta_\gamma(1)\label{id:chern:lin}
\end{eqnarray}
where $\lambda>0$ is as in \eqref{eq:U}. This was shown by Chern-Lin \cite{CL5}, and we include it for the sake of completeness. Here and in the sequel, $\nu_i$ denotes the $i^{th}$ coordinate of the direct outward normal vector on the boundary of the relevant domain (for instance, on $\partial\rnp$, we have that $\nu_i=-\delta_{1i}$). We write
\begin{eqnarray}
&&\int_{\tilde{B}_{\eps^{-1}\delta,+}} \partial_1U(x',\nabla U)\, dx = \sum_{j\geq 2}\int_{\tilde{B}_{\eps^{-1}\delta,+}} \partial_1U (x')^j\partial_j U\, dx\nonumber\\
&&= \sum_{j\geq 2}\int_{\tilde{B}_{\eps^{-1}\delta,+}} \partial_1U \partial_j\left(\frac{|x'|^2}{2}\right)\partial_j U\, dx\nonumber\\
&&= \sum_{j\geq 2}\int_{\partial(\tilde{B}_{\eps^{-1}\delta,+})} \partial_1U \frac{|x'|^2}{2}\partial_j U\nu_j\, d\sigma-\sum_{j\geq 2}\int_{\tilde{B}_{\eps^{-1}\delta,+}} \frac{|x'|^2}{2}\partial_j \left(\partial_1U \partial_j U\right)\, dx\nonumber\\
&&= \sum_{j\geq 2}\int_{\partial\rnp\cap \tilde{B}_{\eps^{-1}\delta}} \partial_1U \frac{|x'|^2}{2}\partial_j U\nu_j\, d\sigma+ O\left(\int_{\rnp\cap \partial\tilde{B}_{\eps^{-1}\delta}} |x'|^2|\nabla U|^2(x)\, d\sigma\right)\nonumber\\
&&-\sum_{j\geq 2}\int_{\tilde{B}_{\eps^{-1}\delta,+}} \frac{|x'|^2}{2}\left(\partial_{1j}U \partial_j U+\partial_1 U\partial_{jj}U\right)\, dx\label{eq:pf:10}
\end{eqnarray}
Since $U(0,x')=0$ for all $x'\in\rr^{n-1}$, using the upper-bound \eqref{bnd:U} and writing $\nabla^\prime=(\partial_2,\dots,\partial_n)$,  we get that
\begin{eqnarray}
&&\int_{\tilde{B}_{\eps^{-1}\delta,+}} \partial_1U(x',\nabla U)\, dx \nonumber\\
&&= -\sum_{j\geq 2}\int_{\tilde{B}_{\eps^{-1}\delta,+}} \frac{|x'|^2}{2}\left(\partial_{1j}U \partial_j U+\partial_1 U\partial_{jj}U\right)\, dx+ \Theta_\gamma(1)\nonumber\\
&&= -\int_{\tilde{B}_{\eps^{-1}\delta,+}} \frac{|x'|^2}{4}\partial_1\left(|\nabla' U|^2\right)\, dx\nonumber\\
&&+ \int_{\tilde{B}_{\eps^{-1}\delta,+}} \frac{|x'|^2}{2}\partial_1U \left(-\Delta U+\partial_{11}U\right)\, dx+ \Theta_\gamma(1)\nonumber\\
&&= -\int_{\partial(\tilde{B}_{\eps^{-1}\delta,+})} \frac{|x'|^2|\nabla' U|^2}{4}\nu_1\, dx+ \int_{\tilde{B}_{\eps^{-1}\delta,+}} \frac{|x'|^2}{2}\partial_1U (-\Delta U)\, dx\nonumber\\
&&+\int_{\tilde{B}_{\eps^{-1}\delta,+}} \partial_1\left(\frac{|x'|^2(\partial_1 U)^2}{4}\right)\, dx +\Theta_\gamma(1)\label{eq:pf:11}
\end{eqnarray}
Using again that $U$ vanishes on $\partial \rnp$ and the bound \eqref{bnd:U}, we get that
\begin{eqnarray}
\int_{\tilde{B}_{\eps^{-1}\delta,+}} \partial_1U(x',\nabla U)\, dx &=& \int_{\tilde{B}_{\eps^{-1}\delta,+}} \frac{|x'|^2}{2}\partial_1U (-\Delta U)\, dx
+\int_{\partial\rnp\cap \tilde{B}_{\eps^{-1}\delta}} \frac{|x'|^2(\partial_1 U)^2}{4}\nu_1\, dx \nonumber\\&& +O\left(\int_{\partial(\tilde{B}_{\eps^{-1}\delta})\cap\rnp} |x'|^2|\nabla U|^2\, dx\right)+ \Theta_\gamma(1)\nonumber\\
&&=\int_{\tilde{B}_{\eps^{-1}\delta,+}} \frac{|x'|^2}{2}\partial_1U (-\Delta U)\, dx\nonumber\\
&&-\int_{\partial\rnp\cap \tilde{B}_{\eps^{-1}\delta}} \frac{|x'|^2(\partial_1 U)^2}{4}\, dx + \Theta_\gamma(1)\label{eq:pf:12}
\end{eqnarray}
as $\eps\to 0$. Now use equation \eqref{eq:U} to get that
\begin{equation}\label{eq:proofread}
\int_{\tilde{B}_{\eps^{-1}\delta,+}} \frac{|x'|^2}{2}\partial_1U (-\Delta U)\, dx=\int_{\tilde{B}_{\eps^{-1}\delta,+}} \frac{|x'|^2}{2}\partial_1U\left(\lambda \frac{U^{\crits-1}}{|x|^s}+\gamma \frac{U}{|x|^2}\right)\, dx.
\end{equation}
Integrating by parts, using that $U$ vanishes on $\partial\rnp$ and the upper-bound \eqref{bnd:U}, for $\sigma\in [0,2]$, 
 we get that 
 \begin{eqnarray}
&&\int_{\tilde{B}_{\eps^{-1}\delta,+}} |x'|^2\partial_1U \frac{U^{\crit(\sigma)-1}}{|x|^\sigma}\, dx=\int_{\tilde{B}_{\eps^{-1}\delta,+}} |x'|^2|x|^{-\sigma}\partial_1\left(\frac{U^{\crit(\sigma)}}{\crit(\sigma)}\right)\, dx\nonumber\\
&&=  \int_{\partial(\tilde{B}_{\eps^{-1}\delta,+})} |x'|^2|x|^{-\sigma}\frac{U^{\crit(\sigma)}}{\crit(\sigma)}\nu_1\, dx-\int_{\tilde{B}_{\eps^{-1}\delta,+}} \partial_1(|x'|^2|x|^{-\sigma})\left(\frac{U^{\crit(\sigma)}}{\crit(\sigma)}\right)\, dx\nonumber\\
&&= O\left(\int_{\rnp\cap \partial\tilde{B}_{\eps^{-1}\delta,+}} |x|^{2-\sigma}U^{\crit(\sigma)}\, d\sigma\right)+\frac{\sigma}{\crits}\int_{\tilde{B}_{\eps^{-1}\delta,+}} \frac{|x'|^2x_1}{|x|^{\sigma+2}}U^{\crit(\sigma)}\, dx\nonumber\\
&&= \frac{\sigma}{\crits}\int_{\tilde{B}_{\eps^{-1}\delta,+}} \frac{|x'|^2x_1}{|x|^{\sigma+2}}U^{\crit(\sigma)}\, dx+\Theta_\gamma(1)\quad \hbox{as $\eps\to 0$.} \label{eqref:pf:13}
\end{eqnarray}
Putting together \eqref{eq:pf:12} to \eqref{eqref:pf:13} yields \eqref{id:chern:lin}. \\

\medskip\noindent{\bf Estimate of $J^\Omega_{\gamma,s}(\ue)$:} 

Since $U\in \dundeuxr$, It follows from \eqref{eq:U} that
$$\int_{\rnp}\left(|\nabla U|^2-\frac{\gamma}{|x|^2}U^2\right)\, dx=\lambda \int_{\rnp}\frac{U^{\crits}}{|x|^s}\, dx.$$
This equality, combined with \eqref{id:nabla:bis} and \eqref{id:s:bis} gives
\begin{eqnarray}
J^\Omega_{\gamma,s}(\ue)&=& \frac{\int_\Omega\left(|\nabla \ue|^2-\frac{\gamma}{|x|^2}\ue^2\right)\, dx}{\left(\int_\Omega\frac{|\ue|^{\crits}}{|x|^s}\, dx\right)^{\frac{2}{\crits}}}\nonumber\\
&=& \frac{\int_{\rnp}\left(|\nabla U|^2-\frac{\gamma}{|x|^2}U^2\right)\, dx}{\left(\int_{\rnp}\frac{|U|^{\crits}}{|x|^s}\, dx\right)^{\frac{2}{\crits}}}\left(1+\epsilon \frac{h_\Omega(0)}{(n-1)\lambda \int_{\rnp}\frac{|U|^{\crits}}{|x|^s}\, dx} C_\eps+\Theta_\gamma(\eps)\right)\label{dev:J:1}
\end{eqnarray}
where for all $\eps>0$, 
\begin{eqnarray}
C_\eps&:=& -2\int_{\tilde{B}_{\eps^{-1}\delta,+}} \partial_1U(x',\nabla U)\, dx +\gamma \int_{\tilde{B}_{\eps^{-1}\delta,+}} \frac{|x'|^2x_1}{|x|^{2}}\frac{U^{2}}{|x|^2}\, dx\nonumber\\
&&+\lambda \frac{s}{\crits}\int_{\tilde{B}_{\eps^{-1}\delta,+}} \frac{|x'|^2x_1}{|x|^{2}}\frac{U^{\crits}}{|x|^s}\, dx\nonumber
\end{eqnarray}
 The identity \eqref{id:chern:lin} then yields
$$C_\eps=\int_{\partial\rnp\cap \tilde{B}_{\eps^{-1}\delta}} \frac{|x'|^2(\partial_1 U)^2}{2}\, dx + \Theta_\gamma(1)$$
as $\eps\to 0$. Therefore, \eqref{dev:J:1} yields that as $\eps\to 0$,
\begin{eqnarray}\label{dev:J:2}
\quad \quad J^\Omega_{\gamma, s}(\ue)&=& \mu_{\gamma,s}(\rnp)\left(1+\epsilon \frac{h_\Omega(0)\int_{\partial\rnp\cap \tilde{B}_{\eps^{-1}\delta}} |x'|^2(\partial_1 U)^2\, dx' }{2(n-1)\lambda \int_{\rnp}\frac{|U|^{\crits}}{|x|^s}\, dx} +\Theta_\gamma(\eps)\right).
\end{eqnarray}

We now distinguish two cases:

\medskip\noindent{\bf Case 1': $\gamma<\frac{n^2-1}{4}$.} The bound \eqref{bnd:U} yields $|x'|^2|\partial_1 U|^2=O(|x'|^{2-2\aps})$ when $|x'|\to +\infty$. Since $\partial\rnp=\rr^{n-1}$, we then get that $x'\mapsto |x'|^2|\partial_1 U(x')|^2$ is in $L^1(\partial\rnp)$, and therefore, \eqref{dev:J:2} yields
\begin{eqnarray}\label{dev:J:3}
\hbox{$J^\Omega_{\gamma,s}(\ue)=  \mu_{\gamma,s}(\rnp)\left(1+C_0\cdot h_\Omega(0)\cdot\epsilon  +o(\eps)\right)$ as $\eps\to 0$, }
\end{eqnarray}
with
$$C_0:=\frac{\int_{\partial\rnp} |x'|^2(\partial_1 U)^2\, dx'} {2(n-1)\lambda \int_{\rnp}\frac{|U|^{\crits}}{|x|^s}\, dx}>0.$$ 

\medskip\noindent{\bf Case 1'': $\gamma=\frac{n^2-1}{4}$.} It follows from \eqref{asymp:U}, Lemma \ref{lem:deriv:1} and a Kelvin transform that
$$\lim_{|x'|\to +\infty}|x'|^{\aps}|\partial_1 U(0,x')|=K_2>0.$$
Since $2\aps-2=n-1$, we get that
$$\int_{\partial\rnp\cap \tilde{B}_{\eps^{-1}\delta}} |x'|^2(\partial_1 U)^2\, dx'=\omega_{n-1}K_2^2\ln\frac{1}{\eps}+o\left(\ln\frac{1}{\eps}\right)$$
as $\eps\to 0$. Therefore, \eqref{dev:J:2} yields
\begin{eqnarray}\label{dev:J:4}
\hbox{$J^\Omega_{\gamma,s}(\ue)=  \mu_{\gamma,s}(\rnp)\left(1+ C_0' h_\Omega(0)\eps \ln\frac{1}{\eps}+o\left(\ln\frac{1}{\eps}\right)\right)$ \quad as $\eps\to 0$, }
\end{eqnarray}
where
$$C_0':=\frac{\omega_{n-1}K_2^2 }{2(n-1)\lambda \int_{\rnp}\frac{|U|^{\crits}}{|x|^s}\, dx}>0.$$

\smallskip\noindent Cases 1 and 2 prove Proposition \ref{prop:test:fct} when $\gamma\leq \frac{n^2-1}{4}$.

\medskip\noindent{\bf Case 2: $\gamma>\frac{n^2-1}{4}$. } In this case, the test-functions are more subtle. First,  use Theorem \ref{def:mass} to obtain $H\in C^2(\overline{\Omega}\setminus\{0\})$ such that \eqref{carac:H} hold and 
\begin{equation}\label{exp:H}
\hbox{$H(x)=\frac{d(x,\partial\Omega)}{|x|^{\aps}}+m_\gamma(\Omega)\frac{d(x,\partial\Omega)}{|x|^{\ams}}+ o\left(\frac{d(x,\partial\Omega)}{|x|^{\ams}}\right)$ \quad when $x\to 0$. }
\end{equation}
As above, we fix $\eta\in C^\infty_c(\rn)$ such that $\eta(x)=1$ for all $x\in \tilde{B}_\delta$, $\eta(x)=0$ for all $x\not\in \tilde{B}_{2\delta}$. We then define $\beta$ such that
$$H(x)=\left(\eta\frac{x_1}{|x|^{\aps}}\right)\circ \varphi^{-1}(x)+\beta(x)\quad \hbox{for all $x\in \Omega$. }$$
Here $\varphi$ is as in \eqref{def:chart:1} to \eqref{def:chart:6}. Note that $\beta\in \dundeux$ and
\begin{equation}\label{asymp:beta}
\hbox{$\beta(x)=m_\gamma(\Omega)\frac{d(x,\partial\Omega)}{|x|^{\ams}}+ o\left(\frac{d(x,\partial\Omega)}{|x|^{\ams}}\right)$\quad as $x\to 0$. }
\end{equation}
Indeed, since $\aps-\ams<1$, an essential point underlying all this subsection is that
$$\hbox{$|x|=o\left(|x|^{\aps-\ams}\right)$\quad as $x\to 0$. }$$
We choose $U$ as in \eqref{eq:U}. Up to multiplication by a constant, we assume that $K_2=1$, that is
\begin{equation}\label{asymp:U:bis}
U(x)\sim_{x\to 0}K_1\frac{x_1}{|x|^{\ams}}\hbox{ and }U(x)\sim_{|x|\to +\infty}\frac{x_1}{|x|^{\aps}}.
\end{equation}
 Now define
\begin{equation}\label{def:ue:mass}
u_\epsilon(x):=\left(\eta\epsilon^{-\frac{n-2}{2}}U(\epsilon^{-1}\cdot)\right)\circ \varphi^{-1}(x)+\epsilon^{\frac{\aps-\ams}{2}}\beta(x)\hbox{ for }x\in\Omega\hbox{ and }\epsilon>0.
\end{equation}
We start by showing that for any $k\geq 0$
\begin{equation}\label{lim:ue:H}
\lim_{\eps\to 0}\frac{\ue}{\eps^{\frac{\aps-\ams}{2}}}=H\hbox{ in }C^k_{loc}(\overline{\Omega}\setminus\{0\}).
\end{equation}
Indeed, the convergence in $C^0_{loc}(\overline{\Omega}\setminus\{0\})$ is a consequence of the definition of $\ue$, the choice $K_2=1$ and the asymptotic behavior \eqref{asymp:U:bis}. For convergence in $C^k$, we need in addition that $\nabla^i(U-x_1|x|^{-\aps})=o\left(|x|^{1-\aps-i}\right)$ as $x\to +\infty$ for all $i\geq 0$. This estimate follows from \eqref{asymp:U:bis} and Lemma \ref{lem:deriv:1}.

\smallskip\noindent  In the sequel, we adopt the following notation: $\theta_c^\eps$ denote any quantity such that there exists $\theta:\rr\to\rr$ such that
$$\lim_{c\to 0}\lim_{\eps\to 0}\theta_c^\eps=0.$$
\medskip\noindent  We first claim that for any $c>0$, we have that
\begin{eqnarray}
&&\int_{\Omega\setminus \varphi(B_c(0)_+)}\left(|\nabla \ue|^2-\frac{\gamma}{|x|^2}\ue^2\right)\,dx \nonumber\\
&&=\eps^{\aps-\ams}\left((\aps-1)c^{n-2\aps}\frac{\omega_{n-1}}{2n}+m_\gamma(\Omega)\frac{(n-2)\omega_{n-1}}{2n}\right)+\theta_c^\eps\eps^{\aps-\ams}.\label{eq:mass:1}
\end{eqnarray}
Indeed, it follows from \eqref{lim:ue:H} that
\begin{equation}\label{eq:mass:6}
\lim_{\eps\to 0}\frac{\int_{\Omega\setminus \varphi(B_c(0)_+)}\left(|\nabla \ue|^2-\frac{\gamma}{|x|^2}\ue^2\right)\, dx}{\eps^{\aps-\ams}}=\int_{\Omega\setminus \varphi(B_c(0)_+)}\left(|\nabla H|^2-\frac{\gamma}{|x|^2}H^2\right)\, dx.
\end{equation}
Since $H$ vanishes on $\partial\Omega\setminus\{0\}$ and satisfies $-\Delta H-\frac{\gamma}{|x|^2}H=0$, integrating by parts yields
\begin{eqnarray}
\int_{\Omega\setminus \varphi(B_c(0)_+)}\left(|\nabla H|^2-\frac{\gamma}{|x|^2}H^2\right)\, dx&=&-\int_{\varphi(\rnp\cap \partial B_c(0))}H\partial_\nu H\, d\sigma\nonumber\\
&=&-\int_{\rnp\cap \partial B_c(0)}H\circ\varphi\,\partial_{\varphi_\star\nu} (H\circ\varphi)\, d(\varphi^\star\sigma), \label{eq:mass:5}
\end{eqnarray} 
where in the two last inequalities, $\nu(x)$ is the outer normal vector of $B_c(0)$ at $x\in\partial B_c(0)$. 

\medskip\noindent We now estimate $H\circ\varphi\,\partial_{\varphi_\star\nu} H\circ\varphi$. It follows from \eqref{exp:H} that
\begin{equation}\label{eq:mass:4}
H\circ\varphi(x)=\frac{x_1}{|x|^{\aps}}+m_\gamma(\Omega)\frac{x_1}{|x|^{\ams}}+o\left(\frac{x_1}{|x|^{\ams}}\right)\quad \hbox{as $x\to 0$.}
\end{equation}
 It follows from elliptic theory and \eqref{asymp:beta} that for any $i=1,...,n$, we have that
\begin{equation}\label{est:nabla:beta}
\hbox{$\partial_i(\beta\circ\varphi)= \partial_i\left(m_\gamma(\Omega)\frac{x_1}{|x|^{\ams}}\right)+o\left(|x|^{-\ams}\right)$\quad as $x\to 0$. }
\end{equation}
Therefore,
\begin{eqnarray}
\partial_i(H\circ\varphi)&=&\delta_{i1}|x|^{-\aps}-\aps x_1x_i |x|^{-\aps-2}\nonumber\\
&&+m_\gamma(\Omega)\left(\delta_{i1}|x|^{-\ams}-\ams x_1x_i |x|^{-\ams-2}\right) +o\left(|x|^{-\ams}\right)\label{eq:mass:2}
\end{eqnarray}
as $x\to 0$. Moreover, $\varphi_\star\nu(x)=\frac{x}{|x|}+O(|x|)$ as $x\to 0$. Therefore, the estimate \eqref{eq:mass:2} yields
\begin{equation}\label{eq:mass:3}
\partial_{\varphi_\star\nu}(H\circ\varphi)=-(\aps-1)\frac{x_1}{|x|^{\aps+1}}-(\ams-1)m_\gamma(\Omega)\frac{x_1}{|x|^{\ams+1}}+o\left(|x|^{-\ams}\right)
\end{equation}
as $x\to 0$. By using that $\aps+\ams=n$ and $\aps-\ams<1$, \eqref{eq:mass:4} and \eqref{eq:mass:3} yield
\begin{eqnarray}\nonumber
\hbox{$-H\circ\varphi\partial_{\varphi_\star\nu}(H\circ\varphi)=\frac{(\aps-1)x_1^2}{|x|^{2\aps+1}}+(n-2)m_\gamma(\Omega)\frac{x_1^2}{|x|^{n+1}}+o\left(|x|^{1-n}\right)$\quad as $x\to 0$. }
\end{eqnarray}
Integrating this expression on $B_c(0)_+=\rnp\cap \partial B_c(0)$ and plugging into \eqref{eq:mass:5} yield
\begin{eqnarray}
\int_{\Omega\setminus \varphi(B_c(0)_+)}\left(|\nabla H|^2-\frac{\gamma}{|x|^2}H^2\right)\, dx&=&\frac{(\aps-1)c^{n-2\aps}\omega_{n-1}}{2n}+(n-2)m_\gamma(\Omega)\frac{\omega_{n-1}}{2n}+\theta_c\nonumber
\end{eqnarray} 
where $\lim_{c\to 0}\theta_c=0$. Here, we have used that 
$$\int_{\mathbb{S}^{n-1}_+}x_1^2\, d\sigma=\frac{1}{2}\int_{\mathbb{S}^{n-1}}x_1^2\, d\sigma=\frac{1}{2n}\int_{\mathbb{S}^{n-1}}|x|^2\, d\sigma=\frac{\omega_{n-1}}{2n},\; \omega_{n-1}:=\int_{\mathbb{S}^{n-1}}d\sigma.$$
This equality and \eqref{eq:mass:6} prove \eqref{eq:mass:1}. 

\medskip\noindent  We now claim that
 \begin{eqnarray}
\int_{\Omega}\left(|\nabla \ue|^2-\frac{\gamma}{|x|^2}\ue^2\right)\, dx&=&\lambda \int_{\rnp}\frac{U^{\crits}}{|x|^s}\, dx\nonumber\\
&&\hbox{$+m_\gamma(\Omega)\frac{(n-2)\omega_{n-1}}{2n}\eps^{\aps-\ams}+o\left(\eps^{\aps-\ams}\right)$\quad as $\eps\to 0$.}\label{eq:mass:7}
\end{eqnarray}
Indeed, define $U_\eps(x):=\eps^{-\frac{n-2}{2}}U(\eps^{-1}x)$ for all $x\in\rnp$. The definition \eqref{def:ue:mass} of $\ue$ rewrites as:
\begin{equation}\nonumber
\hbox{$\ue\circ\varphi(x)=U_\epsilon(x)+\eps^{\frac{\aps-\ams}{2}}\beta\circ\varphi(x)$\quad for all $x\in \rnp\cap \tilde{B}_\delta$. }
\end{equation}
Fix $c\in (0,\delta)$ that we will eventually let go to $0$. Since $d\varphi_0$ is an isometry, we get that
\begin{eqnarray}
&&\qquad \qquad \int_{\varphi(B_c(0)_+)}\left(|\nabla \ue|^2-\frac{\gamma}{|x|^2}\ue^2\right)\, dx\\
&&= \int_{B_c(0)_+}\left(|\nabla (\ue\circ\varphi)|_{\varphi^\star\eucl}^2-\frac{\gamma}{|\varphi(x)|^2}(\ue\circ\varphi)^2\right)|\hbox{Jac}(\varphi)|\, dx\nonumber\\
&&= \int_{B_c(0)_+}\left(|\nabla U_\eps|_{\varphi^\star\eucl}^2-\frac{\gamma}{|\varphi(x)|^2}U_\eps^2\right)|\hbox{Jac}(\varphi)|\, dx\nonumber\\
&&\qquad +2\eps^{\frac{\aps-\ams}{2}} \int_{B_c(0)_+}\left((\nabla U_\eps,\nabla (\beta\circ\varphi))_{\varphi^\star\eucl}-\frac{\gamma}{|\varphi(x)|^2}U_\eps(\ue\circ\varphi)\right)|\hbox{Jac}(\varphi)|\, dx\nonumber\\
&&\qquad +\eps^{\aps-\ams} \int_{B_c(0)_+}\left(|\nabla (\beta\circ\varphi)|_{\varphi^\star\eucl}^2-\frac{\gamma}{|\varphi(x)|^2}(\beta\circ\varphi)^2\right)|\hbox{Jac}(\varphi)|\, dx\label{eq:mass:8}\nonumber
\end{eqnarray}
Since $\varphi^\star\eucl=\eucl+O(|x|)$, $|\varphi(x)|=|x|+O(|x|^2)$ and $\beta\in \dundeux$, we get that
\begin{eqnarray}
&&\int_{\varphi(B_c(0)_+)}\left(|\nabla \ue|^2-\frac{\gamma}{|x|^2}\ue^2\right)\, dx= \int_{B_c(0)_+}\left(|\nabla U_\eps|_{\eucl}^2-\frac{\gamma}{|x|^2}U_\eps^2\right)|\, dx\\
&&+O\left(\int_{B_c(0)_+}|x|\left(|\nabla U_\eps|_{\eucl}^2+\frac{U_\eps^2}{|x|^2}\right)|\, dx\right)\nonumber\\
&&+2\eps^{\frac{\aps-\ams}{2}} \int_{B_c(0)_+}\left((\nabla U_\eps,\nabla (\beta\circ\varphi))_{\eucl}-\frac{\gamma}{|x|^2}U_\eps(\beta\circ\varphi)\right)\, dx\nonumber\\
&&+O\left(\eps^{\frac{\aps-\ams}{2}} \int_{B_c(0)_+}|x|\left(|\nabla U_\eps|\cdot |\nabla (\beta\circ\varphi)|+\frac{U_\eps |\beta\circ\varphi|}{|x|^2}\right)\, dx\right)+\eps^{\aps-\ams}\theta_c^\eps \nonumber
\end{eqnarray}
as $\eps\to 0$. The pointwise estimates \eqref{asymp:U:bis} and \eqref{est:nabla:beta} yield
\begin{eqnarray}
&&\int_{\varphi(B_c(0)_+)}\left(|\nabla \ue|^2-\frac{\gamma}{|x|^2}\ue^2\right)\, dx= \int_{B_c(0)_+}\left(|\nabla U_\eps|_{\eucl}^2-\frac{\gamma}{|x|^2}U_\eps^2\right)\, dx\nonumber\\
&&+2\eps^{\frac{\aps-\ams}{2}} \int_{B_c(0)_+}\left((\nabla U_\eps,\nabla (\beta\circ\varphi))_{\eucl}-\frac{\gamma}{|x|^2}U_\eps(\beta\circ\varphi)\right)\, dx\nonumber\\
&&+\eps^{\aps-\ams}\theta_c^\eps\nonumber
\end{eqnarray}
as $\eps\to 0$. Integrating by parts yields
\begin{eqnarray}
&&\int_{\varphi(B_c(0)_+)}\left(|\nabla \ue|^2-\frac{\gamma}{|x|^2}\ue^2\right)\, dx\nonumber\\
&&= \int_{B_c(0)_+}\left(-\Delta U_\eps-\frac{\gamma}{|x|^2}U_\eps\right)U_\eps\, dx+\int_{\partial (B_c(0)_+)} U_\eps\partial_\nu U_\eps \, d\sigma\nonumber\\
&&+2\eps^{\frac{\aps-\ams}{2}} \left(\int_{B_c(0)_+}\left(-\Delta U_\eps-\frac{\gamma}{|x|^2}U_\eps\right)\beta\circ\varphi\, dx +\int_{\partial (B_c(0)_+)} \beta\circ\varphi\partial_\nu U_\eps \, d\sigma\right)\nonumber\\
&&+\eps^{\aps-\ams}\theta_c^\eps\nonumber
\end{eqnarray}
as $\eps\to 0$. Since both $U$ and $\beta\circ\varphi$ vanish on $\partial\rnp\setminus\{0\}$, we get that
\begin{eqnarray}
&&\qquad \int_{\varphi(B_c(0)_+)}\left(|\nabla \ue|^2-\frac{\gamma}{|x|^2}\ue^2\right)\, dx\\
&&= \int_{B_c(0)_+}\left(-\Delta U_\eps-\frac{\gamma}{|x|^2}U_\eps\right)U_\eps\, dx+\int_{\rnp\cap \partial B_c(0)} U_\eps\partial_\nu U_\eps \, d\sigma\nonumber\\
&&+2\eps^{\frac{\aps-\ams}{2}} \left(\int_{B_c(0)_+}\left(-\Delta U_\eps-\frac{\gamma}{|x|^2}U_\eps\right)\beta\circ\varphi\, dx +\int_{\rnp\cap \partial B_c(0)} \beta\circ\varphi\partial_\nu U_\eps \, d\sigma\right)\nonumber\\
&&+\eps^{\aps-\ams}\theta_c^\eps\label{eq:mass:9}\nonumber
\end{eqnarray}
as $\eps\to 0$. 
The asymptotic estimate \eqref{asymp:U:bis} of $U$ and Lemma \ref{lem:deriv:1} yield (after a Kelvin transform)
$$\partial_\nu U_\eps=-(\aps-1)\eps^{\frac{\aps-\ams}{2}}x_1|x|^{-\aps-1}+o\left(\eps^{\frac{\aps-\ams}{2}}|x|^{-\aps}\right)$$
as $\eps\to 0$ uniformly on compact subsets of $\overline{\rnp}\setminus\{0\}$. We then get that
\begin{equation}\nonumber
\beta\circ\varphi\partial_\nu U_\eps=\eps^{\frac{\aps-\ams}{2}}\left(-m_\gamma(\Omega)(\aps-1)x_1^2|x|^{-n-1}+o\left(|x|^{1-n}\right)\right)
\end{equation}
and
\begin{equation}\nonumber
U_\eps\partial_\nu U_\eps=\eps^{\aps-\ams}\left(-(\aps-1)x_1^2|x|^{-2\aps-1}+o\left(|x|^{1-2\aps}\right)\right)
\end{equation}
as $\eps\to 0$ uniformly on compact subsets of $\overline{\rnp}\setminus\{0\}$. Plugging these identities in \eqref{eq:mass:9} and using equation \eqref{eq:U} yield, as $\eps\to 0$, 
\begin{eqnarray}
\int_{\varphi(B_c(0)_+)}\left(|\nabla \ue|^2-\frac{\gamma}{|x|^2}\ue^2\right)\, dx
&=& \int_{B_c(0)_+}\lambda \frac{U_\eps^{\crits}}{|x|^s}\, dx-(\aps-1)\frac{\omega_{n-1}}{2n}c^{n-2\aps}\eps^{\aps-\ams}\nonumber\\
&&+2\eps^{\frac{\aps-\ams}{2}} \int_{B_c(0)_+}\lambda \frac{U_\eps^{\crits-1}}{|x|^s}\beta\circ\varphi\, dx\nonumber\\
&& -(\aps-1)\frac{\omega_{n-1}}{n}m_\gamma(\Omega)\eps^{\aps-\ams} 
+\eps^{\aps-\ams}\theta_c^\eps\label{eq:mass:10}
\end{eqnarray}
Note that as $\eps\to 0$, 
\begin{eqnarray}
\int_{B_c(0)_+}\lambda \frac{U_\eps^{\crits}}{|x|^s}\, dx&=&\int_{\rnp}\lambda \frac{U_\eps^{\crits}}{|x|^s}\, dx+O\left(\int_{\rnp\setminus B_c(0)_+} \frac{U_\eps^{\crits}}{|x|^s}\, dx\right)\nonumber\\
&&=\int_{\rnp}\lambda \frac{U_\eps^{\crits}}{|x|^s}\, dx+o\left(\eps^{\aps-\ams}\right)\label{eq:mass:11}
\end{eqnarray}
 The expansion \eqref{asymp:beta} and the change of variable $x:=\eps y$ yield as $\eps\to 0$, 
\begin{equation}\label{eq:mass:12}
\int_{B_c(0)_+}\lambda \frac{U_\eps^{\crits-1}}{|x|^s}\beta\circ\varphi\, dx=\lambda m_\gamma(\Omega)\eps^{\frac{\aps-\ams}{2}}\int_{\rnp}\frac{U^{\crits-1}}{|y|^s}\frac{y_1}{|y|^{\ams}}\, dy+  \eps^{\frac{\aps-\ams}{2}}\theta_\eps^c
\end{equation} 
Integrating by parts, and using the asymptotics \eqref{asymp:U:bis} for $U$ yield
\begin{eqnarray}
\lambda \int_{\rnp}\frac{U^{\crits-1}}{|y|^s}\frac{y_1}{|y|^{\ams}}\, dy&=& \lim_{R\to +\infty}\int_{B_R(0)_+}\lambda\frac{U^{\crits-1}}{|y|^s}\frac{y_1}{|y|^{\ams}}\, dy\nonumber\\
&=& \lim_{R\to +\infty}\int_{B_R(0)_+}\left(-\Delta U-\frac{\gamma}{|y|^2}U\right)\frac{y_1}{|y|^{\ams}}\, dy\nonumber\\
&=& \lim_{R\to +\infty}\int_{B_R(0)_+} U\left(-\Delta -\frac{\gamma}{|y|^2}\right)\left(\frac{y_1}{|y|^{\ams}}\right)\, dy\nonumber\\
&&-\int_{\partial B_R(0)_+}\partial_\nu U \frac{y_1}{|y|^{\ams}}\, d\sigma\nonumber\\
&=& \left(\aps-1\right)\frac{\omega_{n-1}}{2n}.\label{eq:mass:13}
\end{eqnarray}
Putting together \eqref{eq:mass:11}, \eqref{eq:mass:12} and \eqref{eq:mass:13} yield
\begin{eqnarray}
\int_{\Omega}\left(|\nabla \ue|^2-\frac{\gamma}{|x|^2}\ue^2\right)\, dx&=&\lambda \int_{\rnp}\frac{U^{\crits}}{|x|^s}\, dx\nonumber\\
&& +m_\gamma(\Omega)\frac{(n-2)\omega_{n-1}}{2n}\eps^{\aps-\ams}+o\left(\eps^{\aps-\ams}\right)\nonumber
\end{eqnarray}
as $\eps\to 0$. This finally yields \eqref{eq:mass:7}.

\medskip\noindent We finally claim that
\begin{eqnarray}
\int_\Omega \frac{\ue^{\crits}}{|x|^s}\, dx&=&\int_{\rnp}\frac{U^{\crits}}{|x|^s}\, dx+\frac{\crits}{\lambda}m_\gamma(\Omega)\frac{(\aps-1)\omega_{n-1}}{2n}\eps^{\aps-\ams}\nonumber\\
&&+o\left(\eps^{\aps-\ams}\right)\hbox{\qquad as $\eps\to 0$. }\label{eq:mass:14}
\end{eqnarray}
Indeed, fix $c>0$. Due to estimates \eqref{asymp:beta} and \eqref{asymp:U:bis}, we have that
\begin{eqnarray}
\int_\Omega \frac{\ue^{\crits}}{|x|^s}\, dx&= &\int_{\varphi(B_c(0)_+)} \frac{\ue^{\crits}}{|x|^s}\, dx+o\left(\eps^{\aps-\ams}\right)\nonumber\\
&= &\int_{B_c(0)_+} \frac{|U_\eps+\eps^{\frac{\aps-\ams}{2}}\beta\circ\varphi|^{\crits}}{|\varphi(x)|^s}|\hbox{Jac}(\varphi)|\, dx+o\left(\eps^{\aps-\ams}\right)\nonumber\\
&= &\int_{B_c(0)_+} \frac{|U_\eps+\eps^{\frac{\aps-\ams}{2}}\beta\circ\varphi|^{\crits}}{|x|^s}|(1+O(|x|))\, dx+o\left(\eps^{\aps-\ams}\right)\nonumber
\end{eqnarray}
as $\eps\to 0$. One can easily check that there exists $C>0$ such that for all $X,Y\in\rr$, 
\begin{equation}\label{ineq:not:numb}
||X+Y|^{\crits}-|X|^{\crits}-\crits|X|^{\crits-2}XY|\leq C\left(|X|^{\crits-2}|Y|^2+|Y|^{\crits}\right)
\end{equation}
Therefore, using the asymptotics \eqref{asymp:beta} and \eqref{asymp:U:bis} of $U$ and $\beta$, we get that
\begin{eqnarray}
\int_\Omega \frac{\ue^{\crits}}{|x|^s}\, dx&= &\int_{B_c(0)_+} \frac{U_\eps^{\crits}}{|x|^s}|(1+O(|x|))\, dx\nonumber\\
&&+\crits\eps^{\frac{\aps-\ams}{2}}\int_{B_c(0)_+} \frac{U_\eps^{\crits-1}}{|x|^s}\beta\circ\varphi(1+O(|x|))\, dx\nonumber\\
&&+\eps^{\frac{\aps-\ams}{2}}\theta_\eps^c\nonumber\\
&= &\int_{B_c(0)_+} \frac{U_\eps^{\crits}}{|x|^s}\, dx+\crits\eps^{\frac{\aps-\ams}{2}}\int_{B_c(0)_+} \frac{U_\eps^{\crits-1}}{|x|^s}\beta\circ\varphi\, dx\nonumber\\
&&+\eps^{\frac{\aps-\ams}{2}}\theta_\eps^c \quad \hbox{as $\eps\to 0$.}\nonumber
\end{eqnarray}
 Then \eqref{eq:mass:14} follows from this latest identity, combined with \eqref{eq:mass:11}, \eqref{eq:mass:12}, and \eqref{eq:mass:13}.

\medskip\noindent We can finally use \eqref{eq:mass:7} and \eqref{eq:mass:14}, and the fact that 
$$\int_{\rnp}(|\nabla U|^2-\frac{\gamma}{|x|^2}U^2)\, dx=\lambda \int_{\rnp}\frac{U^{\crits}}{|x|^s}\, dx,$$
 to get 

\begin{equation}\nonumber
J^\Omega_{\gamma,s}(\ue)=J^{\rnp}_{\gamma,s}(U)\left(1-\frac{\left(\aps-\frac{n}{2}\right)\omega_{n-1}}{n\lambda \int_{\rnp}\frac{U^{\crits}}{|x|^s}\, dx}m_\gamma(\Omega)\eps^{\aps-\ams}+o\left(\eps^{\aps-\ams}\right)\right)\quad \hbox{as $\eps\to 0$, }
\end{equation}
which proves \eqref{asymp:J:gamma:large}. This ends the proof of Proposition \ref{prop:test:fct}, and therefore, as already mentioned, of Theorem \ref{main}.

\section{Examples of domains with positive mass}\label{sec:ex:mass}
We now assume that $\gamma\in (\frac{n^2-1}{4},\frac{n^2}{4})$. We have seen in Proposition \ref{prop:mono} that the mass is negative when $\Omega\subset \rnp=T_0\partial\Omega$. In particular, $m_\gamma(\Omega)<0$ if $\Omega$ is convex and $\gamma<\gamma_H(\Omega)$. In this section, we give examples of domains $\Omega$ with positive mass.

\medskip\noindent For any $x_0\in \rn$, we define the inversion
$$i_{x_0}(x):=x_0+|x_0|^2\frac{x-x_0}{|x-x_0|^2}$$
for all $x\in\rn\setminus\{x_0\}$. The inversion $i_{x_0}$ is the indentity map on  $\partial B_{|x_0|}(x_0)$ (the ball of center $x_0$ and of radius $|x_0|$), and in particular $i(0)=0$. 


\begin{definition} We shall say that a domain $\Omega\subset\rn$ ($0\in\partial\Omega$) is {\it smooth at infinity} if there exists $x_0\not\in \overline{\Omega}$ such that $i_{x_0}(\Omega)$  is a smooth bounded domain of $\rn$ having both $0$ and $x_0$ being on its boundary $ \partial(i_{x_0}(\Omega))$. 
\end{definition}
\noindent One can easily check that $\rnp$ is a smooth domain at infinity (take $x_0:=(-1,0,\dots,0)$). \\
 We now state and prove three propositions that are crucial for the constructions that follow. The first one indicates that the notion of mass defined in Theorem \ref{def:mass} extends to unbounded domains that are smooth at infinity.

\begin{proposition}\label{prop:pert:mass} Let $\Omega$ be a domain that is smooth at infinity such that $0\in\partial\Omega$. We assume that $\gamma_H(\Omega)>\frac{n^2-1}{4}$ and fix $\gamma\in \left(\frac{n^2-1}{4}, \gamma_H(\Omega)\right)$. Then, up to a multiplicative constant, there exists a unique function $H\in C^2(\overline{\Omega}\setminus\{0\})$ such that
\begin{equation}\label{def:H:2}
\left\{\begin{array}{ll}
-\Delta H-\frac{\gamma}{|x|^2}H=0&\hbox{\rm in }\Omega\\
\hfill H>0&\hbox{\rm  in }\Omega\\
\hfill H=0&\hbox{\rm on }\partial\Omega\setminus\{0\}\\
\hfill H(x)\leq C|x|^{1-\ap}&\hbox{\rm for all }x\in\Omega.
\end{array}\right.
\end{equation}
Moreover, there exists $c_1>0$ and $c_2\in\rr$ such that
$$H(x)=c_1\frac{d(x,\partial\Omega)}{|x|^{\ap}}+c_2\frac{d(x,\partial\Omega)}{|x|^{\am}}+o\left(\frac{d(x,\partial\Omega)}{|x|^{\am}}\right) \quad \hbox{as $x\to 0$. }$$
We define the mass $m_\gamma(\Omega):=\frac{c_2}{c_1}$, which is independent of the choice of $H$ in \eqref{def:H:2}.
\end{proposition}
\noindent With this notion of mass, we will be in position to prove the following continuity result.
\begin{proposition}\label{prop:mass:positiv}
Let $\Omega\subset \rn$ be  that is smooth at infinity such that $0\in\partial\Omega$. We assume that $\gamma_H(\Omega)>\frac{n^2-1}{4}$, and fix $\gamma\in \left(\frac{n^2-1}{4}, \gamma_H(\Omega)\right)$. For any $R>0$, let $D_R$ be a smooth domain of $\rn$ such that
\begin{itemize}
\item $B_{R}(x_0)\subset D_R\subset B_{2R}(x_0)$,
\item $\Omega\cap D_R$ is a smooth domain of $\rn$.
\end{itemize}
Let $\Phi\in C^\infty(\rr\times\rn,\rn)$ be such that
\begin{itemize}
\item $\Phi_t:=\Phi(t,\cdot)$ is a smooth diffeomorphism of $\rn$,
\item $\Phi_t(x)=x$ for all $|x|>1/2$ and all $t\in\rr$,
\item $\Phi_t(0)=0$ for all $t\in\rr$,
\item $\Phi_0=Id_{\rn}$.
\end{itemize}
Set $\Omega_{t,R}:=\Phi_t(\Omega)\cap D_R$. Then 
$$\liminf_{t\to 0,\, R\to +\infty}\gamma_H(\Omega_{t,R})\geq \gamma_H(\Omega).$$
Therefore, for $t\to 0$, $R\to +\infty$, we have that $\gamma_H(\Omega_{t,R})>\frac{n^2-1}{4}$ and $m_\gamma(\Omega_{t,R})$ is well defined. In addition, 
$$\lim_{t\to 0,\, R\to +\infty}m_\gamma(\Omega_{t,R})=m_\gamma(\Omega).$$
\end{proposition}
As a consequence of the above, we shall be able to construct smooth bounded domains with positive or negative mass with any behavior at $0$.

\begin{proposition}\label{prop:ex:any:behave} Let $\omega$ be a smooth open set of $\rn$. Then, there exist $\Omega_+,\Omega_-$ two smooth bounded domains of $\rn$ with Hardy constants $>\frac{n^2-1}{4}$, and there exists $r_0>0$  such that
$$\Omega_+\cap B_{r_0}(0)=\Omega_-\cap B_{r_0}(0)=\omega\cap B_{r_0}(0),$$
and for any $\gamma \in (\frac{n^2-1}{4},\min\{\gamma_H(\Omega_+),\gamma_H(\Omega_-)\})$, we have that
$$m_\gamma(\Omega_+)>0>m_\gamma(\Omega_-).$$
\end{proposition}

\medskip\noindent The remainder of this section will be devoted to the proof of these three propositions. As a preliminary remark, we claim that if $\Omega$ is a domain of $\rn$ such that $0\in\partial\Omega$ and $\Omega$ is smooth at infinity, then 
\begin{equation}\label{lim:gamma:t}
\liminf_{t\to 0,R\to \infty}\gamma_H(\Omega_{t,R})\geq \gamma_H(\Omega),
\end{equation}
where $\Omega_{t,R}$ are defined as in Proposition \ref{prop:mass:positiv}. Indeed, by definition, $\gamma_H(\Omega_{t,R})\geq \gamma_H(\Omega_{t})=\gamma_H(\Phi_t(\Omega))$. Inequality \eqref{lim:gamma:t} then follows from the limit \eqref{lim:gamma:h} of Lemma \ref{lem:ex:gamma}. 
We shall proceed in 7 steps.

\smallskip\noindent{\bf Step 1: Reformulation via the inversion.} For convenience, up to a rotation and a dilation, we can assume that $x_0:=(-1,0,\dots,0)\in\rn$ and we define the inversion 
$$i(x):=x_0+\frac{x-x_0}{|x-x_0|^2} \quad \hbox{for all $x\in\rn\setminus\{x_0\}$.} $$
Define $\tOmega:=i(\Omega)$, $\tPhi(t,x):=i\circ\Phi(t,i(x))$ for $(t,x)\in\rr\times\rn$, and $\tD_r:=\rn\setminus i(D_{r^{-1}})$ (i.e., the complement in $\rn$). Here, note that $R\to+\infty$ in Proposition \ref{prop:mass:positiv} is equivalent to $r\to 0$ in here. We then have that
$$0,x_0\in\partial\tOmega\hbox{ and }\tOmega\hbox{ is a smooth bounded domain of }\rn.$$
Note  that $\tPhi\in C^\infty(\rr\times\rn,\rn)$ is such that
\begin{itemize} 
\item For any $t\in (-2,2)$, $\tPhi_t:=\tPhi(t,\cdot)$ is a $C^\infty-$diffeomorphism onto its open image $\tPhi_t(\rn)$.
\item $\tPhi_0=\hbox{Id}$,
\item $\tPhi_t(0)=0$ for all $t\in (-2,2)$,
\item $\tPhi_t(x)=x$ for all $t\in (-2,2)$ and all $x\in B_{2\delta}(x_0)$ with $\delta<1/4$. 
\end{itemize}
We define
$$\tOmega_t:=\tPhi_t(\tOmega).$$
The sets $\tD_r$ satisfy the following properties:
\begin{itemize}
\item $B_{r/2}(x_0)\subset \tD_r\subset B_{r}(x_0)$,
\item $\tOmega_{t,r}:=\tOmega_t\setminus \tD_r$ is a smooth domain of $\rn$.
\end{itemize}
In particular, we have that
$$\tOmega_{t,r}=i(\Omega_{t,r^{-1}}).$$
Let $u\in C^2(\overline{\Omega_{t,r}}\setminus \{0\})$ be such that
\begin{equation}\label{carac:H:prime}
-\Delta u-\frac{\gamma}{|x|^2}u=0\hbox{ in }\Omega_{t,r}\; ,\; u>0\hbox{ in }\Omega_{t,r}\; ,\; u=0\hbox{ on }\partial\Omega_{t,r}.
\end{equation}
The existence of $u$ follows from Theorem \ref{def:mass}. Consider the Kelvin transform of $u$, that is 
$$\tilde{u}(x):=|x-x_0|^{2-n}u(i(x))\hbox{ for all }x\in\tOmega_{t,r}.$$
We then get that
$$-\Delta \tilde{u}-V\tilde{u}=0\hbox{ in }\tOmega_{t,r},$$
where
\begin{equation}\label{def:V}
V(x):=\frac{\gamma}{|x|^2|x-x_0|^2}\hbox{ for }x\in\rn\setminus\{0,x_0\}.
\end{equation}
It is easy to check that
\begin{equation}\nonumber
V(x)=\frac{\gamma+O(|x|)}{|x|^2}\hbox{ as }x\to 0\hbox{ and }V(x)=\frac{\gamma+O(|x-x_0|)}{|x-x_0|^2}\hbox{ as }x\to x_0
\end{equation}
The coercivity of $-\Delta-\gamma|x|^{-2}$ on $\Omega$ (since $\gamma<\gamma_H(\Omega)$) yields the coercivity of $-\Delta-V$ on $\tOmega$, that is there exists $c_0>0$ such that
\begin{equation}\nonumber
\int_{\tOmega}\left(|\nabla u|^2-V(x) u^2\right)\, dx\geq c_0\int_{\tOmega} |\nabla u|^2\, dx\hbox{ for all }u\in D^{1,2}(\tOmega).
\end{equation}
From now on, we should be able to transfer the analysis to the bounded domain $\tOmega$. 

\medskip\noindent{\bf Step 2: Perturbation of the domain via the two singular points $0$ and $x_0$.}\par
\noindent We shall need the following.
\begin{proposition}\label{prop:sol:pert}
For any $t\in (-1,1)$, there exists $u_t\in C^2(\overline{\tOmega_t}\setminus\{0,x_0\})$ such that
\begin{equation}\label{eq:ut}
\left\{\begin{array}{ll}
\hfill -\Delta u_t-V u_t=0 &\hbox{\rm in }\tOmega_t\\
\hfill u_t>0 & \hbox{\rm in }\tOmega_t\\
\hfill u_{t}=0 & \hbox{\rm on }\partial\tOmega_t\setminus\{0, x_0\}\\
\hfill u_t(x)\leq C|x|^{1-\ap}+C|x-x_0|^{1-\am}&\hbox{ \rm for }x\in\tOmega_t.
\end{array}\right.
\end{equation} 
Moreover, we have that
\begin{equation}\label{asymp:ut}
u_t(x)=\frac{d(x,\partial\tOmega_t)}{|x|^{\ap}}(1+O(|x|^{\ap-\am})) 
\end{equation}
as $x\to 0$, uniformly wrt $t\in (-1,1)$.
\end{proposition}

\smallskip\noindent{\it Proof of Proposition \ref{prop:sol:pert}.} We construct approximate singular solutions as in Section \ref{sec:sub:super}. For all $t\in (-2,2)$, there exists a chart $\varphi_t$ that satisfies \eqref{def:chart:1} to \eqref{def:chart:6} for $\tOmega_t$. Without restriction, we assume that $\lim_{t\to 0}\varphi_t=\varphi_0$ in $C^k(\tilde{B}_{2\delta},\rn)$. We define a cut-off function $\eta_\delta$ such that $\eta_\delta(x)=1$ for $x\in \tilde{B}_\delta$ and $\eta_\delta(x)=0$ for $x\not\in \tilde{B}_{2\delta}$. As in \eqref{def:ua}, we define  $u_{\ap, t}\in C^2(\overline{\tOmega_t}\setminus\{0\})$ with compact support in $\varphi_t(\tilde{B}_{2\delta})$ such that
\begin{equation}\label{def:ua:bis}
u_{\aps,t}\circ\varphi_t(x_1,x'):=\eta_\delta(x_1,x')x_1|x|^{-\aps}(1+\Theta_t(x))\hbox{ for all }(x_1,x')\in \tilde{B}_{2\delta}\setminus \{0\},
\end{equation}
where 
$\Theta_t(x_1,x'):=e^{-\frac{1}{2}x_1 H_t(x')}-1$ 
for all $x=(x_1,x')\in \tilde{B}_{2\delta}$ and all $t\in (-2,2)$. Here, $H_t(x')$ is the mean curvature of $\partial\tOmega_t$ at the point $\varphi_t(0,x')$. Note that $\lim_{t\to 0}\Theta_t=\Theta_0$ in $C^k(U)$. Arguing as is Section \ref{sec:sub:super}, we get that
\begin{equation}\nonumber
\left\{\begin{array}{lll}
\hfill (-\Delta-V)u_{\aps, t}&=O(d(x,\partial\tOmega_t)|x|^{-\ap-1}) &\hbox{ in }\tOmega_t\cap \tilde{B}_{\delta}\\
\hfill u_{\aps, t}&>0 & \hbox{ in }\tOmega_t\cap \tilde{B}_{\delta}\\
\hfill u_{\aps, t}&=0 & \hbox{ on }\partial\tOmega_t\setminus\{0\},
\end{array}\right.
\end{equation} 
and
\begin{equation}\nonumber
u_{\aps, t}(x)=\frac{d(x,\partial\tOmega_t)}{|x|^{\ap}}(1+O(|x|)\hbox{ as }x\to 0.
\end{equation}
The construction in Section \ref{sec:sub:super} also yields
\begin{equation}\label{cv:uap:t}
\lim_{t\to 0}u_{\aps,t}\circ\Phi_t=u_{\aps,0}\hbox{ in }C^2_{loc}(\overline{\tOmega}\setminus\{0\}).
\end{equation}
Note also that all these estimates are uniform in $t\in (-1,1)$. In particular, defining 
\begin{equation}\label{eq:ft}
f_t:=-\Delta u_{\aps, t}-Vu_{\aps, t},
\end{equation}
then there exists $C>0$ such that
\begin{equation}\label{bnd:ft}
|f_t(x)|\leq Cd(x,\partial\tOmega_t)|x|^{-\ap-1}\leq C|x|^{-\ap}
\end{equation} 
for all $t\in (-1,1)$ and all $x\in\tOmega_t\cap \tilde{B}_{\delta}$. Therefore, since $\gamma>\frac{n^2-1}{4}$, it follows from \eqref{cv:uap:t} and this pointwise control that $f_t\in L^{\frac{2n}{n+2}}(\tOmega_t)$ for all $t\in (-1,1)$ and that
\begin{equation}\label{lim:ft}
\lim_{t\to 0}\Vert f_t\circ\Phi_t-f_0\Vert_{L^{\frac{2n}{n+2}}(\tOmega)}=0.
\end{equation}
For any $t\in (-1,1)$, we let $v_t\in D^{1,2}(\tOmega_t)$ be such that
\begin{equation}\label{def:vt}
-\Delta v_t-V v_t=f_t\hbox{ weakly in }D^{1,2}(\tOmega_t). 
\end{equation}
The existence follows from the coercivity of $-\Delta-V$ on $\tOmega_t$, which follows itself from the coercivity on $\tOmega=\tOmega_0$. We then get from \eqref{lim:ft} and the uniform coercivity on $\tOmega_t$ that
\begin{equation}\nonumber
\lim_{t\to 0}v_t\circ\Phi_t=v_0\hbox{ in }D^{1,2}(\tOmega)\hbox{ and }C^1_{loc}(\overline{\tOmega}\setminus\{0,x_0\}).
\end{equation}
It follows from the construction of the mass in Section \ref{sec:mass} (see the proof of Theorem \ref{def:mass}) that around $0$, $|v_t(x)|$ is bounded by $|x|^{1-\am}$. Around $x_0$, $-\Delta v_t-V v_t=0$ and the regularity Theorem \ref{th:hopf} yields a control by $|x-x_0|^{1-\am}$. These controls are uniform with respect to $t\in (-1,1)$. Therefore, there exists $C>0$ such that
$$|v_t(x)|\leq Cd(x,\partial\tOmega_t)\left(|x|^{-\am}+|x-x_0|^{-\am}\right)$$
for all $t\in (-1,1)$ and all $x\in\tOmega_t$. Now define 
\begin{equation}\nonumber
u_t(x):=u_{\aps,t}(x)-v_t(x)
\end{equation}
for all $t\in (-1,1)$ and $x\in \tOmega_t$. This function satisfies all the requirements of Proposition \ref{prop:sol:pert}. 
\hfill$\Box$

\smallskip\noindent{\bf Step 3: Chopping off a neighborhood of $x_0$:} We now study $\tOmega_{t,r}=\tOmega_t\setminus \tD_r$. For $r\in (0,\delta/2)$, note that $\tOmega_{t,r}\cap B_\delta(0)=\tOmega\cap B_\delta(0)$. We shall now define a mass associated  to the potential $V$, and prove its continuity.

\smallskip\noindent{\bf Step 3.1: } The function $f_t:\tOmega_t\to\rr$ defined in \eqref{eq:ft} has compact support in $B_{2\delta}(0)$, therefore, it is well-defined also on $\tOmega_{t,r}$. We define $v_{t,r}\in D^{1,2}(\tOmega_{t,r})$ such that
\begin{equation}\label{eq:vtr}
-\Delta v_{t,r}-V v_{t,r}=f_{t}\hbox{ weakly in }D^{1,2}(\tOmega_{t,r}). 
\end{equation}
Since the operator $-\Delta-V$ is uniformly coercive on $\tOmega_t$, it is also uniformly coercive of $\tOmega_{t,r}$ with respect to $(t,r)$, so that the definition of $v_{t,r}$ via \eqref{eq:vtr} makes sense. The uniform coercivity and \eqref{eq:ft}-\eqref{bnd:ft} yield the existence if $C>0$ such that $\Vert v_{t,r}\Vert_{D^{1,2}(\tOmega_{t,r})}\leq C$ for all $t,r$. Since $x_0\not\in \tOmega_{t,r}$, \eqref{eq:ft}-\eqref{bnd:ft} and regularity theory yield $v_{t,r}\in C^1(\overline{\tOmega_{t,r}}\setminus\{0\})$ and for all $\rho>0$, there exists $C(\rho)>0$ such that 
\begin{equation}\label{bnd:v:1}
\Vert v_{t,r}\Vert_{C^1(\tOmega_{t,r}\setminus (B_\rho(0)\cup B_\rho(x_0)))}\leq C(\rho).
\end{equation}

\medskip\noindent{\bf Step 3.2: } We claim that there exists $C>0$ such that
\begin{equation}\label{asymp:vtr}
|v_{t,r}(x)|\leq Cd(x,\partial\Omega_{t})\left(|x|^{-\am}+|x-x_0|^{-\am}\right)
\end{equation}
for all $t\in (-1,1)$ and all $x\in\tOmega_{t,r}$. Indeed, around $0$, $\tOmega_{t,r}$ coincides with $\tOmega_{t}$, and the proof of the control goes as in the construction of the mass in Section \ref{sec:mass} (see the proof of Proposition \ref{def:mass}). The argument is different around $x_0$. We let $r_0>0$ be such that $\tOmega_{t}\cap B_{2r_0}(x_0)=\tOmega\cap B_{2r_0}(x_0)$. Therefore, for $r\in (0,r_0)$, we have that 
$$\tOmega_{t,r}\cap B_{2r_0}(x_0)=(\tOmega\setminus \tD_r)\cap B_{2r_0}(x_0).$$
Arguing as in the proof of Proposition \ref{prop:sub:super}, there exists $\tilde{u}_{\ams}\in C^\infty(\overline{\tOmega}\setminus\{0\})$ and $\tau'>0$ such that
$$\left\{\begin{array}{lll}
\hfill \tilde{u}_{\ams}&>0 &\hbox{ in }\tOmega\cap B_{2r_0}(x_0)\\
\hfill \tilde{u}_{\ams}&=0 &\hbox{ in }(\partial\tOmega)\cap B_{2r_0}(x_0)\\
\hfill -\Delta \tilde{u}_{\ams}-V\tilde{u}_{\ams}&>0&\hbox{ in }\tOmega\cap B_{2r_0}(x_0).\end{array}
\right.$$
Moreover, we have that
\begin{equation}\label{asymp:tu}
\tilde{u}_{\ams}(x)=\frac{d(x,\partial\tOmega)}{|x-x_0|^{\ams}}(1+O(|x-x_0|)) \quad\hbox{as $x\to x_0$, $x\in \tOmega$. }
\end{equation}
Therefore, since $v_{t,r}$ vanishes on $B_{2r_0}(x_0)\cap \partial(\tOmega\setminus \tD_r)$, it follows from \eqref{bnd:v:1} and the properties of $\tilde{u}_{\ams}$ that there exists $C>0$ such that
$$v_{t,r}\leq C \tilde{u}_{\ams}\hbox{ on }\partial\left((\tOmega\cap \tD_r)\cap B_{2r_0}(x_0)\right).$$
Since in addition $(-\Delta-V)v_{t,r}=0<(-\Delta-V)(C \tilde{u}_{\ams})$, it follows from the comparison principle that $v_{t,r}\leq C \tilde{u}_{\ams}$ in $(\tOmega\setminus \tD_r)\cap B_{2r_0}(x_0)$. Arguing similarly with $-v_{t,r}$ and using the asymptotic \eqref{asymp:tu}, we get \eqref{asymp:vtr}. 

\medskip\noindent{\bf Step 3.3: } We claim that
\begin{equation}\label{cv:vtr}
\lim_{t,r\to 0}v_{t,r}\circ\Phi_t=v_0\hbox{ in }D^{1,2}(\tOmega)_{loc, \{x_0\}^c}\cap C^1_{loc}(\overline{\tOmega}\setminus\{0,x_0\}),
\end{equation}
where $v_0$ was defined in \eqref{def:vt} and convergence in $D^{1,2}(\tOmega)_{loc, \{x_0\}^c}$ means that $\lim_{t,r\to 0}\eta v_{t,r}\circ\Phi_t=\eta v_0$ in $D^{1,2}(\tOmega)$ for all $\eta\in C^\infty(\rn)$ vanishing around $x_0$. 
Indeed, $v_{t,r}\circ\Phi_t\in D^{1,2}(\tOmega\setminus \tD_r)\subset D^{1,2}(\tOmega)$. Uniform coercivity yields weak convergence in $D^{1,2}(\tOmega)$ to $\tilde{v}\in D^{1,2}(\tOmega)$. Passing to the limit, one gets $(-\Delta-V)\tilde{v}=f_0$, so that $\tilde{v}=v_0$. Uniqueness then yields convergence in 
$C^1_{loc}(\overline{\tOmega}\setminus\{0,x_0\})$. With a change of variable, equation \eqref{eq:vtr} yields an elliptic equation for $v_{t,r}\circ\Phi_t$. Multiplying this equation by $\eta^2\cdot (v_{t,r}\circ\Phi_t-v_0)$ for $\eta\in C^\infty(\rn)$ vanishing around $x_0$, one gets convergence of $\eta v_{t,r}\circ\Phi_t$ to $\eta v_0$ in $D^{1,2}(\tOmega)$. This proves the claim.

\medskip\noindent It follows from the construction of the mass (see Theorem \ref{def:mass}) and the regularity Theorem \ref{th:hopf} that there exists $K_0\in \rr$ and for all $(t,r)$ small, there exists $K_{t,r}\in\rr$ such that
\begin{equation}\label{mass:Ktr}
v_{t,r}(x)=K_{t,r}\frac{d(x,\partial\tOmega_t)}{|x|^{\am}}+o\left(\frac{d(x,\partial\tOmega_t)}{|x|^{\am}}\right)
\end{equation}
and
\begin{equation}\label{mass:K}
v_{0}(x)=K_{0}\frac{d(x,\partial\tOmega)}{|x|^{\am}}+o\left(\frac{d(x,\partial\tOmega)}{|x|^{\am}}\right)
\end{equation}
as $x\in\tOmega$ goes to $0$. Note that around $0$, $\tOmega_{t,r}$ coincides with $\tOmega_t$. 

\medskip\noindent{\bf Step 3.4: } We claim that
\begin{equation}\label{cv:ktr}
\lim_{t,r\to 0}K_{t,r}=K_0.
\end{equation}
We only give a sketch. Noting $\tilde{v}_{t,r}:=v_{t,r}\circ\Phi_t$, the proof relies on \eqref{cv:vtr} and the fact that 
$$-\Delta_{\Phi_t^\star\eucl}\tilde{v}_{t,r}-V\circ\Phi_t\tilde{v}_{t,r}=f_t\circ\Phi_t\hbox{ in }\tOmega\cap B_\delta(0).$$
The comparison principle and the definitions \eqref{mass:Ktr} and \eqref{mass:K} then yield \eqref{cv:ktr}.

\medskip\noindent{\bf Step 4: Proof of Proposition \ref{prop:pert:mass}.} We define $\tH_0(x):=u_{\ap,0}(x)-v_0(x)$ for all $x\in\overline{\tOmega}\setminus\{0,x_0\}$, and consider its Kelvin transform
\begin{equation}\label{def:H:0:ter}
H_0(x):=|x-x_0|^{2-n}\tH_0(i(x))=|x-x_0|^{2-n}\left(u_{\ap,0}-v_0\right)(i(x))
\end{equation}
for all $x\in\Omega$. It follows from the definitions of $u_{\ap,0}$ and $v_0$ that $H_0$ satisfies the following properties:
\begin{equation}\label{eq:mathH}
\left\{\begin{array}{lll}
\hfill -\Delta H_0-\frac{\gamma}{|x|^2}H_0&=0& \hbox{ in }\Omega\\
\hfill H_0&>0& \hbox{ in }\Omega\\
\hfill H_0&=0& \hbox{ in }\partial\Omega\setminus\{0\}.
\end{array}\right.
\end{equation}
Concerning the pointwise behavior, we have that
\begin{equation}\label{exp:mathH}
H_0(x)=\frac{d(x,\partial\Omega)}{|x|^{\aps}}-K_{0}\frac{d(x,\partial\Omega)}{|x|^{\ams}}+o\left(\frac{d(x,\partial\Omega)}{|x|^{\ams}}\right)
\end{equation}
as $x\to 0$, $x\in \Omega$, and
\begin{equation}\label{bnd:H:1}
H_0(x)\leq C|x|^{1-\aps}\hbox{ for all }x\in\Omega,\, |x|>1.
\end{equation}
This proves the existence part in Proposition \ref{prop:pert:mass}. We now deal with the uniqueness. We let $H\in C^2(\overline{\Omega}\setminus\{0\})$ be as in Proposition \ref{prop:pert:mass}, and consider its Kelvin transform $\tilde{H}(x):=|x-x_0|^{2-n}H(i(x))$ for all $x\in\overline{\tOmega}\setminus \{0,x_0\}$. We then have that
\begin{equation}\label{eq:36}
\left\{\begin{array}{lll}
\hfill -\Delta \tH-V\tH&=0& \hbox{ in }\tOmega\\
\hfill \tH&>0& \hbox{ in }\tOmega\\
\hfill \tH&=0& \hbox{ in }\partial\tOmega\setminus\{0, x_0\}.
\end{array}\right.
\end{equation}
Moreover, we have that $\tH(x)\leq C|x|^{1-\ap}+C|x-x_0|^{1-\am}$ for all $x\in\tOmega$. It then follows from Theorem \ref{th:classif} that there exists $C_1,C_2>0$ such that
\begin{equation}\label{est:tH}
\tH(x)\sim_{x\to 0}C_1\frac{d(x,\partial\tOmega)}{|x|^{\alpha}}\hbox{ and }\tH(x)\sim_{x\to x_0}C_2\frac{d(x,\partial\tOmega)}{|x-x_0|^{\am}},
\end{equation}
where $\alpha\in\{\am,\ap\}$. We claim that $\alpha=\ap$. Indeed, otherwise, we would have $\tH\in D^{1,2}(\tOmega)$ (see Theorem \ref{th:classif}) and then \eqref{eq:36} and coercivity would yield $\tH\equiv 0$, which is a contradiction. Therefore $\alpha=\ap$. In particular, the estimates \eqref{est:tH} hold for $\tH_0$ (with different constants $C_1,C_2$). Arguing as in the proof of Theorem \ref{def:mass}, we get that there exists $\lambda>0$ such that $\tH=\lambda \tH_0$, and therefore $H=\lambda H_0$. This proves uniqueness and completes the proof of Proposition \ref{prop:pert:mass}.

\medskip\noindent As a consequence of \eqref{exp:mathH}, the mass $m_\gamma(\Omega)$ is well-defined and we have that
\begin{equation}\label{eq:mass:k0}
m_\gamma(\Omega)=-K_0.
\end{equation}

\medskip\noindent{\bf Step 5: convergence of the mass:} We claim that
\begin{equation}\label{cv:mass:k0}
\lim_{t\to 0, R\to\infty}m_\gamma(\Omega_{t,R})=m_\gamma(\Omega).
\end{equation}
We define $\tH_{t,r}:=u_{\aps,t}-v_{t,r}$ in such a way that
$$-\Delta \tH_{t,r}-V\tH_{t,r}=0\hbox{ in }\tOmega_{t,r}.$$
It follows from \eqref{asymp:ut} and \eqref{mass:Ktr} that $\tH_{t,r}>0$ around $0$. From the maximum principle, we deduce that  $\tH_{t,r}>0$ on $\tOmega_{t,r}$ and that it vanishes on $\partial\tOmega_{t,r}\setminus\{0,x_0\}$. 

\medskip\noindent It follows from \eqref{asymp:ut} and \eqref{mass:Ktr} that
\begin{equation}\nonumber
\tH_{t,r}(x)=\frac{d(x,\partial\tOmega_{t,r})}{|x|^{\aps}}-K_{t,r}\frac{d(x,\partial\tOmega_{t,r})}{|x|^{\ams}}+o\left(\frac{d(x,\partial\tOmega_{t,r})}{|x|^{\ams}}\right)
\end{equation}
as $x\to 0$, $x\in \tOmega_{t,r}$. Coming back to $\Omega_{t,R}$ with $R=r^{-1}$ via the inversion $i$ with $H_{t,R}(x):=|x-x_0|^{2-n}\tH_{t,r}(i(x))$ for all $x\in \Omega_{t,R}$, we get that
\begin{equation}\nonumber
\left\{\begin{array}{lll}
\hfill -\Delta H_{t,R}-\frac{\gamma}{|x|^2}H_{t,R}&=0& \hbox{ in }\Omega_{t,R}\\
\hfill  H_{t,R}&>0& \hbox{ in }\Omega_{t,R}\\
\hfill H_{t,R}&=0& \hbox{ in }\partial\Omega_{t,R}\setminus\{0\}
\end{array}\right.
\end{equation}
and
\begin{equation}\nonumber
H_{t,R}(x)=\frac{d(x,\partial\Omega_{t,R})}{|x|^{\aps}}-K_{t,r}\frac{d(x,\partial\Omega_{t,R})}{|x|^{\ams}}+o\left(\frac{d(x,\partial\Omega_{t,R})}{|x|^{\ams}}\right)
\end{equation}
as $x\to 0$, $x\in \Omega_{t,R}$. Therefore, it follows from the definition of the mass (see Theorem \ref{def:mass}) that $m_\gamma(\Omega_{t,R})=-K_{t,r}$ for all $t,r$, $R=r^{-1}$. Claim \eqref{cv:mass:k0} then follows from \eqref{cv:ktr} and \eqref{eq:mass:k0}.

\medskip\noindent This ends the proofs of Propositions \ref{prop:pert:mass} and \ref{prop:mass:positiv}. 

\medskip\noindent{\bf Step 6:} In order to prove Proposition \ref{prop:ex:any:behave}, we need to exhibit prototypes of  unbounded domains with either positive or negative mass.

\begin{proposition}\label{prop:ex:pos:unbounded} Let $\Omega$ be a domain such that $0\in\partial\Omega$ and $\Omega$ is smooth at infinity. Assume that $\gamma_H(\Omega)>\frac{n^2-1}{4}$ and fix $\gamma\in \left(\frac{n^2-1}{4},\gamma_H(\Omega)\right)$. Then $m_\gamma(\Omega)>0$ if $\rnp\subsetneq \Omega$, and $m_\gamma(\Omega)<0$ if $\Omega\subsetneq \rnp$.
\end{proposition}

\medskip\noindent{\it Proof of Proposition \ref{prop:ex:pos:unbounded} :} With $H_0$ defined as in \eqref{def:H:0:ter}, we set
$$\mathcal{U}(x):=H_0(x)-x_1|x|^{-\aps}\hbox{ for all }x\in\Omega.$$
We first assume that $\rnp\subset \Omega$ while being distinct. We then have that
\begin{equation}\label{eq:mathU}
\left\{\begin{array}{ll}
-\Delta \mathcal{U}-\frac{\gamma}{|x|^2}\mathcal{U}=0& \hbox{ in }\rnp\\
\mathcal{U}\gneqq 0& \hbox{ in }\partial\rnp\setminus\{0\}.
\end{array}\right.
\end{equation}
\noindent We claim that
\begin{equation}\label{bnd:U:nrj}
\int_{\rnp}|\nabla \mathcal{U}|^2\, dx<+\infty.
\end{equation}
Indeed, at infinity, this is the consequence of the fact that $|\nabla \mathcal{U}|(x)\leq C|x|^{-\aps}$ for all $x\in\rnp$ large, this latest bound being a consequence of \eqref{bnd:H:1} combined with elliptic regularity theory. At zero, the argument is different. Indeed, one first notes that $d(x,\partial\Omega')=x_1+O(|x|^2)$ for $x\in\rnp$ close to $0$, and therefore, $\mathcal{U}(x)=O(|x|^{1-\ams})$ for $x\to 0$. The control on the gradient $|\nabla \mathcal{U}|(x)\leq C|x|^{-\ams}$ at $0$ follows from the construction of $\tH_0$. This yields integrability at $0$ and proves \eqref{bnd:U:nrj}.

\medskip\noindent We claim that $\mathcal{U}>0$ in $\rnp$. Indeed, 
it follows from \eqref{eq:mathU} and \eqref{bnd:U:nrj} that $\mathcal{U}_-\in D^{1,2}(\rnp)$. Multiplying equation \eqref{eq:mathH} by $\mathcal{U}_-$, integrating by parts on $(B_{R}(0)\setminus B_{\epsilon}(0))\cap\rnp$, and letting $\epsilon\to 0$ and $R\to +\infty$ by using \eqref{bnd:U:nrj}, one gets $\mathcal{U}_-\equiv 0$, and then $\mathcal{U}\geq 0$. The result follows from Hopf's maximum principle.

\medskip\noindent We now claim that
\begin{equation}\label{sign:K}
m_\gamma(\Omega)>0.
\end{equation}
Indeed, since $\mathcal{U}>0$ in $\rnp$, there exists $c_0>0$ such that $\mathcal U(x)\geq c_0 x_1|x|^{-\ams}$ for all $x\in \partial(B_1(0)_+)$. It then follows from \eqref{bnd:U:nrj}, \eqref{eq:mathU} and the comparison principle that $\mathcal U(x)\geq c_0 x_1|x|^{-\ams}$ for all $x\in B_1(0)_+$. The expansion \eqref{exp:mathH} then yields $-K_0\geq c_0>0$. This combined with \eqref{eq:mass:k0}  proves the claim.

\medskip\noindent When $\Omega\subset\rnp$, the argument is similar except that one works on $\Omega$ (and not $\rnp$) and that $\mathcal{U}\lneqq 0$ in $\partial\Omega\setminus\{0\}$. This ends the proof of Proposition \ref{prop:ex:pos:unbounded}.

\medskip\noindent{\bf Step 7: Proof of Proposition \ref{prop:ex:any:behave}:} Let $\omega$ be a smooth domain of $\rn$ such that $0\in\partial\Omega$. Up to a rotation, there exists $\varphi\in C^\infty(\rr^{n-1})$ such that $\varphi(0)=0$, $\nabla\varphi(0)=0$ and there exists $\delta_0>0$ such that
$$\omega\cap B_{\delta_0}(0)=\{x_1>\varphi(x')/\, (x_1,x')\in B_{\delta_0}(0)\}.$$
Let $\eta\in C^\infty_c(B_{\delta_0}(0))$ be such that $\eta(x)=1$ for all $x\in B_{\delta_0/2}(0)$, and define
$$\Phi_t(x):=\left(x_1+\eta(x)\frac{\varphi(tx')}{t}, x'\right)\hbox{ for all }t> 0\hbox{ and }x\in\rn,$$ 
and $\Phi_0:=Id_{\rn}$. It is easy to see that $\Phi_t$ satisfies the hypothesis of Proposition \ref{prop:mass:positiv}. Moreover, for $0<t<1$, we have that 
$$\frac{\omega}{t}\cap \Phi_t(B_{\delta_0/2}(0))=\Phi_t(\rnp\cap B_{\delta_0/2}(0)).$$
We let $\Omega$ be a smooth domain at infinity such that
\begin{equation}\label{hyp:omega:ex}
\Omega\cap B_1(0)=\rnp\cap B_1(0)\hbox{ and }\gamma_H(\Omega)>\frac{n^2-1}{4}.
\end{equation}
(for example, $\rnp$), and let $\Omega_{t,R}$ be as in Proposition \ref{prop:mass:positiv}. It is easy to see  that
$$\omega\cap t\Phi_t(B_{\delta_0/2}(0))=t\Omega_{t,R}\cap t\Phi_t(B_{\delta_0/2}(0)).$$
Therefore, for $t>0$ small enough, we have that
$$\omega\cap B_{t\delta_0/3}(0) =t\Omega_{t,R}\cap  B_{t\delta_0/3}(0).$$
Moreover, $\gamma_H(t\Omega_{t,R})=\gamma_H(\Omega_{t,R})>(n^2-1)/4$ as $t\to 0$ and $R\to +\infty$ (see \eqref{lim:gamma:t}). Concerning the mass, we have that
$$t^{\ap-\am}m_\gamma(t\Omega_{t,R})=m_\gamma(\Omega_{t,R})\to m_\gamma(\Omega)\hbox{ as }t\to 0,R\to +\infty.$$
We now choose $\Omega$ appropriately. 

\noindent To get a negative mass, we choose $\Omega$ smooth at infinity such that $\Omega\cap B_1(0)=\rnp\cap B_1(0)$ and $\Omega\subsetneq \rnp$. Then $\gamma_H(\Omega)=n^2/4$, \eqref{hyp:omega:ex} holds and Proposition \ref{prop:ex:pos:unbounded} yields $m_\gamma(\Omega)<0$. 

\noindent To get a positive mass, we choose $\rnp\subsetneq\Omega$ such that  \eqref{hyp:omega:ex} holds (this is possible for any value of $\gamma_H(\Omega)$ arbitrarily close to $\frac{n^2}{4}$, see point (5) of Proposition \ref{prop:gamma}). Then Proposition \ref{prop:ex:pos:unbounded} yields $m_\gamma(\Omega)>0$. This proves Proposition \ref{prop:ex:any:behave}.

\section{The remaining cases corresponding to  $s=0$ and $n=3$\label{sec:n:3}}
\noindent According to Proposition \ref{prop:s0:gamma:neg} and Theorem \ref{main}, there remains the situation, when $s=0$, $n=3$ and $\gamma \in (0, \frac{n^2}{4})$. Note first, that if $\gamma \geq \gamma_H(\Omega)$, we have from Proposition \ref{prop:pptes:inf} and Theorem \ref{tool} that $\mu_{\gamma, 0}(\Omega)\leq 0<\mu_{\gamma, 0}(\rnp)$ and the existence of extremals is guaranteed. Another situation is when $\mu_{\gamma, 0}(\rnp)$ does have an extremal $U$. In this case, Proposition \ref{prop:test:fct} provides sufficient conditions for $\mu_{\gamma, 0}(\Omega)<\mu_{\gamma, 0}(\rnp)$, and hence there are extremals by again using Theorem \ref{tool}. 
The rest of this section addresses the remaining case, that is when $\gamma\in (0,\gamma_H(\Omega))$ and when $\mu_{\gamma, 0}(\rnp)$ has no extremal, and therefore  $\mu_{\gamma,0}(\rr^3_+)=K(3,2)^{-2}$ according to Proposition \ref{prop:ext:rnp}.

\medskip\noindent We first define the ``interior" mass in the spirit of Schoen-Yau \cite{SY}.
\begin{proposition}\label{prop:def:mass:int} Let $\Omega\subset \rr^3$ be an open smooth bounded domain such that $0\in\partial\Omega$. 
If $\gamma\in (0,\gamma_H(\Omega))$, then the equation 
$$\left\{\begin{array}{ll}
-\Delta G-\frac{\gamma}{|x|^2} G=0 &\hbox{\rm in }\Omega\setminus\{x_0\}\\
\hfill G>0 &\hbox{\rm in }\Omega\setminus\{x_0\}\\
\hfill G=0 &\hbox{\rm on }\partial\Omega\setminus\{0\}
\end{array}\right.$$
has a solution $G$, that is unique up to multiplication by a constant.  
Moreover,  for any $x_0\in \Omega$, there exists a unique $R_\gamma(x_0)\in\rr$ independent of the choixe of $G$ and $c_G>0$ such that
$$G(x)=c_G\left(\frac{1}{|x-x_0|^{n-2}}+R_\gamma(x_0)\right)+o(1)\hbox{ as }x\to x_0.$$

\end{proposition}

\smallskip\noindent{\it Proof of Proposition \ref{prop:def:mass:int}.} Since $\gamma<\gamma_H(\Omega)$, the operator $-\Delta-\gamma|x|^{-2}$ is coercive and we can consider $G$ to be its Green's function  on $\Omega$ with Dirichlet boundary condition.  In particular, for any $\varphi\in C^\infty_c(\Omega)$, we have that
$$\varphi(x)=\int_\Omega G_x(y)\left(-\Delta \varphi(y)-\gamma\frac{\varphi(y)}{|y|^2}\right)\, dy \quad\hbox{for $x\in\Omega$,}$$
where $G_x:=G(x,\cdot)$. Fix $x_0\in\Omega$ and let $\eta\in C^\infty_c(\Omega)$ be such that $\eta(x)=1$ around $x_0$. Define the distribution $\beta_{x_0}:\Omega\to\rr$ as 
$$\hbox{$G_{x_0}(x)=\frac{1}{\omega_2}\left(\frac{\eta(x)}{|x-x_0|}+\beta_{x_0}(x)\right)$ \quad for all $x\in\Omega$. }$$
Set 
$$\hbox{$f(x):=-\left(-\Delta-\frac{\gamma}{|x|^2}\right)\left(\frac{\eta(x)}{|x-x_0|}\right)$ \quad for all $x\neq x_0$. }$$
In particular,
$$\left(-\Delta-\frac{\gamma}{|x|^2}\right)\beta_{x_0}=f\hbox{ in the distribution sense.}$$
On can easily see that there exists $C>0$ such that
$$\hbox{$|f(x)|\leq C|x-x_0|^{-1}$\quad for all $x\in\Omega$.}$$
 Therefore $f\in L^2(\Omega)$ and, by uniqueness of the Green's function (since the operator is coercive), we have that $\beta_{x_0}\in \dundeux$. It follows from standard elliptic theory that $\beta_{x_0}\in C^\infty(\overline{\Omega}\setminus\{0,x_0\})\cap C^{0,\theta}(\overline{\Omega}\setminus B_\delta(0))$ for all $\theta\in (0,1)$ and $\delta>0$. In addition, for any $\theta\in (0,1)$ and $\delta>0$, there exists $C_\theta>0$ such that
\begin{equation}\label{bnd:der:beta:1}
\hbox{$|\nabla\beta_{x_0}(x)|\leq C_\theta|x|^{\theta-1}$\quad for all $x\in \Omega\setminus B_\delta(0)$.}
\end{equation}
Since $f$ vanishes around $0$, it follows from Theorem \ref{th:hopf} and Lemma \ref{lem:deriv:1} that
\begin{equation}\label{bnd:der:beta:2}
\hbox{$\beta_{x_0}(x)=O(|x|^{1-\am})$ \quad and \quad $|\nabla \beta_{x_0}(x)|=O(|x|^{-\am})$\quad when $x\to 0$. }
\end{equation}
We can therefore define the {\it mass of $\Omega$} at $x_0$ associated to the operator $L_\gamma$ by 
$$R_\gamma(\Omega, x_0):=\beta_{x_0}(x_0).$$
One can easily check that this quantity is independent of the choice of $\eta$, which finishes the proof of Proposition \ref{prop:def:mass:int}.

\begin{lem}\label{lem:dim3} Let $\Omega\subset \rr^3$ be an open smooth bounded domain such that $0\in\partial\Omega$ and $x_0\in\Omega$. Assume that $\gamma\in (0,\gamma_H(\Omega))$ and that $\mu_{\gamma,0}(\rr^3_+)=K(3,2)^{-2}$. Then, there exists a family $(u_\epsilon)_\epsilon$ in $D^{1,2}(\Omega)$ such that
\begin{equation}\label{asymp:lim:app}
\hbox{$J^\Omega_{\gamma,0}(\ue)=\frac{1}{K(n,2)^2}\left(1-\frac{\omega_2 R_\gamma(x_0)}{3\int_{\rr^3}U^{\crit}\, dx}\eps+o(\eps)\right)$ as $\eps\to 0$,}
\end{equation}
where $U(x):=(1+|x|^2)^{-1/2}$ for all $x\in\rr^3$ and $\crit=\crit(0)=\frac{2n}{n-2}$.

\end{lem} 
\smallskip\noindent{\it Proof of Lemma \ref{lem:dim3}:} We proceed as in Schoen \cite{schoen1} (see Druet \cites{d2,d4} and Jaber \cite{jaber}). The computations are similar to the case $\gamma>\frac{n^2-1}{4}$ performed in Section \ref{sec:best:cst}. For $\eps>0$, define the functions
$$
\hbox{$\ue(x):=\eta(x)\left(\frac{\eps}{\eps^2+|x-x_0|^2}\right)^{\frac{1}{2}}+\eps^{\frac{1}{2}}\beta_{x_0}(x)$\quad for all $x\in\Omega$.}$$
 One can easily check that $\ue\in \dundeux$. We now estimate $J^\Omega_{\gamma,0}(\ue)$.

\smallskip\noindent In the sequel, $\Theta_c(\eps)$ denote any quantity such that
$$\lim_{c\to 0}\lim_{\eps\to 0}\frac{\Theta_c(\eps)}{\eps}=0.$$
We first claim that
\begin{equation}\label{lim:step:1}
\int_{\Omega\setminus B_c(x_0)}\left(|\nabla \ue|^2-\gamma\frac{\ue^2}{|x|^2}\right)\, dx=\omega_2c^{-1}\eps+\omega_2 R_\gamma(x_0)\eps+\Theta_c(\eps).
\end{equation}
Indeed, it is clear that $\eps^{-\frac{1}{2}}\ue\to G'_{x_0}:=\omega_2 G_{x_0}$ in $C^2_{loc}(\overline{\Omega}\setminus\{0,x_0\})$. Therefore, Lebesgue's dominated convergence theorem yields
\begin{equation}\label{lim:step:1:2}
\lim_{\eps\to 0}\frac{\int_{\Omega\setminus B_c(x_0)}\left(|\nabla \ue|^2-\gamma\frac{\ue^2}{|x|^2}\right)\, dx}{\eps}=\int_{\Omega\setminus B_c(x_0)}\left(|\nabla G'_{x_0}|^2-\gamma\frac{(G'_{x_0})^2}{|x|^2}\right)\, dx.
\end{equation}
Integrating by parts and using \eqref{bnd:der:beta:2}, as $\delta>0$ goes to $0$, we have that
\begin{eqnarray*}
&&\int_{\Omega\setminus B_c(x_0)}\left(|\nabla G'_{x_0}|^2-\gamma\frac{(G'_{x_0})^2}{|x|^2}\right)\, dx\\
&&= \int_{\Omega\setminus (B_c(x_0)\cup B_\delta(0)}\left(|\nabla G'_{x_0}|^2-\gamma\frac{(G'_{x_0})^2}{|x|^2}\right)\, dx+o(1)\\
&&= \int_{\Omega\setminus (B_c(x_0)\cup B_\delta(0)}\left(-\Delta G'_{x_0}-\gamma\frac{G'_{x_0}}{|x|^2}\right)G'_{x_0}\, dx-\int_{\partial B_c(x_0)}G'_{x_0}\partial_\nu G'_{x_0}\, d\sigma\\
&&-\int_{\partial B_\delta(0)}G'_{x_0}\partial_\nu G'_{x_0}\, d\sigma +o(1)\\
&&= -\int_{\partial B_c(x_0)}G'_{x_0}\partial_\nu G'_{x_0}\, d\sigma+O(\delta^{n-1}\delta^{1-\am}\delta^{-\am})+o(1) \quad\hbox{as $\delta\to 0$.}
\end{eqnarray*}
Since $\am<n/2$, we then have that
$$\int_{\Omega\setminus B_c(x_0)}\left(|\nabla G'_{x_0}|^2-\gamma\frac{(G'_{x_0})^2}{|x|^2}\right)\, dx =-\int_{\partial B_c(x_0)}G'_{x_0}\partial_\nu G'_{x_0}\, d\sigma.$$
With the definition $R_\gamma(x_0)$ and \eqref{bnd:der:beta:1}, we have that $G'_{x_0}=c^{-1}+R_\gamma(x_0)+O(c^\theta)$ and  $\partial_\nu G'_{x_0}(x)=-c^{-2}+O(c^{\theta-1})$ on $\partial B_c(x_0)$ as $c\to 0$. Therefore
$$\hbox{$-\int_{\partial B_c(x_0)}G'_{x_0}\partial_\nu G'_{x_0}\, d\sigma=\omega_2c^{-1}+\omega_2R_\gamma(x_0)+O(c^\theta)$\quad as $c\to 0$. }$$
Combined with \eqref{lim:step:1:2}, this proves \eqref{lim:step:1}.

\noindent Now define for each $\eps>0$, the function 
$$\hbox{$U_\eps(x):=\left(\frac{\eps}{\eps^2+|x-x_0|^2}\right)^{\frac{1}{2}}$ \quad for all $x\in\rr^3$,}$$
 and set $U(x):=(1+|x|^2)^{-1/2}$ for all $x\in\rr^3$. It is clear that $\Delta U=3 U^{\crit-1}$. 
\noindent
We claim that
\begin{equation}\label{lim:step:2}
\int_{B_c(x_0)}\left(|\nabla \ue|^2-\gamma\frac{\ue^2}{|x|^2}\right)\, dx=3\int_{\rr^3}U^{\crit}\, dx-\omega_2c^{-1}\eps+\Theta_c(\eps).
\end{equation}
Indeed, note first $|\ue(x)|\leq C\sqrt{\eps}|x-x_0|^{-1}$ for all $\eps>0$ and all $x\in\Omega$ close to $x_0$. Therefore, for $c>0$ small enough, we have that
\begin{equation}
\int_{B_c(x_0)}\frac{\ue^2}{|x|^2}\, dx=\Theta_c(\eps).
\end{equation}
Using that $\beta_{x_0}\in \dundeux$ and integrating by parts, we get that
\begin{eqnarray}
\int_{B_c(x_0)}|\nabla \ue|^2\, dx&=&\int_{B_c(x_0)}|\nabla (U_\eps+\sqrt{\eps}\beta_{x_0})|^2\, dx\nonumber\\
&=& \int_{B_c(x_0)}|\nabla U_\eps|^2\, dx+2\sqrt{\eps}\int_{B_c(x_0)}\nabla U_\eps\nabla\beta_{x_0}\, dx\nonumber\\
&&+\eps  \int_{B_c(x_0)}|\nabla \beta_{x_0}|^2\, dx\nonumber\\
&=& \int_{B_c(x_0)}U_\eps(-\Delta U_\eps)\, dx+\int_{\partial B_c(x_0)}U_\eps\partial U_\eps\, d\sigma\nonumber\\
&&+2\sqrt{\eps}\int_{B_c(x_0)}\beta_{x_0}(-\Delta U_\eps)\, dx\nonumber\\
&&+2\sqrt{\eps}\int_{\partial B_c(x_0)}\beta_{x_0}\partial_\nu U_\eps\, d\sigma+\Theta_c(\eps)\label{eq:app:1}
\end{eqnarray}
Since $\eps^{-1/2}U_\eps\to |\cdot-x_0|^{-1}$ in $C^1_{loc}(\rr^3\setminus\{0\})$, we get that
\begin{equation}\label{eq:app:2}
\int_{\partial B_c(x_0)}U_\eps\partial_\nu U_\eps\, d\sigma=-\omega_2 c^{-1}\eps+o(\eps) \hbox{\quad as $\eps\to 0$. }
\end{equation}
Using in addition that $\beta_{x_0}\in C^{0,\theta}$ around $x_0$, we get as $\eps\to 0$ and for $c>0$ small, that
\begin{equation}\label{eq:app:3}
\int_{\partial B_c(x_0)}\beta_{x_0}\partial_\nu U_\eps\, d\sigma=-\sqrt{\eps}\omega_2 R_\gamma(x_0)+O(c^\theta\sqrt{\eps})
\end{equation}
Plugging \eqref{eq:app:2} and \eqref{eq:app:3} into \eqref{eq:app:1} yields
\begin{eqnarray}\label{eq:app:4}
\qquad \int_{B_c(x_0)}|\nabla \ue|^2\, dx&=& 3 \int_{B_{c/\eps}(0)}U^{\crit}\, dx-\omega_2 c^{-1}\eps+o(\eps)\\
&&+2\sqrt{\eps}\int_{B_c(x_0)}\beta_{x_0} (-\Delta U_\eps)\, dx-2\eps\omega_2 R_\gamma(x_0)+\Theta_c(\eps).\nonumber 
\end{eqnarray}
It is easy to check that $\int_{B_{c/\eps}(0)}U^{\crit}\, dx=\int_{\rr^3}U^{\crit}\, dx+o(\eps)$ as $\eps\to 0$. For $\theta\in (1/2,1)$ we have that 
\begin{equation}\label{eq:app:5}
\hbox{$\int_{B_c(x_0)}|\Delta U_\eps|\cdot |x-x_0|^\theta\, dx=o(\eps)$\quad as $\eps\to 0$.}
\end{equation}
Integrating by parts and using that $\eps^{-1/2}U_\eps(x)\to |x-x_0|^{-1}$ in $C^1_{loc}(\rr^3\setminus\{0\})$, we get that as $\eps\to 0$, 
\begin{eqnarray}
\int_{B_c(x_0)}-\Delta U_\eps\, dx&=& -\int_{\partial B_c(x_0)}\partial_\nu U_\eps\, d\sigma\nonumber\\
&=& -\sqrt{\eps}\int_{\partial B_c(x_0)}\partial_\nu |x-x_0|^{-1}\, d\sigma+o(\eps)\nonumber\\
&=& \omega_2\sqrt{\eps}+o(\eps)\label{eq:app:6}
\end{eqnarray}
Plugging \eqref{eq:app:5} and \eqref{eq:app:6} into \eqref{eq:app:4} and using that $\beta_{x_0}(x)=R_\gamma(x_0)+O(|x-x_0|^\theta)$, we get \eqref{lim:step:2}.

\medskip\noindent Putting together \eqref{lim:step:1} and \eqref{lim:step:1:2} yields
\begin{equation}\label{lim:nabla}
\hbox{$\int_{\Omega}\left(|\nabla \ue|^2-\gamma\frac{\ue^2}{|x|^2}\right)\, dx=3\int_{\rr^3}U^{\crit}\, dx+\omega_2 R_\gamma(x_0)\eps+o(\eps)$\quad as $\eps\to 0$.}
\end{equation}
\medskip\noindent We now claim that
\begin{equation}\label{lim:denom}
\hbox{$\int_\Omega\ue^{\crit}\, dx =\int_{\rr^3}U^{\crit}\, dx+\frac{\crit}{3}\omega_2 R_\gamma(x_0)\eps+o(\eps)$\quad as $\eps\to 0$. }
\end{equation}
Using \eqref{ineq:not:numb}, the boundedness of $\beta_{x_0}$ around $x_0$, and the above computations, we get that
\begin{eqnarray*}
\int_\Omega\ue^{\crit}\, dx &=& \int_{B_c(x_0)}\ue^{\crit}\, dx +o(\eps)\\
&=& \int_{B_c(x_0)}|U_\eps+\sqrt{\eps}\beta_{x_0}|^{\crit}\, dx +o(\eps)\\
&=& \int_{B_c(x_0)}U_\eps^{\crit}\, dx +\crit\sqrt{\eps}\int_{B_c(x_0)}\beta_{x_0}U_\eps^{\crit-1}\, dx\\&&+O\left(\int_\Omega \left(\eps U_\eps^{\crit-2}\beta_{x_0}^2+|\sqrt{\eps}\beta_{x_0}|^{\crit}\right)\, dx\right)+o(\eps)\\
&=& \int_{B_{c/\eps}(0)}U^{\crit}\, dx +\frac{\crit}{3}\omega_2 R_\gamma(x_0)\eps+o(\eps)\\
&=& \int_{\rr^3}U^{\crit}\, dx +\frac{\crit}{3}\omega_2 R_\gamma(x_0)\eps+o(\eps) \hbox{\quad as $\eps\to 0$, }
\end{eqnarray*}
which proves \eqref{lim:denom}.

\medskip\noindent Putting together \eqref{lim:nabla} and \eqref{lim:denom} yields
\begin{equation}
\frac{\int_{\Omega}\left(|\nabla \ue|^2-\gamma\frac{\ue^2}{|x|^2}\right)\, dx}{\left(\int_\Omega\ue^{\crit}\, dx \right)^{\frac{2}{\crit}}}=\frac{3\int_{\rr^3}U^{\crit}\, dx }{\left(\int_{\rr^3}U^{\crit}\, dx \right)^{\frac{2}{\crit}}}\left(1-\frac{\omega_2 R_\gamma(x_0)}{3\int_{\rr^3}U^{\crit}\, dx}\eps+o(\eps)\right)
\end{equation}
as $\eps\to 0$. Since $\Delta U=3 U^{\crit-1}$ and $U$ is an extremal for the Sobolev inequality $\mu_{0,0}(\rr^3)$, we have that
\begin{equation*}
\hbox{$J^\Omega_{\gamma,0}(\ue)=\frac{1}{K(n,2)^2}\left(1-\frac{\omega_2 R_\gamma(x_0)}{3\int_{\rr^3}U^{\crit}\, dx}\eps+o(\eps)\right)$\quad as $\eps\to 0$.}
\end{equation*}
This proves Lemma \ref{lem:dim3}.\hfill$\Box$

\medskip\noindent We finally get the following.

\begin{theorem}\label{th:3}  Let $\Omega$ be a bounded smooth domain of $\rr^3$ such that $0\in \partial \Omega$. 
\begin{enumerate}
\item If $\gamma \geq \gamma_H(\Omega)$, then there are extremals for $\mu_{\gamma,0}(\Omega)$.
\item If $\gamma\leq 0$, then there are no extremals for $\mu_{\gamma,0}(\Omega)$.
\item If $0<\gamma <\gamma_H(\Omega)$ and there are extremals for $\mu_{\gamma, 0}(\rnp)$, then there are extremals for $\mu_{\gamma,0}(\Omega)$ under either one of the following conditions:
\begin{itemize}
\item  $\gamma\leq\frac{n^2-1}{4}$ and the mean curvature of $\partial \Omega$ at $0$ is negative.
\item  $\gamma>\frac{n^2-1}{4}$ and the mass $m_\gamma(\Omega)$ is positive.
\end{itemize}  
\item If $0<\gamma <\gamma_H(\Omega)$ and there are no extremals for $\mu_{\gamma, 0}(\rnp)$, then there are extremals for $\mu_{\gamma,0}(\Omega)$ if there exists $x_0\in\Omega$ such that $R_\gamma(\Omega, x_0)>0$.
\end{enumerate}
\end{theorem}

\smallskip\noindent{\it Proof of Theorem \ref{th:3}:} The two first points
of the theorem follow from Proposition \ref{prop:s0:gamma:neg} and Theorem \ref{tool}. The third point follows from Proposition \ref{prop:test:fct}. For the fourth point, in this situation, it follows from Proposition \ref{prop:ext:rnp} that $\mu_{\gamma, 0}(\rnp)=\frac{1}{K(3,2)^2}$, and then Lemma \ref{lem:dim3} yields $\mu_{\gamma,0}(\Omega)<\mu_{\gamma,0}(\rnp)$, which yields the existence of extremals by Theorem \ref{tool}. This proves Theorem \ref{th:3}.\hfill$\Box$




\section{Appendix 1: Existence of extremals for $\mu_{\gamma, s}(\rnp)$}\label{sec:app:1}
The following result is used frequently throughout this paper. It has been noted in other contexts by Bartsche-Peng-Zhang \cite{BPZ} and Lin-Wang \cite{CL5}. We sketch an independent proof for the convenience of the reader.

\begin{proposition}\label{prop:ext:rnp} Fix $\gamma<\frac{n^2}{4}$ and $s\in [0,2)$ with $n\geq 3$. Then, 
\begin{enumerate}
\item If $\{s>0\}$ or $\{s=0$, $\gamma>0$ and $n\geq 4\}$, then there are extremals for $\mu_{\gamma, s}(\rnp)$.
\item If $\{s=0$ and $\gamma\leq 0\}$, there are no extremals for $\mu_{\gamma,0}(\rnp)$.
\item If there are no extremals for $\mu_{\gamma, 0}(\rnp)$, then 
\begin{equation}
\mu_{\gamma, 0}(\rnp)=\frac{1}{K(n,2)^2}:=\inf_{u\in D^{1,2}(\rn)\setminus\{0\}}\frac{\int_{\rn}|\nabla u|^2\, dx}{\left(\int_{\rn}|u|^{\crit}\, dx\right)^{\frac{2}{\crit}}}.
\end{equation}
\end{enumerate}

\end{proposition}
As a consequence, the only unknown situation is when $\{s=0$, $n=3$ and $\gamma>0\}$.

\smallskip\noindent{\it Proof of Proposition \ref{prop:ext:rnp}:} This goes as the proof of the existence of extremals for the classical Sobolev inequalities using Lions's concentration-compactness Lemmae (\cite{lions1,lions2}, see also Struwe \cite{st} for an exposition in book form). Here, we follow the notations of the proof performed in Filippucci-Pucci-Robert (\cite{FiPuRo}, Theorem 4 of Appendix A).

\medskip\noindent Point (2) is a direct consequence of Proposition \ref{prop:s0:gamma:neg}. In particular, it follows from the computations of Proposition \ref{prop:s0:gamma:neg} that
\begin{equation}\label{bnd:mu}
\mu_{\gamma, 0}(\rnp)\leq K(n,2)^{-1}\hbox{ for all }\gamma<\frac{n^2}{4}.
\end{equation}


\medskip\noindent For points (1) and (3), we consider a minimizing sequence $(\tilde{u}_k)_k\in \dundeuxrnp$  for $\mu_{\gamma, s}(\rnp)$ such that
$$\int_{\rnp}\frac{|\tilde{u}_k|^{\crits}}{|x|^s}\, dx=1\hbox{ and }\lim_{k\to +\infty}\int_{\rnp}\left(|\nabla \tilde{u}_k|^2-\frac{\gamma}{|x|^2}\tilde{u}_k^2\right)\, dx=\mu_{\gamma, s}(\rnp).$$
We use a concentration compactness argument  in the spirit of Lions \cite{lions1,lions2}. For any $k$, there exists $r_k>0$ such that $\int_{B_{r_k}(0)_+}\frac{|\tilde{u}_k|^{\crits}}{|x|^s}\, dx=1/2$. We define $u_k(x):=r_k^{\frac{n-2}{2}}u_k(r_k x)$ for all $x\in\rnp$. We then have that
\begin{equation}\label{mes:l}
\lim_{k\to +\infty}\int_{\rnp}\left(|\nabla u_k|^2-\frac{\gamma}{|x|^2}u_k^2\right)\, dx=\mu_{\gamma, s}(\rnp),
\end{equation}
and
\begin{equation}\label{mes:n}
\int_{\rnp}\frac{|u_k|^{\crits}}{|x|^s}\, dx=1 \; , \; \int_{B_1(0)_+}\frac{|u_k|^{\crits}}{|x|^s}\, dx=\frac{1}{2}.
\end{equation}
Up to a subsequence, there exists $u\in \dundeuxrnp$ such that $u_k\rightharpoonup u$ weakly in $\dundeuxrnp$ as $k\to +\infty$. In particular for any $1\leq q<\frac{2n}{n-2}$, the sequence $(u_k)_k$ goes to $u$ in $L^q_{loc}(\rn)$ when $k\to +\infty$. We define the measure $\lambda_k$ on $\rn$ as
$$\lambda_k:=\left(|\nabla u_k|^2-\frac{\gamma}{|x|^2}u_k^2\right){\bf 1}_{\rnp}\, dx.$$
We define the measure $\lambda$ and $\nu$ on $\rn$ as 
\begin{equation*}
\lambda:=\lim_{k\to +\infty}\lambda_k\hbox{ and }\nu:=\lim_{k\to +\infty}\frac{|u_k|^{\crits}}{|x|^s}{\bf 1}_{\rnp}\, dx,
\end{equation*}
both limit being taken in the sense of measures. Arguing as in Lions \cite{lions1,lions2} (see also Struwe \cite{st} or Filippucci-Pucci-Robert \cite{FiPuRo} for the present notations), we get that one and only one of the two following situations occur:
\begin{enumerate}
\item Either $\lambda=\left(|\nabla u|^2-\frac{\gamma}{|x|^2}u^2\right)\, dx$, $\nu=\frac{|u|^{\crits}}{|x|^s}\, dx$, and $u$ is an extremal for $\mu_{\gamma, s}(\rnp)$; 
 \item or there exists $x_0\neq 0$  such that $\lambda=\mu_{\gamma, s}(\rnp)\delta_{x_0}$ and $\nu= \delta_{x_0}.$
\end{enumerate}
In the sequel, we assume that there are no extremals for $\mu_{\gamma, s}(\rnp)$. We therefore have that (2) holds. This yields $u\equiv 0$. 

\medskip\noindent We shall now show that $s=0$. \\
Indeed, if not, then  $\crits<\frac{2n}{n-2}$, and since $x_0\neq 0$, we then get that for $\delta<|x_0|/2$, we have $\lim_{k\to +\infty}\int_{B_\delta(x_0)_+}\frac{|u_k|^{\crits}}{|x|^s}=0$, but this contradicts that $\nu=\delta_{x_0}$, which proves our claim that $s=0$.

\medskip\noindent We now prove  that for all $\delta>0$, we have that
\begin{equation}
\lim_{k\to +\infty}\int_{\rnp\setminus B_\delta(x_0)}|\nabla u_k|^2\, dx=0.
\end{equation}
Indeed,  let $\varphi\in C^\infty(\rn)$ be such that $0\leq\varphi\leq 1$, $\varphi(x)=1$ for all $x\in B_\delta(x_0)$ and $\varphi(x)=0$ for all $x\in\rn\setminus B_{2\delta}(x_0)$. We have that $(1-\varphi)u_k\in \dundeuxrnp$ for all $k$, and therefore, integrating by parts and using that $u_k\to 0$ in $L^2_{loc}(\rnp)$ as $k\to +\infty$, we get that
\begin{eqnarray*}
\int_{\rnp}\left(|\nabla(1-\varphi)u_k|^2-\frac{\gamma}{|x|^2}((1-\varphi)u_k)^2\right)\, dx&=& \int_{\rnp}(1-\varphi)^2\left(|\nabla u_k|^2-\frac{\gamma}{|x|^2}u_k^2\right)\, dx\\
&&+o(1)\\
&=& \int_{\rnp}d\lambda_k+\langle\lambda_k,\varphi(\varphi-2)\rangle+o(1)\\
&= &\mu_{\gamma, s}(\rnp)\\
&&+\mu_{\gamma, s}(\rnp)(\varphi(x_0)(\varphi(x_0)-2))+o(1)\\
&=&o(1) \quad \hbox{as $k\to +\infty$. }
\end{eqnarray*}
Since $\gamma<\frac{n^2}{4}$, the Hardy inequality \eqref{ineq:hardy:rk} then yields $\int_{\rnp}|\nabla(1-\varphi)u_k|^2\, dx\to 0$ as $k\to +\infty$. Therefore $\int_{\rnp\setminus B_{2\delta}(x_0)}|\nabla u_k|^2\, dx\to 0$ as $k\to +\infty$. \\
Now we show that
\begin{equation}\label{lim:l2}
\lim_{k\to +\infty}\int_{\rnp}\frac{u_k^2}{|x|^2}\, dx=0.
\end{equation}
Indeed, it follows from the analysis above and Hardy's inequality, that for any $\delta>0$, we have  $\int_{\rnp}\frac{((1-\varphi)u_k)^2}{|x|^2}\, dx\to 0$ when $k\to +\infty$, and therefore $\int_{\rnp\setminus B_{2\delta}(x_0)}\frac{u_k^2}{|x|^2}\, dx\to 0$ as $k\to +\infty$. If we take $\delta<|x_0|/2$, since $u_k\to 0$ in $L^2_{loc}(\rn)$, then $\int_{\rnp\cap B_{2\delta}(x_0)}\frac{u_k^2}{|x|^2}\, dx\to 0$ as $k\to +\infty$, which  proves \eqref{lim:l2}.\\

\medskip\noindent We show that 
\begin{equation}\label{eq:Kn2}
\mu_{\gamma, 0}(\rnp)=\frac{1}{K(n,2)^2}.
\end{equation}
Indeed, it follows from \eqref{mes:l}, \eqref{mes:n} and \eqref{lim:l2} that $\mu_{\gamma, s}(\rnp)\geq\frac{1}{K(n,2)^2}$, which implies \eqref{eq:Kn2} since \eqref{bnd:mu} holds.

\medskip\noindent Finally, we show that if $\gamma>0$, then $n=3$. Now, consider the family $u_\eps$ as in the proof of Proposition \ref{prop:s0:gamma:neg}. Well known computations by Aubin \cite{aubin} yield that 
$$J_{\gamma,s}^{\rnp}(u_\eps)=K(n,2)^{-2}-\gamma|x_0|^{-2}c\theta_\eps+o(\theta_\eps)\hbox{\, as $\eps\to 0$,}$$
where $c>0$, $\theta_\eps=\eps^2$ if $n\geq 5$ and $\theta_\eps=\eps^2\ln\eps^{-1}$ if $n=4$. It follows that if $\gamma>0$ and $n\geq 4$, then $\mu_{\gamma, s}(\rnp)<K(n,2)^{-1}$, hence contradicting \eqref{eq:Kn2}. Therefore $n=3$ and we are done. 
 
\medskip\noindent As a conclusion, we have proved that if there are no extremals, then $s=0$ and \eqref{eq:Kn2} holds. If in addition $\gamma>0$, then $n=3$. This proves Proposition \ref{prop:ext:rnp}. 

\section{Appendix 2: Symmetry of the extremals for $\mu_{\gamma, s}(\rnp)$}
The symmetry of the nonnegative solutions to the Euler-Lagrange equation for $\mu_{\gamma,s}(\rnp)$ is proved in Chern-Lin \cite{CL5} for $\gamma<(n-2)^2/4$. The proof of the symmetry carried out by Ghoussoub-Robert \cite{gr2} in the case $\gamma=0$ extends immediately to the case $0\leq \gamma<n^2/4$. For the convenience of the reader, we give here a general and complete proof inspired by Chern-Lin \cite{CL5}. We consider nontrivial solutions $u\in D^{1,2}(\rnp)$ to the problem
\begin{equation}\label{app:eq:1}
\left\{\begin{array}{ll}
-\Delta u-\frac{\gamma}{|x|^2}u=\frac{u^{\crits-1}}{|x|^s}&\hbox{ weakly in }D^{1,2}(\rnp)\\
u\geq 0 &\hbox{ in }\rnp\\
u= 0 &\hbox{ on }\partial\rnp\\
\end{array}\right.
\end{equation}
and prove the following. Here, $\gamma<n^2/4$, $s\in [0,2)$ and $\crits:=\frac{2(n-s)}{n-2}$. 

\begin{theorem}\label{th:sym} We let $u\in D^{1,2}(\rnp)$ be a solution to \eqref{app:eq:1}. Then  $u\circ\sigma=u$ for 
all isometry of $\rn$ such that $\sigma(\rnp)=\rnp$. In particular, there 
exists $v\in C^\infty((0,+\infty)\times \rr)$ such 
that for all $x_1>0$ and all $x'\in\rr^{n-1}$, we have that $u(x_1,x')=v(x_1,|x'|)$.
\end{theorem}
\noindent{\bf Remark:} Unlike the case of the extremals for the full space $\rn$, there is no symmetry-breaking phenomenon in the case of the half-space $\rnp$. However, the price to pay is that the best constant when restricted to the functions with best possible symmetry is unknown, contrary to the case of $\rn$. We refer to the historical reference Catrina-Wang \cite{CW} and to Dolbeault-Esteban-Loss-Tarantello \cite{DELT} for disussions and developments on the symmetry-breaking phenomenon.

\medskip\noindent We adapt the moving-plane method of Chern-Lin \cite{CL5} that was made in the case $\gamma<\frac{(n-2)^2}{4}$. Given any $\theta\in [0,\frac{\pi}{2}]$, we define the hyperplane and the half space:
$$P_\theta:=\{x\in\rn/\, x_1\cos\theta=x_2\sin\theta\},$$
$$P_{\theta}^-:=\{x\in\rn/\, x_1\cos\theta < x_2\sin\theta\}.$$
We define $s_\theta:\rn\to\rn$ as the orthogonal symmetry with respect to $P_\theta$. As one checks, we have that
\begin{equation}\label{app:exp:s:theta}
s_\theta(x)=\left(\begin{matrix} -x_1\cos(2\theta)+x_2\sin(2\theta)\\x_1\sin(2\theta)+x_2\cos(2\theta)\\ x_i\, (i\geq 2)
\end{matrix}\right)\hbox{ and }s_\theta(x)-x=2\left(x_2\sin\theta-x_1\cos\theta\right)\left(\begin{matrix} \cos\theta\\-\sin\theta\\ 0
\end{matrix}\right)
\end{equation}
Note that it follows from Theorem \ref{th:hopf} that there exists $K_1>0$ such that 
\begin{equation}\label{app:eq:6}
u(x)\sim_{x\to 0}K_1\frac{x_1}{|x|^{\am}}.
\end{equation}
The proof of Theorem \ref{th:sym} relies on two main Lemmae:

\begin{lem}\label{lem:diff} For all $j=1,n...,n$, we have that
\begin{equation}\label{app:eq:5}
\lim_{x\to 0}\left(|x|^{\am}\partial_j u(x)-K_1\left(\delta_{j,1}-\am\frac{x_1x_j}{|x|^2}\right)\right)=0,
\end{equation}
and 
\begin{equation}\label{app:eq:5:bis}
|x|^{\am+1}\Vert d^2u_x\Vert\leq C\hbox{ for all }x\in\rnp, \; |x|<1.
\end{equation}
\end{lem}
\noindent{\it Proof of Lemma \ref{lem:diff}:} We proceed by contradiction and assume that there exists $(x_k)_k\in\rnp$ such that $x_k\to 0$ and
\begin{equation}\label{app:eq:8} 
 \left(|x_k|^{\am}\partial_j u(x_k)-K_1\left(\delta_{j,1}-\am\frac{x_{k,1} x_{k,j}}{|x_k|^2}\right)\right)\not\to 0
\end{equation}
as $k\to +\infty$. We define $u_k(x):=|x_k|^{\am-1}u(|x_k|x)$ for all $x\in\rnp$. It follows from \eqref{app:eq:6} that 
\begin{equation}\label{app:eq:7}
\lim_{k\to +\infty}u_k(x)= K_1\frac{x_1}{|x|^{\am}}\hbox{ for all }x\in\overline{\rnp}\setminus\{0\}.
\end{equation}
Moreover, this convergence holds in $C^0_{loc}(\overline{\rnp}\setminus\{0\})$. The equation \eqref{app:eq:1} rewrites
$$-\Delta u_k-\frac{\gamma}{|x|^2}u_k=|x_k|^{(\crits-2)\left(\frac{n}{2}-\am\right)}\frac{u_k^{\crits-1}}{|x|^s}\hbox{ in }\rnp$$ 
for all $k$, and $u_k$ vanishes on $\partial\rnp$. It then follows from elliptic theory that the convergence in \eqref{app:eq:7} holds in $C^2_{loc}(\overline{\rnp}\setminus\{0\})$. Therefore, 
$$\lim_{k\to +\infty}\partial_j u_k\left(\frac{x_k}{|x_k|}\right)=\partial_j(K_1 x_1|x|^{-\am})(X_\infty)$$
where $X_\infty:=\lim_{k\to +\infty}\frac{x_k}{|x_k|}$. Coming back to $u_k$ contradicts \eqref{app:eq:8}. This proves \eqref{app:eq:5}. The proof of \eqref{app:eq:5:bis} is similar. This ends the proof of Lemma \ref{lem:diff}.

\medskip\noindent The second Lemma is a general analysis of the difference $u(s_\theta(x))-u(x)$.
\begin{lem}\label{lem:asymp:gene} We let $(\theta_i)_i\in\rr$ and $(x_i)\in\rnp$ be such that $x_i\in \rnp\cap P_{\theta_i}^-$ for all $i\in\nn$. We assume that $\theta_i\to \theta_\infty$ and $x_i\to x_\infty$ as $i\to +\infty$, and that
\begin{equation}\label{hyp:o}
s_{\theta_i}(x_i)- x_i=o(|x_i|)\hbox{ as }i\to +\infty.
\end{equation}
Then,
\begin{itemize}
\item If $x_\infty\neq 0$ then
\begin{equation}\label{case:1}
\lim_{i\to +\infty}\frac{u(s_{\theta_i}(x_i))-u(x_i)}{2(x_{i,2}\sin\theta_i-x_{i,1}\cos\theta_i)}=\cos(\theta_\infty)\partial_1u(x_\infty)-\sin(\theta_\infty)\partial_2u(x_\infty).
\end{equation}
\item If $x_\infty= 0$, then
\begin{equation}\label{case:2}
\lim_{i\to +\infty}\frac{u(s_{\theta_i}(x_i))-u(x_i)}{2(x_{i,2}\sin\theta_i-x_{i,1}\cos\theta_i)|x_i|^{-\am}}=K_1\cos(\theta_\infty).
\end{equation}
\end{itemize}
\end{lem}
\noindent{\it Proof of Lemma \ref{lem:asymp:gene}:} Taylor's formula yields
\begin{equation}
\left| u(s_{\theta_i}(x_i))-u(x_i)-du_{x_i}(s_{\theta_i}(x_i)-x_i)\right|\leq \Vert s_{\theta_i}(x_i)-x_i\Vert^2\sup_{t\in [0,1]}\Vert d^2u_{x_i+t(s_{\theta_i}(x_i)-x_i)}\Vert
\end{equation}
for all $i$. It follows from \eqref{app:eq:5:bis} that
\begin{eqnarray*}
\sup_{t\in [0,1]}\Vert d^2u_{x_i+t(s_{\theta_i}(x_i)-x_i)}\Vert &\leq& C \sup_{t\in [0,1]}\Vert x_i+t(s_{\theta_i}(x_i)-x_i)\Vert^{-(1+\am)}\\
&=&C \left\Vert \frac{x_i+s_{\theta_i}(x_i)}{2}\right\Vert^{-(1+\am)},
\end{eqnarray*}
and therefore
\begin{equation}
u(s_{\theta_i}(x_i))=u(x_i)+du_{x_i}(s_{\theta_i}(x_i)-x_i)+O\left(\frac{\Vert s_{\theta_i}(x_i)-x_i\Vert^2}{\Vert x_i+s_{\theta_i}(x_i)\Vert^{1+\am}}\right),
\end{equation}
and then, we get with \eqref{hyp:o} that
\begin{equation}
u(s_{\theta_i}(x_i))=u(x_i)+du_{x_i}(s_{\theta_i}(x_i)-x_i)+o\left(\frac{\Vert s_{\theta_i}(x_i)-x_i\Vert }{|x_i|^{\am}}\right),
\end{equation}
as $i\to +\infty$. With the expression \eqref{app:exp:s:theta}, we get that 
\begin{equation}\label{eq:norm}
\Vert s_{\theta_i}(x_i)-x_i\Vert=2(x_{i,2}\sin\theta_i-x_{i,1}\cos\theta_i)>0
\end{equation}
and that
\begin{equation}
\frac{u(s_{\theta_i}(x_i))-u(x_i)}{2(x_{i,2}\sin\theta_i-x_{i,1}\cos\theta_i)}=\partial_1 u(x_i)\cos\theta_i-\partial_2 u(x_i)\sin\theta_i+o\left(|x_i|^{-\am}\right)
\end{equation}
as $i\to +\infty$. If $x_\infty\neq 0$, then we get \eqref{case:1} and we are done. We assume that $x_\infty=0$. It then follows from Lemma \ref{lem:diff} that
\begin{eqnarray*}
\frac{|x_i|^{\am}(u(s_{\theta_i}(x_i))-u(x_i))}{\Vert s_{\theta_i}(x_i)-x_i\Vert}&=&K_1\left(1-\am\left(\frac{x_{i,1}}{|x_i|}\right)^2\right)\cos\theta_i\\
&&+K_1\am \frac{x_{i,1}x_{i,2}}{|x_i|^2}\sin\theta_i+o(1)\\
&=& K_1\left[\cos\theta_i +\am\frac{x_{i,1}}{|x_i|}\left(\frac{x_{i,2}\sin\theta_i-x_{i,1}\cos\theta_i}{|x_i|}\right)\right]\\
&&+o(1)\hbox{\, as $i\to +\infty$.}
\end{eqnarray*}
Using \eqref{hyp:o} and \eqref{eq:norm}, we get that
\begin{equation}
\lim_{i\to +\infty}\frac{|x_i|^{\am}(u(s_{\theta_i}(x_i))-u(x_i))}{\Vert s_{\theta_i}(x_i)-x_i\Vert}=K_1\cos\theta_\infty.
\end{equation}
This ends the proof of Lemma \ref{lem:asymp:gene}.\hfill$\Box$

\medskip\noindent We are now in position to initiate the moving plane method.
\begin{proposition}\label{prop:1} There exists $\theta_0>0$ such that
\begin{equation}\label{res:claim:1}
\hbox{for all }\theta\in (0,\theta_0),\hbox{ then }u(s_\theta(x))>u(x)\hbox{ for all }x\in P_\theta^-\cap\rnp
\end{equation}
\end{proposition}
\smallskip\noindent{\it Proof of Proposition \ref{prop:1}:} We argue by contradiction and we assume that there exists $(\theta_i)_i\in (0,+\infty)$, there exists $x_i\in P_{\theta_i}^-\cap\rnp$ such that
\begin{equation}\label{app:eq:2}
\lim_{i\to +\infty}\theta_i=0\hbox{ and }u(s_{\theta_i}(x_i))\leq u(x_i)\hbox{ for all }i.
\end{equation}

\medskip\noindent We first claim that without loss of generality, we can assume that $(x_i)_i$ is bounded in $\rn$. For that we define the Kelvin transform $\tilde{u}(x):=|x|^{2-n}u(x/|x|^2)$ for all $x\in\rnp$. As one checks, $\tilde{u}\in D^{1,2}(\rnp)$ satisfies \eqref{app:eq:1} and \eqref{app:eq:2} rewrites $\tilde{u}(s_{\theta_i}(\tilde{x}_i))\leq \tilde{u}(\tilde{x}_i)$ for all $i$, where $\tilde{x}_i:=x_i/|x_i|^2\in P_{\theta_i}^-\cap\rnp$. Therefore, up to changing $u$ into $\tilde{u}$, we can assume that $(x_i)_i$ is bounded. This proves the claim.

\medskip\noindent Now define $\lim_{i\to +\infty }x_i=x_\infty$. We claim that 
\begin{equation}\label{claim:lim}
x_{i,1}=o(x_{i,2})\hbox{ as }i\to +\infty\hbox{ and }x_\infty\in\partial\rnp.
\end{equation}
Indeed, since $x_i\in P_{\theta_i}^-\cap\rnp$, we have that $x_{i,1}>0$ and $x_{i,1}\cos\theta_i < x_{i,2}\sin\theta_i$ for all $i$. Letting $i\to \infty$ yields $x_{i,1}=o(x_{i,2})$ as $i\to +\infty$, and therefore $x_\infty\in\partial\rnp$. 

\medskip\noindent We now show that 
\begin{equation}\label{claim:lim:2}
s_{\theta_i}(x_i)- x_i=o(|x_i|)\hbox{ as }i\to +\infty.
\end{equation}
Indeed, it suffices to note that the expression \eqref{eq:norm} and \eqref{claim:lim} yield $s_{\theta_i}(x_i)- x_i=o(|x_{i,2}|)=o(|x_i|)$ as $i\to +\infty$. 

\medskip\noindent We now conclude the proof of Proposition \ref{prop:1}. If $x_\infty=0$, it follows from \eqref{case:2} that $u(s_{\theta_i}(x_i)) -u(x_i)>0$ for $i\to +\infty$, contradicting \eqref{app:eq:2}. If $x_\infty\neq 0$, it follows from \eqref{case:1} and \eqref{app:eq:2} that $\partial_1 u(x_\infty)\leq 0$: this contradicts Hopf's strong maximum principle since $x_\infty\in\partial\rnp$. This ends the proof of Proposition \ref{prop:1}.
\hfill$\Box$\\

\medskip\noindent Define now 
\begin{equation*}
\theta_0:=\sup\left\{0< \theta\leq \frac{\pi}{2}/\, u(s_t(x))>u(x)\hbox{ for all }x\in P_t^-\cap\rnp\hbox{ and all }0<t<\theta\right\}
\end{equation*}
It follows from Proposition \ref{prop:1} that $\theta_0>0$ exists. Our objective is to prove that $\theta_0=\frac{\pi}{2}$. We argue by contradiction and assume that
\begin{equation}\label{app:hyp:theta}
0<\theta_0<\frac{\pi}{2}.
\end{equation}
For any $\theta\geq 0$, we define 
$$v_\theta(x):=u(s_\theta(x))-u(x)$$	
for all $x\in P_\theta^-\cap\rnp$. Since $s_\theta$ is an isometry for all $\theta\geq 0$, we have that
\begin{equation}\label{app:eq:v}
-\Delta v_\theta-\frac{\gamma}{|x|^2}v_\theta=c_\theta(x) v_\theta
\end{equation}
where $c_\theta(x)=|x|^{-s}\frac{u(s_\theta(x))^{\crits-1}-u(x)^{\crits-1}}{u(s_\theta(x))-u(x)}$ if $u(s_\theta(x))\neq u(x)$, and $c_\theta(x)=|x|^{-s}(\crits-1)u(x)^{\crits-2}$ otherwise. In particular, $c_\theta>0$. It follows from the definition of $\theta_0$ that $v_{\theta_0}\geq 0$. It then follows from \eqref{app:eq:v} and Hopf's maximum principle that either $v_{\theta_0}> 0$ in $P_{\theta_0}^-\cap\rnp$ or $v_{\theta_0}\equiv 0$ in $P_{\theta_0}^-\cap\rnp$. In the latest case, taking points on $\partial\rnp$, we would get that $u(x)=0$ on $P_{2\theta_0}\cap\rnp$: this is impossible since $\theta_0<\frac{\pi}{2}$ and $u>0$. Therefore
\begin{equation}\label{app:pos:v}
v_{\theta_0}> 0\hbox{ in }P_{\theta_0}^-\cap\rnp. 
\end{equation}
It follows from the definition of $\theta_0$ that there exists $(\theta_i)_i\in (\theta_0,+\infty)$ such that
\begin{equation}\label{app:ppty:theta:i}
\lim_{i\to +\infty}\theta_i=\theta_0\hbox{ and }\forall i\hbox{ there exists }x_i\in P_{\theta_i}^-\cap\rnp\hbox{ such that }v_{\theta_i}(x_i)\leq 0.
\end{equation}
Arguing as in Step 1 of the proof of Proposition \ref{prop:1}, we can assume with no loss of generality that $(x_i)_i$ is bounded, and, up to a subsequence, that there exists $x_\infty\in \rn$ such that $\lim_{i\to +\infty}x_i=x_\infty$.

\medskip\noindent We claim that 
\begin{equation}\label{app:x:in:P}
x_\infty\in P_{\theta_0}\cap\overline{\rnp}.
\end{equation}
Indeed, it follows from \eqref{app:ppty:theta:i} that $x_\infty\in \overline{P_{\theta_0}^-\cap\rnp}$ and $v_{\theta_0}(x_\infty)\leq 0$. It then follows from \eqref{app:pos:v} that $x_\infty\in \partial{P_{\theta_0}^-\cap\rnp}=(P_{\theta_0}\cap\overline{\rnp})\cup (\partial\rnp\cap \overline{P_{\theta_0}^-}$ and $v_{\theta_0}(x_\infty)=0$. We argue by contradiction and assume that \eqref{app:x:in:P} does not hold. Therefore, $x_\infty\in\partial\rnp$ and $u(s_{\theta_0}(x_\infty))=v_{\theta_0}(x_\infty)=0$, and then $s_{\theta_0}(x_\infty)\in\partial\rnp$. We then get with \eqref{app:exp:s:theta} and \eqref{app:hyp:theta} that $s_{\theta_0}(x_\infty)=x_\infty$ and then $x_\infty\in P_{\theta_0}$, which contradicts our initial hypothesis. This proves \eqref{app:x:in:P} and therefore the claim.

\medskip\noindent We claim that 
\begin{equation}\label{app:eq:13}
s_{\theta_i}(x_i)- x_i=o(|x_i|)\hbox{ as }i\to +\infty.
\end{equation}
It follows from \eqref{app:x:in:P} that $s_{\theta_0}(x_\infty)=x_\infty$, and therefore \eqref{app:eq:13} holds if $x_\infty\neq 0$. We now assume that $x_\infty=0$. Dividing \eqref{app:eq:13} by $|x_i|$ and passing to the limit $i\to +\infty$, one gets that \eqref{app:eq:13} is equivalent to proving that $s_{\theta_0}(X_\infty)=X_\infty$ where $X_\infty:=\lim_{i\to +\infty}\frac{x_i}{|x_i|}$. Since $x_i\in P_{\theta_i}^-$, we have that $x_{i,2}\sin\theta_i>x_{i,1}\cos\theta_i$ for all $i$. Dividing by $|x_i|$ and passing to the limit $i\to +\infty$ yields 
\begin{equation}\label{app:eq:11}
X_{\infty,2}\sin\theta_0\geq X_{\infty,1}\cos\theta_0.
\end{equation}
Since $u(s_{\theta_i}(x_i))\leq u(x_i)$, the asymptotic \eqref{app:eq:7} yields
$$K_1\frac{(s_{\theta_i}(x_i))_1}{|s_{\theta_i}(x_i)|^{\am}}\leq (1+o(1))K_1\frac{x_{i,1}}{|x_i|^{\am}}$$
as $i\to +\infty$. Dividing by $|x_i|$ and passing to the limit, we get that 
\begin{equation}\label{app:eq:12}
(s_{\theta_0}(X_\infty)-X_\infty)_1\leq 0.
\end{equation}
Plugging \eqref{app:eq:11} and \eqref{app:eq:12} into \eqref{app:exp:s:theta} yields $s_{\theta_0}(X_\infty)=X_\infty$. As already mentioned, this proves the claim.

\medskip\noindent Here goes the final argument. We apply Lemma \ref{lem:asymp:gene}. If $x_\infty=0$, \eqref{case:2}, \eqref{eq:norm} and \eqref{app:ppty:theta:i} yield $K_1\cos(\theta_0)\leq 0$: a contradiction since $K_1>0$ and $0<\theta_0<\pi/2$. If $x_\infty\neq 0$, \eqref{case:1} and \eqref{app:ppty:theta:i} yield
\begin{equation}\label{ineq:app:1}
\partial_1u(x_\infty)\cos(\theta_0)-\partial_2u(x_\infty)\sin(\theta_0)\leq 0.
\end{equation}
If $x_\infty\in\partial\rnp$, then $\partial_2u(x_\infty)=0$ and $\partial_1 u(x_\infty)>0$ (Hopf's Lemma), contradicting \eqref{ineq:app:1}. So $x_\infty\in P_{\theta_0}\setminus\partial\rnp$. It then follows from \eqref{app:eq:v}, \eqref{app:pos:v}, \eqref{app:x:in:P} and Hopf's Lemma that $\partial_{\vec{N}}v_{\theta_0}(x_\infty)<0$ with $\vec{N}=(\cos\theta_0,-\sin\theta_0, 0)$. However, one can easily see that  $\partial_{\vec{N}}v_{\theta_0}(x_\infty)=-2(\partial_1u(x_\infty)\cos(\theta_0)-\partial_2u(x_\infty)\sin(\theta_0))$, which again contradicts \eqref{ineq:app:1}.

\medskip\noindent In all cases, we get a contradiction, and therefore \eqref{app:hyp:theta} is not valid, which means that $\theta_0=\frac{\pi}{2}$. It follows that 
$$u(x_1,-x_2,...)\geq u(x_1, x_2,...)\hbox{ for all }x\in\rnp,\, x_2>0.$$
Since the equation satisfied by $u$ is invariant under the action of isometries fixing $\partial\rnp$, we get the reverse inequality and therefore $u(x_1,x_2,...)=u(x_1, -x_2,...)$ for all $x\in\rnp$. So $u$ is invariant under the action of the symmetry wrt $\{x_2=0\}$. This argument works for any hyperplane orthogonal to $\partial\rnp$: then $u$ is invariant under the action of the symmetries fixing $\partial\rnp$. This completes the proof of Theorem \ref{th:sym}.

\end{document}